%% file: embeddings.tex
\documentclass[11pt,a4paper]{amsart}
\usepackage[T1]{fontenc}
\usepackage{amssymb}
\usepackage{enumitem}
\usepackage{ifthen}
\usepackage{tabularx}
\usepackage{calc}
\usepackage[backref]{hyperref}

\makeatletter
 \theoremstyle{plain}    
 \newtheorem{thm}{Theorem}[section]
 \numberwithin{equation}{section} 
 \numberwithin{figure}{section} 
 \usepackage{verbatim}
 \theoremstyle{plain}
 \theoremstyle{definition}
  \newtheorem{example}[thm]{Example}
 \theoremstyle{remark}    
 \newtheorem{notation}[thm]{Notation} 
 \theoremstyle{remark}    
 \newtheorem{claim}[thm]{Claim}
 \theoremstyle{definition}
 \newtheorem{defn}[thm]{Definition}
 \theoremstyle{plain}    
 \newtheorem{prop}[thm]{Proposition} 
 \theoremstyle{plain}    
 \newtheorem{lem}[thm]{Lemma} 
 \newtheorem{sublem}[thm]{Sublemma}
 \theoremstyle{plain}    
 \newtheorem{cor}[thm]{Corollary} 
 
 \theoremstyle{definition}
 \newtheorem*{defn*}{Definition}
 
 \theoremstyle{plain}    
 \newtheorem*{prop*}{Proposition} 
 \newtheorem{question}{Question}
 \theoremstyle{remark}
 \newtheorem{remark}[thm]{Remark}
 
 \newlength{\label@width}
 \newlength{\label@sep}
 \newlength{\left@margin}

 \newenvironment{enumeq}{
 \settowidth{\label@width}{(\thesection.\arabic{equation})}
 \settowidth{\label@sep}{\espc}
 \setlength{\left@margin}{\label@width+\label@sep}
 \begin{enumerate}[label=(\thesection.\arabic*), ref=\thesection.\arabic*, labelwidth=\label@width, 
leftmargin=\left@margin, labelsep=\label@sep]\setcounter{enumi}{\value{equation}}}
 {\setcounter{equation}{\value{enumi}}\end{enumerate}}

\include{macros}

\newcommand{\basel}{\iprod}
\newcommand{\jid}{\tu{\textrm{(JID)}}}
\renewcommand{\mid}{\tu{\textrm{(MID)}}}
\newcommand{\pom}{\tu{\textrm{POM}}}
\renewcommand{\random}{\operatorname{Ran}}
\newcommand{\IFF}{\qquad\mathrm{iff}\qquad}

\newcommand{\upcl}[1]{{{\uparrow}#1}}

\showhyphens{monoid}

\title{Partial order embeddings with convex range}
\author{James Hirschorn}
\date{January 2, 2007.}
\address{Graduate School of Science and Technology, Kobe University, Japan}
\email{\href{mailto:j_hirschorn@yahoo.com}{j\_hirschorn@yahoo.com}}
\thanks{The author acknowledges the generous support of the Japanese Society for the
  Promotion of Science (JSPS Fellowship for Foreign Researchers, ID\# P04301).}
\urladdr{\href{http://www.logic.univie.ac.at/~hirschor/}{http://www.logic.univie.ac.at/\textasciitilde  hirschor/}}

\begin{document}

\begin{abstract}
A careful study is made of embeddings of posets which have a convex range. 
We observe that such embeddings 
share nice properties
with the homomorphisms of more restrictive categories; for example, we show that
every order embedding between two lattices with convex range 
is a continuous lattice homomorphism. 
A number of posets are considered; for example, we prove that
every product order embedding $\sigma:\irr\to\irr$ with convex range is of the
form 
\begin{equation}
  \sigma(x)(n)=\bigl((x\circ g_\sigma)+y_\sigma\bigr)(n)
  \espc\text{if $n\in K_\sigma$},\label{eq:68}
\end{equation}
and $\sigma(x)(n)=y_\sigma(n)$ otherwise, for all $x\in\irr$, where
$K_\sigma\subseteq\N$, $g_\sigma:K_\sigma\to\N$ is a bijection and
$y_\sigma\in\irr$. The most complex poset examined here is the quotient of the 
lattice of Baire measurable functions,
with codomain of the form $\N^I$ for some index set $I$, modulo equality on a
comeager subset of the domain, with its `natural' ordering. 
\end{abstract}

\maketitle

\tableofcontents

\section{Overview}
\label{sec:overview}

The abstract objects of study here are embeddings between arbitrary posets having
preregular (e.g.~convex, see definition~\ref{d-2}) range. 
The main observation is that
they are continuous with respect to a natural poset topology (namely the Scott
topology). This observation is
then applied to investigate embeddings with convex range between 
some specific classes of posets. 

It turns out that the determination of these embeddings for various examples
of posets generalizes known results about homomorphisms in more restrictive
categories. For example, we show in~\Section\ref{sec:new-subsection}
that every order embedding $\sigma:\power(X)\to\power(Y)$ with convex range is
of the form $\sigma(a)=h[a]\cup b$ for some injection $h:X\to Y$ 
and some $b\subseteq Y$; this can be compared with the known fact that every 
continuous Boolean algebra monomorphism $\sigma:\power(X)\to\power(Y)$ 
is of the form $\sigma(a)=h\inv[a]$ for some partial surjection $h:Y\pfn X$. 

By definition, an order embedding is an order preserving map (i.e.~partial order
homomorphism) that is an isomorphism onto its range. Thus for order theoretic
structures, e.g.~lattices and Boolean algebras, where the structure is completely
determined by the ordering, an order embedding is in fact an isomorphism for the
given structure onto its range. Hence characterizing these embeddings involves two
aspects: determining the isomorphisms of the given structure; and determining the
range of these embeddings. For example, it is a special case of the above mentioned
known fact,
that Boolean algebra isomorphisms between $\power(X)$ and $\power(Y)$ are of the
form $\sigma(a)=h[a]$ for some bijection $h:X\to Y$. 
Thus the `new result' in the above characterization concerns the range,
namely that every order embedding $\sigma:\power(X)\to\power(Y)$
with convex range has range equal to the interval
$[\sigma(\emptyset),\sigma(X)]=\{a\subseteq Y:\sigma(\emptyset)\subseteq
a\subseteq\sigma(X)\}$. (Of course, the computation of the range in this example
follows trivially from the definition of convexity, 
but in slightly more complex examples this computation can 
become difficult). However, this is an incomplete view of the situation
where there is additional structure that is not purely order theoretic. 
Indeed many our examples are also monoids. 

Write $\N$ for the set $\{0,1,\dots\}$ of nonnegative integers. The usual (linear)
partial ordering of $\N$ is $0<1<\cdots$.
The example most important to us is the class of partial order embeddings 
between the irrationals $\irr$, i.e.~the set of all functions from $\N$ into~$\N$ 
(see e.g.~\cite{irrationals}), with the product order.
Indeed this paper is the second part of a series, where the third part~\cite{irrationals} is
entitled ``Characterizing the quasi ordering of the irrationals by eventual
dominance'' (and it is the sequel to~\cite{pinning}), where order embeddings of the
irrationals with convex range play a crucial role in characterizing the eventual
dominance ordering (cf.~\Section\ref{sec:irrfin-lefnt}) 
in terms of the product order. 

More generally, we consider embeddings between the product orders
$\irri I$ and $\irri J$ for arbitrary index sets $I$ and $J$. Thus for all
$x,y\in\irri I$,
\begin{equation}
  \label{eq:1}
  x\le y\Iff x(i)\le y(i)\spc\text{for all $i\in I$}.
\end{equation}
In this case we ascertained (cf.~\Section\ref{sec:powers-n}) 
that every order embedding from $\irri I$ into $\irri J$
with convex range is in fact a monoid embedding for coordinatewise addition 
plus  a constant in~$\irri J$. 

Further generalization was desired 
to function spaces with the irrationals as the codomain.
Let $X$ and $Y$ be topological spaces. We let $\baire(X,Y)$ denote the family of
all Baire measurable functions from $X$ into $Y$. 
In \Section\ref{sec:baire-functions} we consider quotients of $\baire(X,\irr)$ (where
$\irr$ has the product topology) over the equivalence relation of ``almost always''
equality $\eqaa$:
\begin{equation}
  \label{eq:67}
  f\eqaa g\If f(z)=g(z)\espc\text{for almost all $z\in X$,}
\end{equation}
or in other words for comeagerly many $z\in X$; where $\baireaa{X,\irr}$ is given
the order $\leaa$ induced by $\le$:
 \begin{equation}
   \label{eq:2}
   [f]\leaa[g]\Iff f(z)\le g(z)\espc\text{for almost all $z\in X$}.
 \end{equation}
The embeddings of $\baireaa{X,\irr}$ into $\baireaa{Y,\irr}$ with convex range
are described precisely, where, as far as we know, 
even the isomorphism structure was not previously known. 
This involves characterizing embeddings with convex range 
of category algebras of topological spaces, 
or equivalently the regular open algebras of these spaces, 
in~\Section\ref{sec:category-algebras}. We note that there are
some subtleties in generalizing to arbitrary index sets, 
i.e.~to $\baireaa{X,\irri I}$, that are addressed there. 

Our interest in quotients of Baire functions stems from set theoretic
forcing. For example, for an index set $I$ in the `reference model'---called the
\emph{ground model} in the terminology of set theoretic forcing, 
it is well known that every member of $\irri I$ in the extension of this model
obtained by Cohen forcing, is determined by a member of $\baireaa{\reals,\irri I}$ 
in the ground model. 
Whilst in the other direction, we shall use set theoretic forcing to prove
that $(\baireaa{X,\irri I},\leaa)$ is a complete semilattice satisfying many of the
same properties as $\irri I$ (theorem~\ref{l-2}). 
This can be generalized extensively; for example, for $S$ in the
ground model, every member of $\power(S)$  in the Cohen model corresponds to a
member of $\baireaa{\reals,\power(S)}$ 
(cf.~\Section\ref{sec:bairex-powersdiveqaa}), and by replacing ``almost
always'' with ``almost everywhere'' in the measure theoretic sense, we obtain a
correspondence with random forcing (cf.~\Section\ref{sec:l0mu-irri-idiveqae}).

In \Section\ref{sec:quasi-order-mono} we study a standard association of a quasi order to
every monoid. We are especially interested in those 
monoids---we call them lattice monoids---where
the associated quasi order is in fact a lattice. Then various infinite distributive
laws are examined for these lattice monoids. The main purpose of this section for
the present paper, is that it
allows us in~\Section\ref{sec:algebraic-domains} to use algebraic methods 
to extend embeddings from a suitably dense subset of some lattice monoid to the
entire lattice (cf.~corollary~\ref{o-34}). This in turn is applied
in~\Section\ref{sec:cont-funct-into} to use our description 
in~\Section\ref{sec:baire-functions} of the embeddings of $\baireaa{X,\irri I}$ into
$\baireaa{Y,\irri J}$ with convex range in order to 
obtain a 
precise description of the class of embeddings with
convex range of $\cts(X,\irri I)$ into $\cts(Y,\irri J)$. 
That is, the family of continuous functions ordered pointwise.  
We also see examples (e.g.~lemma~\ref{l-27}, theorem~\ref{l-3}) of how
algebraic properties of subsets of lattices may entail order theoretic regularity
properties. 

The general theory is in~\Section\ref{sec:homom-with-regul}, including the Scott
continuity of embeddings with preregular range, 
and numerous consequences of this result.
In~\Section\ref{sec:bases-lattices}, we study various notions of
denseness in a lattice, 
and this is applied in~\Section\ref{sec:algebraic-domains} to extend homomorphisms
from suitably dense subsets of a lattice while preserving various desirable properties. 

Further directions for the Boolean algebra $(\pnfin,\subseteqfnt)$ and its close
relative the lattice $(\irrfin,\lefnt)$ are suggested
in~\Section\ref{sec:pnfin-subseteqfnt} and~\Section\ref{sec:irrfin-lefnt}. 
\subsection{Terminology}
\label{sec:terminology}

A \emph{homomorphism} refers to an arrow (i.e.~morphism) 
of the intended category,
while an \emph{isomorphism} is an invertible homomorphism. 
For a concrete category
$(C,U)$, we say that an \emph{embedding} is a homomorphism $h:c\to\nobreak d$ 
for which
there exists some invertible homomorphism $g$ with domain $c$ such that 
$Uh(x)=Ug(x)$ for all $x\in Uc$. Thus in a category where the
homomorphisms are functions, an embedding is a homomorphism $f:X\to Y$ that
is an isomorphism onto its range, i.e.~$\ran(f)$ is an object of the category and
$f':X\to\ran(f)$ given by $f'(x)=f(x)$ is an isomorphism. \emph{Monomorphisms}
(i.e.~monics) refer to homomorphisms that are left cancellative under composition, 
i.e.~$f$ is a monomorphism iff $f\circ g=f\circ h$ implies $g=h$ for all $g,h$ with
codomain equal to the domain of $f$. 
And \emph{epimorphisms} (i.e.~epis) refer to homomorphisms that are right cancellative.
In all of the concrete categories considered below, 
monomorphisms are simply the injective homomorphisms, 
and in all of these categories with the exception of the monoids, 
epimorphisms are precisely the surjective homomorphisms.

A \emph{quasi order} (also often called a \emph{preorder}) is a pair $(O,\le)$ where
$\le$ is a reflexive and transitive relation on $O$. For a quasi order $(O,\le)$, we
write $<$ for the relation defined by $p<q$ if $p\le q$ and $q\nleq p$. 
Note this disagrees with another usage where $p<q$ iff $p\le q$ and $p\ne q$. 
A \emph{poset} (\emph{partial order}) is a quasi order where the relation is also
antisymmetric (i.e.~$p\le q$ and $q\le p$ imply $p=q$). In a poset $(P,\le)$, 
$p<q$ iff $p\le q$ and $p\ne q$. A poset with a minimum element is called a
\emph{pointed poset}. 
Note that if $(O,\le)$ is a quasi order (poset) 
then for every subset $A\subseteq O$, $(A,\le)$ is also a quasi order (poset). We
call it a \emph{quasi suborder} (\emph{subposet}). 
In the category of quasi orders, the homomorphisms are order preserving maps,
i.e.~for two quasi orders $(O,\le)$ and $(Q,\altle)$, $\sigma:O\to Q$ is \emph{order
  preserving} if $p\le q$ implies $\sigma(p)\altle\sigma(q)$ for all $p,q\in
O$. Thus isomorphisms are bijections that are both order preserving and order
reflecting, where $\sigma:O\to Q$ is \emph{order reflecting} if
$\sigma(p)\altle\sigma(q)$ implies $p\le q$ for all $p,q\in O$. 
A subset $A\subseteq O$ is called \emph{bounded above}, or just \emph{bounded}, if
it has an \emph{upper bound}, 
i.e.~some $p\in O$ such that $a\le p$ for all $a\in A$. And we say that $A$ is
\emph{bounded below} if it has a \emph{lower bound}, i.e.~some $q\in O$ such that
$q\le a$ for all $a\in A$. 

The class of posets is viewed as a full subcategory of the quasi orders. 
Notice that for $(P,\le)$ a poset and $(Q,\altle)$ a quasi order, $\sigma:P\to Q$ is
an embedding iff it is both order preserving and reflecting. We write $\spr A$ for
the supremum, i.e.~minimum upper bound, of a subset $A\subseteq P$ 
(which may or may not exist),  and we write $\ifm A$ for the infimum, i.e.~maximum
lower bound. We write $a \lor b$ and $a \land b$ for $\spr\{a,b\}$ and
$\ifm\{a,b\}$, respectively. A function $\sigma:P\to Q$ between two posets is called
\emph{join-preserving} if $\sigma(p\lor q)=\sigma(p)\lor\sigma(q)$ whenever $p\lor
q$ exists, for all $p,q \in P$ (i.e.~if $p\lor q$ exists then so does
$\sigma(p)\lor\sigma(q)$ satisfying the equation). The dual notion is
\emph{meet-preserving}. 

We take a \emph{join semilattice}, which we also just call a \emph{semilattice}, 
to be a poset $(L,\le)$
such that $\lor$ is a binary operation on $L$, i.e.~$p\lor q$ exists for all $p,q\in
L$. And a \emph{meet semilattice} is defined dually. A \emph{lattice} is a poset
that is both a semilattice and a meet semilattice. 
A \emph{subsemilattice} of a semilattice $(L,\le)$ is a semilattice $(A,\le)$, where
$A\subseteq L$, such that $(A,\le)\models\ulc c=a\lor b\urc$ iff $(L,\le)\models\ulc
c=a\lor b\urc$, i.e.~the supremum computed in the order $(A,\le)$ agrees with the
supremum taken in~$L$. \emph{Meet subsemilattices} and \emph{sublattices} are
defined analogously. Note that a subset $A\subseteq L$ may be a lattice as a
subposet of $(L,\le)$, without being a sublattice. 
The homomorphisms of the category of semilattices are the join-preserving functions,
while the homomorphisms of the category of meet semilattices are the meet-preserving
functions. Note that these are both subcategories of the poset category, because
$p\lor q$ and $p\land q$ do not exist when $p\le q$, $q\le p$ and $q\ne p$, and
because join/meet semilattice homomorphisms are order preserving, as $p\le q$ iff
$p\lor q=q$ iff $p\land q=p$. The category of lattices is the intersection of the
categories of join and meet semilattices, i.e.~the homomorphisms are the functions
that are both join and meet-preserving.
It is easy to find a counterexample showing that a quasi
order homomorphism between two lattices need not be a lattice homomorphism
(i.e.~lattices are not a full subcategory of the quasi orders).
On the other hand, since the lattice operations are obviously determined by
the ordering of the lattice, a quasi order isomorphism between two lattices
(equivalently, a poset isomorphism) is in fact a lattice isomorphism. 
Take note that lattice embeddings are the same thing as lattice monomorphisms (and
similarly for join/meet semilattice embeddings),
i.e.~they are the injective lattice homomorphisms (this is not true of the category
of posets). Moreover:

\begin{prop}
\label{p-34}
A lattice homomorphism $\sigma$ is an embedding 
iff it is \emph{strictly order preserving},
i.e.~$p<q$ implies $\sigma(p)<\sigma(q)$.
\end{prop}

By a \emph{complete semilattice} we mean a join semilattice $(L,\le)$ 
such that $\spr A$ exists whenever $A\subseteq L$ is bounded.
Note that a complete semilattice is pointed, 
with minimum element $\spr\emptyset$, which we denote by $0$. 
Notice also that a complete semilattice is in fact
a lattice, and moreover $\ifm A$ exists whenever $A\ne\emptyset$
(see e.g.~\cite{MR1902334}). A \emph{complete lattice} is a lattice such that 
$\spr A$ and $\ifm A$ exist for every subset $A\subseteq L$. By adding a top
element to any complete semilattice one obtains a complete lattice. A lattice is
called \emph{bounded} if it has both a maximum and minimum element; we denote the
maximum element by $1$.

We take a \emph{Boolean algebra} to be a bounded lattice $(B,\le)$
such that every $p\in B$ has a
\emph{complement}, which we write as $-p$, satisfying $p\land-p=0$ and $p\lor-p=1$. 
Homomorphisms in the category of Boolean algebras are lattice homomorphisms that
preserve complements (and thus also preserve $0$ and $1$). 
Since the Boolean algebra operations are completely determined by the order, 
every poset isomorphism between two Boolean
algebras is in fact a Boolean algebra isomorphism. In the category of Boolean
algebras embeddings and monomorphisms coincide.

Recall that a \emph{monoid} $(M,\cdot)$ is a semigroup that has an identity. We do
not need to specify the identity because it is uniquely determined; we denote it by
$e$. A monoid is
\emph{cancellative} if $a\cdot b=a\cdot c$ implies $b=c$ and $b\cdot a=c\cdot a$
implies $b=c$, for all $a,b,c\in M$. The inverse of $a\in M$, when it exists, is
denoted $a\inv$; more generally, we write $a\cdot b\inv$ for the element $c\in M$
such that $a=c\cdot b$, if it exists. It is uniquely determined so long as $M$ is cancellative.
We write $(M,+)$ when dealing with a
commutative monoid; use $-$ to denote the inverse; and we use $0$ to denote the identity.
A monoid homomorphism is map preserving both the monoid
operation and the identity.
In the category of monoids embeddings and monomorphisms both
coincide with injective homomorphisms. When we say that a subset $S$ of a commutative 
monoid is \emph{closed under subtraction} we of course mean that $a-b\in S$ whenever
$a,b\in S$ and $a-b$ exists. Since every cancellative commutative monoid embeds into
an Abelian group (lemma~\ref{l-5}), a submonoid of a cancellative commutative monoid 
that is closed under subtraction can be viewed as a `subgroup'.

\begin{prop}
\label{p-51}
Let $(G,+)$ be an Abelian group, and $M\subseteq G$ be a submonoid.
Then a subsemigroup $S\subseteq M$ is closed under subtraction iff $\<S\>\cap M=S$,
where $\<S\>$ denotes the subgroup generated by $S$. 
\end{prop}
\begin{proof}
Since $S$ is a subsemigroup of an Abelian group, $\<S\>=S-S$. 
It immediately follows that $\<S\>\cap M=S$ iff $S$ is closed under subtraction.
\end{proof}

Recall that a binary relation $R$ on a set $X$ is a \emph{congruence} on some
$n$-ary relation $S$ on $X$, 
if $S(x_0,\dots,x_{n-1})$ iff $S(x'_0,\dots,x'_{n-1})$,
whenever $x_k\rel x'_k$ for all $k=0,\dots,n-\nobreak1$; 
and it is a \emph{congruence} on
some $n$-ary function $f$ on $X$, if $f(x_0,\dots,x_{k-1})\rel
f(x'_0,\dots,x'_{k-1})$ whenever $x_k\rel x'_k$ for all $k=0,\dots,n-\nobreak1$.
A congruence $\sim$ on some quasi order $(O,\le)$, i.e.~a congruence on $\le$, 
that is moreover an equivalence relation on $O$ is
an \emph{orderable} partition of $(O,\le)$
in the terminology of~\cite{MR620665}, and it is called so because
$\le$ determines a well defined ordering 
of the quotient $O\div{\sim}$ via $[p]\le[q]$ if $p\le q$. A \emph{lattice
  congruence} on some lattice $(L,\le)$ means a relation on $L$ that is 
a congruence for both of the binary lattice operations. A lattice congruence $\sim$
that is also an equivalence relation induces a lattice structure on $L\div{\sim}$,
with $[a]\land[b]=[a\land b]$ and $[a]\lor[b]=[a\lor b]$. And a \emph{Boolean
  algebra} congruence is a lattice congruence that is also a congruence for the
unary complement operator. Of course a Boolean algebra congruence that is an equivalence
relation on $(B,\le)$ determines a quotient Boolean algebra $B\div{\sim}$ with
$-[a]=[-a]$. 

Similarly, a \emph{monoid congruence} for a monoid $(M,\cdot)$ is a congruence on
the monoid operation $\cdot$, and an equivalence relation $\sim$ 
that is a congruence on $\cdot$ determines a quotient monoid $(M\div{\sim},\cdot)$.
\section{The quasi ordering of a monoid}
\label{sec:quasi-order-mono}

All monoids have a naturally associated quasi ordering. This is well known,
e.g.~\cite{MR0285644} to give one example. 
In this section we examine some basic properties of this associated quasi ordering,
and introduce the notion of a lattice monoid. Then we focus on distributivity, where
we prove that a large class of monoids satisfy
certain infinitary distributive laws; for example, in corollary~\ref{o-10} we prove
that in every cancellative commutative monoid, addition distributes over arbitrary suprema.
As far as the examples of~\Section\ref{sec:partial-orders} are concerned, 
the results of this section are only
applied to the lattice monoid of~\Section\ref{sec:cont-funct-into}. 
However, lattice monoids will play a bigger role in the sequel to this paper. 

\begin{defn}
A monoid $(M,\cdot)$ has an \emph{associated quasi order} $\len{(M,\cdot)}$
on $M$ defined by
\begin{equation}
\label{eq:60}
x\len{(M,\cdot)} y\If x\cdot a=y\textrm{ for some }a\in M\textrm{.}
\end{equation}
\end{defn}

\begin{prop}
\label{p-45}
If $(M,\cdot)$ is a monoid the the relation defined in equation~\eqref{eq:60} is a
quasi ordering.
\end{prop}
\begin{proof}
It is reflexive because $(M,\cdot)$ has an identity, and it is transitive because
$(M,\cdot)$ is a semigroup.
\end{proof}

\begin{prop}
\label{p-16}
For any monoid, the identity is a minimum element of the associated quasi order.
\end{prop}
\begin{proof}
For all $a\in M$, $e\cdot a=a$ implies $e\len{(M,\cdot)} a$.
\end{proof}

The quasi order associated with a cancellative monoid 
has no maximal elements unless the monoid
is a group, in which case the associated quasi order is the complete quasi order.

\begin{lem}
\label{p-20}
Let $(M,\cdot)$ be a cancellative monoid.
If $(M,\len{(M,\cdot)})$ has a maximal element then $(M,\cdot)$ is a group.
\end{lem}
\begin{proof}
Let $a\in M$ be a maximal element. Take any $b\in M$. 
Then $a\len{(M,\cdot)} a\cdot
b$ implies $a\cdot b\len{(M,\cdot)} a$ and thus $a\cdot b\cdot c=a$ for some $c\in
M$. 
Hence $b\cdot c=e$ by cancellativity. 
The existence of right inverses for all elements entails that $(M,\cdot)$ is a group. 
\end{proof}

\begin{prop}
\label{p-4}
Let $(M,\cdot)$ be a monoid. Then $b\len{(M,\cdot)}c$ implies
$a\cdot b\len{(M,\cdot)}a\cdot c$, for all $a,b,c\in M$.
\end{prop}
\begin{proof}
Assuming $b\len{(M,\cdot)} c$, there exists $d$ such that $b\cdot d=c$.
Thus $(a\cdot b)\cdot d=a\cdot(b\cdot d)=a\cdot c$.
\end{proof}

We are mostly interested in commutative monoids so that we at least have
monotonicity. Indeed, the category of 
\emph{positively quasi ordered commutative monoids} 
is widely studied in Ordered Algebra, see e.g.~\cite{MR1190444} where they are named
\pom's, and the
quasi order $\len{(M,+)}$ associated with a commutative monoid is the minimal
quasi ordering such that $(M,+,\len{(M,+)})$ is a \pom.

\begin{prop}
\label{p-18}
Let $(M,+)$ be a commutative monoid.
Then $a\len{(M,+)} b$ implies $a+c\len{(M,+)} b+c$, for all $a,b,c\in M$.
\end{prop}
\begin{proof}
Write  $a+d=b$. Then $(a+c)+d=a+d+c=b+c$.
\end{proof}

\begin{defn}
A \emph{poset monoid} is a monoid $(M,\cdot)$ whose associated quasi order is in
fact a partial order.
\end{defn}

\begin{prop}
\label{p-36}
Every poset monoid is pointed with $e$ the minimum element.
\end{prop}
\begin{proof}
Proposition~\ref{p-16}.
\end{proof}

\begin{lem}
\label{p-15}
Poset monoids do not have invertible elements besides the identity.
\end{lem}
\begin{proof}
Suppose that $a\in M$ is invertible. 
Then $a\cdot a\inv=e$ implies $a\len{(M,\cdot)} e$, and $e\len{(M,\cdot)}a$ by
proposition~\ref{p-16}. Thus $a=e$ by antisymmetry.
\end{proof}

\begin{cor}
\label{p-21}
Cancellative poset monoids do not have maximal elements, with the exception of the
singleton monoid $\{e\}$. 
\end{cor}
\begin{proof}
Lemmas~\ref{p-20} and~\ref{p-15}.
\end{proof}

\begin{defn}
A monoid $(M,\cdot)$ is called a \emph{semilattice monoid} if the associated quasi order
is actually a semilattice, and it is called a \emph{lattice monoid} if the
associated quasi order is moreover a lattice. 
\end{defn}

\begin{remark}
\label{r-5}
Note that in a semilattice monoid $0=e$ by proposition~\ref{p-36}.
Thus our usage of $0$ as the additive identity in a commutative monoid is consistent
with its usage as the minimum element of a pointed semilattice. 
\end{remark}

We recall the following basic fact and provide a proof.

\begin{lem}
\label{l-5}
A commutative monoid is embeddable \tu(as a monoid\tu) in some
Abelian group iff it is cancellative.
\end{lem}
\begin{proof}
Let $(M,+)$ be commutative and cancellative. Let $M^-$ denote the subgroup of
invertible elements of $M$. Note that commutativity entails that
$N=(M\setminus M^-)\cup\{0\}$ is a submonoid of $(M,+)$.
Now consider the product monoid $M\times N$ modulo the relation defined by
\begin{equation}
  \label{eq:46}
  (a,b)\sim(\bar a,\bar b)\If \exists r,\bar r\in N\spc a+r=\bar a+\bar r
  \and b+r=\bar b+\bar r.
\end{equation}
Note that $\sim$ is an equivalence relation:
$0\in N$ implies reflexivity, symmetry is obvious, 
and we have transitivity because supposing $(a,b)\sim(\bar a,\bar b)$ 
and $(\bar a,\bar b)\sim(\bar{\bar a},\bar{\bar b})$,
$a+r=\bar a+\bar r$ and $\bar a+s=\bar{\bar a}+\bar s$ which with commutativity imply 
$a+r+s=\bar a+\bar r+s=\bar{\bar a}+\bar r+\bar s$, 
and similarly $b+r+s=\bar{\bar b}+\bar r+\bar s$. 
And using commutativity similarly, 
$\sim$ is a congruence for addition (cf.~\Section\ref{sec:terminology}). 
Therefore ${(M\times N)}\div{\sim}$ is a quotient monoid. Observe that
\begin{equation}
  \label{eq:47}
  [b,b]=[0,0]\espc\text{for all $b\in N$,}
\end{equation}
and of course $[0,0]$, i.e.~$[(0,0)]$, is the identity. 

\begin{claim}
\label{c-5}
$({(M\times N)}\div{\sim},+)$ is an Abelian group.
\end{claim}
\begin{proof}
It is commutative because it is the quotient of a commutative monoid. 
Take $[a,b]$ in the quotient. 
We prove that it has an inverse $-[a,b]$. First suppose that $a\in M^-$. 
Then $-[a,b]=[-a+b,0]$ because $[a,b]+[-a+b,0]=[b,b]=[0,0]$
by~\eqref{eq:47}. Otherwise $-[a,b]=[b,a]$ because
$[a,b]+[b,a]=[a+b,a+b]$. 
\end{proof}

\begin{claim}
\label{c-9}
$a\mapsto [a,0]$ is a monomorphism between $M$ and ${(M\times N)}\div{\sim}$.
\end{claim}
\begin{proof}
$a+b$ is mapped to $[a+b,0]=[a,0]+[b,0]$ and hence it a monoid homomorphism.
Now suppose $[a,0]=[b,0]$. Then there are $r$ and $s$ in $N$ such that $a+r=b+s$,
and $0+r=0+s$, i.e.~$r=s$. 
And by the cancellative property $a+r=b+r$ implies $a=b$ as needed.
\end{proof}

Claims~\ref{c-5} and~\ref{c-9} establish one direction of the lemma, and it is
clear that a noncancellative monoid cannot be embedded in a group. 
\end{proof}

This is used to show that cancellative commutative monoids satisfy monotonicity for
subtraction. 

\begin{prop}
\label{p-52}
Let $(M,+)$ be a cancellative commutative monoid. Whenever $a-c$ exists,
$b\len{(M,+)} c$ implies $a-b$ exists and $a-b\gen{(M,+)}a-c$, for all $a,b,c\in M$. 
\end{prop}
\begin{proof}
By lemma~\ref{l-5}, we may assume that $M$ is a submonoid of some Abelian group
$(G,+)$. Now $b\len{(M,+)}c$ implies $b+d=c$ for some $d\in M$, and thus
$(a-c)+d=(a-(b+d))+d=a-b$ as wanted, by the basic properties of a group. 
\end{proof}

\begin{prop}
\label{p-53}
Let $(M,+)$ be a cancellative commutative monoid. Whenever $a-c$ exists,
$a\len{(M,+)} b$ implies $b-c$ exists and $a-c\len{(M,+)}b-c$,
for all $a,b,c\in M$. 
\end{prop}
\begin{proof}
Write $a+d=b$. Then since by lemma~\ref{l-5} we can assume we are working inside an
Abelian group, $(a-c)+d=(a+d)-c=b-c$ by basic properties of an Abelian group.
\end{proof}
\subsection{Distributive laws}
\label{sec:infin-distr-laws}

Recall that a lattice $(L,\lor,\land)$ is \emph{distributive} 
if $a\land (b\lor c)=(a\land b)\lor(a\land c)$ for all $a,b,c\in L$, and that:

\begin{prop}
\label{p-19}
$(L,\lor,\land)$ is distributive iff $a\lor(b\land c)=(a\lor b)\land(a\lor c)$ for
all $a,b,c\in L$.
\end{prop}

\begin{defn}
We say that a monoid $(M,\cdot)$ is 
\emph{left $\cdot$-distributive over $\lor$}, or
simply \emph{left $(\cdot,\lor)$-distributive}, if
\begin{equation}
  \label{eq:43}
  a\cdot(b\lor c)=(a\cdot b)\lor(a\cdot c)\espc\text{whenever $b\lor c$ exists}
\end{equation}
for all $a,b,c\in M$, and it is \emph{right $\cdot$-distributive over $\lor$} if
\begin{equation}
  \label{eq:44}
  (a\lor b)\cdot c=(a\cdot c)\lor(b\cdot c)\espc\text{whenever $a\lor b$ exists}
\end{equation}
for all $a,b,c\in M$, while it is \emph{$\cdot$-distributive over $\lor$} (or just
\emph{$(\cdot,\lor)$-distributive}) if it is
both left and right $\cdot$-distributive over $\lor$.
The definition of $\cdot$-distributivity over $\land$ is exactly analogous.

Let us say that it is \emph{$(\cdot,\lorand)$-distributive} if it is both
$(\cdot,\lor)$-distributive and $(\cdot,\land)$-distributive.
We will say that a lattice monoid is \emph{distributive} if it is both
$(\cdot,\lorand)$-distributive \emph{and} distributive as a lattice.
\end{defn}

Note that while we are mostly interested in $\cdot$-distributivity over $\lor$ in
semilattice monoids, and $\cdot$-distributivity over $\land$ is meet semilattice
monoids, the definition does make sense for arbitrary monoids.

Of course, the notions of left $(+,\lor)$-distributivity, right
$(+,\lor)$-distributivity and $(+,\lor)$-distributivity all coincide for
any commutative monoid $(M,+)$, 
and similarly for $+$-dis\-tribu\-tiv\-ity over $\land$. 

\begin{lem}
\label{p-33}
Let $(M,+)$ be a commutative $(+,\land)$-distributive monoid. 
Then $a\land\nobreak b=0$ implies $a\lor b=a+b$, for all $a,b\in M$. 
\end{lem}
\begin{proof}
By commutativity, $a,b\le a+b=b+a$. Suppose $a\land b=0$ and $d\ge a,b$. 
Let $c$ satisfy $a+c=d$. 
Then by right $(+,\land)$-distributivity, 
$c=(a\land b)+c=(a+c)\land(b+c)=d\land (b+c)\ge b$ implies $d\ge a+b$ by
proposition~\ref{p-4}. 
\end{proof}

\begin{lem}
\label{l-20}
Let $(M,+)$ be a commutative cancellative $+$-distributive monoid over~$\lor$ and
$\land$. Then $a\land c=0$ and $b\land c=0$ imply $(a+b)\land c=0$.
\end{lem}
\begin{proof}
Suppose  that $a\land c=b\land c=0$,  and take $d\le (a+b),c$.
By lemma~\ref{p-33} and $+$-distributivity over $\lor$, 
$a+b+d=a+(b\lor d)=(a+b)\lor(a+d)$, and similarly $a+b+d=(a+b)\lor(b+d)$. Thus
\begin{equation}
  \label{eq:127}
  a+b+d=(a+b)\lor(a+d)\lor(b+d).
\end{equation}
But we also have 
$a+b+2d=(a+d)+(b+d)=(a\lor d)+(b\lor d)=(a+b)\lor(a+\nobreak d)\lor(d+b)\lor (d+d)$ 
by distributivity over $\lor$. Therefore, as $2d\le a+b+d$ since $d\le a+b$, we obtain
$a+b+2d=a+b+d$, and thus $d=0$ by cancellativity.
\end{proof}

\subsubsection{Infinite distributive laws}
\label{sec:infin-distr-laws-2}

\begin{notation}
We extend binary operations to sets in the usual way, 
e.g.~$a\land B=\{a\land b:b\in B\}$, $A\cdot b=\{a\cdot b:a\in A\}$, etc... .
\end{notation}

A lattice $L$ is \emph{join-infinite distributive} \jid\ if
\begin{equation}
  \label{eq:95}
  a\land\left(\spr B\right)=\spr(a\land B)
\end{equation}
whenever the supremum on the left hand side exists. The dual condition is called
\emph{meet-infinite distributive} \mid. Note that distributivity does not imply
either of these properties, even for complete lattices (see e.g.~\cite{MR1902334}). 
Also note that every Boolean
algebra is both join-infinite and meet-infinite distributive (see
e.g.~\cite[Ch.~1, \Section1]{MR991565}). We will need to recall that $a$ can be
replaced with a supremum:

\begin{prop}
\label{p-27}
If $L$ is \jid, then $\bigl(\spr A\bigr)\land\bigl(\spr B\bigr)=\spr(A\land B)$ 
whenever $\spr A$ and $\spr B$ both exist. Similarly for the \mid. 
\end{prop}

\begin{defn}
We shall call a monoid $(M,\cdot)$ 
\emph{infinitely left $\cdot$-distributive over $\lor$} 
or \emph{infinitely left $(\cdot,\lor)$-distributive} if
\begin{alignat}{2}
  \label{eq:50}
  a\cdot\spr B&=\spr(a\cdot B)
  \espc&&\text{for all $a\in M$ and $B\subseteq M$,}
  \intertext{whenever $\spr B$ exists, 
    and call it \emph{infinitely right $(\cdot,\lor)$-distributive} if}
  \label{eq:55}
  \left(\spr A\right)\cdot b&=\spr(A\cdot b)
  \espc&&\text{for all $A\subseteq M$ and $b\in M$,}  
\end{alignat}
whenever $\spr A$ exists,  and we call the monoid \emph{infinitely
  $(\cdot,\lor)$-distributive} if it is both infinitely left and right 
  $(\cdot,\lor)$-distributive. The definition of \emph{infinite
  $(\cdot,\land)$-distributivity} is exactly analogous.
\end{defn}

\begin{prop}
\label{p-32}
If $(M,\cdot)$ is infinitely $(\cdot,\lor)$-distributive, 
then $\bigl(\spr A\bigr)\cdot\bigl(\spr B\bigr)=\spr(A\cdot B)$
whenever $\spr A$ and $\spr B$ both exist;  similarly for infinite
$(\cdot,\land)$-distributivity. 
\end{prop}

\begin{prop}
\label{p-17}
Let $(M,+)$ be a commutative monoid. 
Then infinite left and right $(+,\lor)$-distributivity both coincide with
$(+,\lor)$-distributivity. And infinite left and right 
$(+,\land)$-distributivity both coincide with $(+,\land)$-distributivity. 
\end{prop}

\begin{example}
\label{x-12}
The monoid $(\N,+)$ is a lattice monoid, because $\len{(\N,+)}$ agrees with the
usual ordering of $\N$, which is a lattice with supremum $\max$ and infimum $\min$.
It is a commutative cancellative and distributive lattice monoid, that embeds into
the Abelian group of integers $(\intgr,+)$. It is also a complete semilattice 
satisfying \jid, \mid\ and infinite $(+,\lorand)$-distributivity. 
\end{example}

\begin{example}
\label{x-9}
Let $I$ be an index set, and consider the product monoid $(\irri I,+)$, i.e.~with
coordinatewise addition. Since the quasi order associated with a product monoid is
the corresponding product of the associated quasi orders, $(\irri I,+)$ is a
commutative lattice monoid with $x\len{(\irri I,+)} y$ iff $x(i)\le y(i)$ for all
$i\in I$. It is a complete semilattice that is cancellative, 
\jid, \mid\ and infinitely $+$-distributive over~$\lorand$, 
because all of these properties are preserved under products, i.e.~the
suprema and infima of a product lattice are taken coordinatewise (see also
theorem~\ref{l-7}). It embeds into the generalized Baer--Specker group $(\intgr^{I},+)$. 
\end{example}

We shall start simplifying the notation by writing $\le$ instead of
$\len{(M,\cdot)}$ and $\len{(M,+)}$. 

\begin{lem}
\label{l-11}
Every cancellative monoid $(M,\cdot)$ is infinitely left $(\cdot,\lor)$-distributive. 
\end{lem}
\begin{proof}
Suppose that $a\in M$ and $\spr B$ exists.
$a\cdot b\le a\cdot\spr B$ for all $b\in B$ by proposition~\ref{p-4}, 
and thus $a\cdot\spr B$ is an upper bound of $a\cdot B$.
On the other hand, suppose that $r$ is an upper bound of $a\cdot B$. 
Then for all $b\in B$, $r\ge a\cdot b$ implies that there exists $c_b\in M$ with
\begin{equation}
  \label{eq:58}
  a\cdot b\cdot c_b=r.
\end{equation}
Then by the cancellative property, 
there exists $d\in M$ with $b\cdot c_b=d$ for all $b\in B$. 
Now since $d\ge\spr B$, $r=a\cdot d\ge a\cdot\spr B$ by proposition~\ref{p-4},
completing the proof that $\spr(a\cdot B)=a\cdot\spr B$.
\end{proof}

\begin{cor}
\label{o-10}
Every commutative cancellative monoid $(M,+)$ is infinitely $+$-distributive over
$\lor$. 
\end{cor}
\begin{proof}
Lemma~\ref{l-11} and proposition~\ref{p-17}.
\end{proof}

\begin{thm}
\label{l-7}
Every commutative cancellative semilattice monoid $(M,+)$ is both infinitely
$(+,\lor)$-distributive and $(+,\land)$-distributive.
\end{thm}
\begin{proof}
Infinite $(+,\lor)$-distributivity is by corollary~\ref{o-10}. 
Suppose then that  $\ifm B$ exists. Then $a+b\ge a+\ifm B$ for all $b\in B$, and
thus $a+\ifm B$ is a lower bound for $a+B$. To prove that it is in fact
the greatest lower bound, 
since $M$ is a semilattice it will suffice take a lower bound $r\ge a+\ifm B$ and
show that $r= a+\ifm B$. Now since $r\ge a$, $r-a\in M$. Therefore, 
for all $b\in B$, $r\le a+b$ implies that $r-a\le (a+b)-a=b$ by
proposition~\ref{p-53}.
Thus $r-a\le\ifm B$, and hence $r=a+(r-a)\le a+\ifm B$ as required.
\end{proof}
\section{Order embeddings with preregular range}
\label{sec:homom-with-regul}

The main result of this section (theorem~\ref{l-12}) is that order embeddings
between posets, with some regularity property imposed on their ranges, 
are continuous with respect to the Scott topology on these posets. 
The other major theme of this section appears in a series of results 
(lemma~\ref{l-19}--corollary~\ref{o-34}) on extending a given
continuous homomorphism from a suitably `dense' subset of some lattice to the entire
lattice.
\subsection{Preregularity}
\label{sec:preregularity}

We introduce the notion of preregularity and observe (lem\-ma~\ref{l-6}) that convex
subsets of lattices are preregular. We also prove (lemma~\ref{l-27}) that every dense
`subgroup' of a lattice monoid is regular. 

\begin{notation}
For a poset $(P,\le)$ and a subposet $A\subseteq P$ we write $\spr^A B$ for the
supremum of $B\subseteq A$ taken in the poset $(A,\le)$, and similarly for infima.
\end{notation}

\begin{prop}
\label{p-13}
Suppose $(P,\le)$ is a poset and $A\subseteq P$ is a subposet. 
Then for all $B\subseteq A$\tu: 
if $\spr^A B$ and $\spr^P B$ both exist then $\spr^P B\le \spr^A B$\tu;
and if $\ifm^A B$ and $\ifm^P B$ both exist then $\ifm^P B\ge\ifm^A B$.
\end{prop}

\begin{defn}
Let $(P,\le)$ be a poset. We call a subset $A\subseteq P$ \emph{order closed} if
for every nonempty $\emptyset\ne B\subseteq A$: 
if $p=\spr^P B$ exists in $(P,\le)$ then $p\in A$; and if
$q=\ifm^P B$ exists in $(P,\le)$ then $q\in A$. 
\end{defn}

We will want to consider a weaker property than order closed, and also to separate
upwards and downwards closedness. 

\begin{defn}
\label{d-1}
Let us call a subset $A\subseteq P$ \emph{upwards boundedly order closed} 
if for every
nonempty $B\subseteq A$ that has an upper bound in $A$, if $p=\spr^P B$ exists
then $p\in A$; and we call $A$ \emph{downwards boundedly order closed} if 
for every nonempty $B\subseteq A$ with a lower bound in $A$,
if $p=\ifm^P B$ exists then $p\in A$. And we call $A\subseteq P$ \emph{boundedly
  order closed} if it is both upwards and downwards boundedly order closed. 

We can separate the order closed property into \emph{upwards order closed} and
\emph{downwards order closed} analogously.
\end{defn}

\begin{remark}
\label{r-8}
Note that an upwards order closed subset should be distinguished from an
\emph{upwards closed} subset $A$ of $P$, which of course refers to a set $A$
satisfying $\upcl A=A$ where $\upcl A=\{p\in P:p\ge a$ for some $a\in A\}$.
To avoid possible confusion we follow~\cite{MR1365749} and call them 
\emph{upper sets}. 
Analogous remarks are made for \emph{downwards closed} subsets, which we call
\emph{lower sets}. 
\end{remark}

\begin{defn}
\label{d-10}
For a subposet $A$ of $(P,\le)$, 
the \emph{upwards \tu(downwards\tu) order closure} of $A$ is the
smallest upwards (downwards) order closed $Q\subseteq P$ containing $A$. 
\end{defn}

\noindent Note that the upwards (downwards) order closure always exists because 
the family of upwards (downwards) order closed sets is closed under intersections.
Indeed, they have explicit descriptions.

\begin{prop}
\label{p-43}
The upwards order closure of $A\subseteq P$ is given by
\begin{equation}
  \label{eq:130}
  \left\{\spr B:B\subseteq A\textup{ and the supremum
        exists}\right\},
\end{equation}
and thus the downwards order closure has the dual description. 
\end{prop}

\begin{prop}
\label{p-28}
For any $B\subseteq A\subseteq P$, if $\spr^P B$ exists and is an element of $A$,
then $\spr^A B=\spr^P B$. Similarly for infima.
\end{prop}

\noindent Hence, for example:

\begin{prop}
\label{p-10}
If $A\subseteq P$ is order closed then for all $\emptyset\ne B\subseteq A$, 
$\spr^A B=\spr^P B$ if
the supremum exists in $P$. Similarly for infima. 
\end{prop}

On the other hand, the supremum may exist in the subposet $A$ but not in $P$. This
leads to the following definition (see e.g.~\cite[Ch.~1, \Section1]{MR991565}).

\begin{defn}
A subset $A\subseteq P$ is \emph{upwards regular} if for all $B\subseteq A$,
if $a=\spr^A B$ exists in $(A,\le)$ then $\spr^P B=a$ in $(P,\le)$, and it is
\emph{downwards regular} if for all $B\subseteq A$,
if $b=\ifm^A B$ exists in $A$ then $\ifm^P B=b$. A subset is \emph{regular} if it is
both upwards and downwards regular.
\end{defn}

However, we want to avoid the singularity of empty suprema and infima, 
and thus we instead define the following concept.

\begin{defn}
\label{d-2}
A subset $A\subseteq P$ is \emph{upwards preregular} 
if for all nonempty $\emptyset\ne B\subseteq A$: 
if $a=\spr^A B$ exists then $\spr^P B=a$, and it is \emph{downwards preregular} 
if for all $\emptyset\ne B\subseteq A$:
if $b=\ifm^A B$ exists then $\ifm^P B=b$.
It is \emph{preregular} if it is both upwards and downwards preregular. 
\end{defn}

We point out that preregularity is a transitive property.

\begin{prop}
\label{p-54}
Let $(P,\le)$ be a poset.
If $A\subseteq P$ is preregular, and $B\subseteq A$ is preregular as a subset of the
poset $(A,\le)$, then $B$ is preregular as a subset of the poset $P$. 
\end{prop}

\begin{lem}
\label{l-22}
Every boundedly order closed subset of a complete semilattice is preregular. 
\end{lem}
\begin{proof}
Given a complete semilattice $L$,
let $A\subseteq L$ be boundedly order closed, and $\emptyset\ne B\subseteq A$.
If $\spr^A B$ exists, then $B$ has an upper bound in $A$,
and thus $a=\spr^L B$ exists, and $a\in A$. Now $\spr^A B=a$ by
proposition~\ref{p-28}. And if $\ifm^A B$ exists then $B$ has a lower bound in $A$,
and $\ifm^L B$ exists because $B\ne\emptyset$; hence, 
we can conclude that $\ifm^L B=\ifm^A B$. 
\end{proof}

\begin{defn}
\label{def:convex}
Let $(O,\le)$ be a quasi order. We call $A\subseteq O$ \emph{convex} if $p,q\in A$
implies $[p,q]\subseteq A$ where $[p,q]$ is the interval $\{r\in O:p\le r\le q\}$. 
\end{defn}

\begin{lem}
\label{l-6}
A convex subset of any lattice is preregular.
\end{lem}
\begin{proof}
Suppose $(L,\le)$ is a lattice, $A\subseteq L$ is convex, 
and $\emptyset\ne B\subseteq A$, say with $b\in B$. 
Assume $a=\spr^A B$ exists. Then if $p\in L$ is an upper bound of $B$, 
so is $p\land a$. However, 
$[b,a]\subseteq A$ and thus $b\le p\land a\le a$ implies $p\land a\in A$, 
which entails that $p\ge a$ as needed. Same for the infimum.
\end{proof}

\begin{lem}
\label{l-27}
Let $(M,+)$ be a cancellative commutative lattice monoid. 
Then every dense \tu(cf.~definition~\tu{\ref{d-11})} submonoid of $M$ that is closed
under subtraction is regular.
\end{lem}
\begin{proof}
Let $D\subseteq M$ be a dense submonoid closed under subtraction. 
Suppose $A\subseteq D$ and $a=\spr^D A$. 
Assume towards a contradiction that $a\ne \spr^M A$. Then there is an upper bound
$b\in M$ of $A$ such that $a\nleq b$, and thus $a\land b$ is an upper bound of $A$,
and $a\land b<a$, which by proposition~\ref{p-36} means there is a $c\in M^+$ 
with $(a\land b)+c=a$. By density, there exists $d\le c$ in $D^+$. 
Now $a-c=a\land b$, and thus $a-d$ exists and $a-d\ge a\land b$ by
proposition~\ref{p-52}. Thus $a-d$ is an upper bound of $A$.  
But $a\in D$, and hence $a-d\in D$ because it is closed under subtraction.
Since $a-d<a$ this is contrary to $a$ being the supremum in $D$, thereby completing
the proof of upwards regularity.

Suppose $A'\subseteq D$ and $a'=\ifm^D A'$. Assuming towards a contradiction that
$a'$ is not the infimum in $M$, there exists a lower bound $b'>a'$ of $A'$. Writing
$a'+c'=b'$, there exists $d'\le c'\in D^+$ by density. But since $D$ is a
subsemigroup, $a'+d'\in D$, and it is a lower bound of $A'$ since $a'+d'\le b'$.
This is contrary to $a'$ being the infimum since $a'<a'+d'$, thereby completing the
proof of downwards regularity.
\end{proof}

\begin{cor}
\label{o-35}
Let $(M,+)$ be a cancellative commutative lattice monoid. Then every dense submonoid
that is closed under subtraction forms a sublattice. 
\end{cor}
\subsection{Order reflecting homomorphisms}
\label{sec:order-refl-homom}

\begin{notation}
\label{sec:infin-distr-laws-1}
For a relation $(S,\le)$, we write $S^0$ for the set of all $\le$-minimal elements
of $S$, and $S^+=S\setminus S^0$. 
\end{notation}

\noindent Note that a lattice has at most one minimal element in which case it is the
$0$ of the lattice. 

\begin{defn}
\label{d-7}
Let $(O,\le)$ be a quasi order.
An element of $O$ is an \emph{atom} if it is a nonminimal element 
that cannot be split, where we say
that $p$ can be \emph{split} if there exists $q,r\le p$ in $O^+$ such that $q$ is
\emph{incompatible} with $r$, written $q\incompat r$
(i.e.~there is no common extension of $q$ and $r$ in $O^+$).
Write $\atom(O,\le)$ for the collection of atoms of $O$. 
$O$ is \emph{atomless} if $O^+$ has no atoms, 
whereas $O$ is \emph{atomic} if $\atom(O,\le)$ is dense in $O$.

A mapping $\sigma:O\to Q$ between two quasi orders is said 
to \emph{preserve
  \tu(non\tu)atoms} 
if $\sigma(a)$ is a (non)atom whenever $a$ is a (non)atom. 
\end{defn}

\begin{prop}
\label{p-5}
In a Boolean algebra $(B,\le)$, $a\in B$ is an atom 
iff $a\ne0$ and there is no $0<b<a$.
\end{prop}

\begin{remark}
\label{r-4}
In the literature (e.g.~\cite{MR1902334}), 
sometimes the definition of an atom of a lattice is as in the
characterization of proposition~\ref{p-5}, and this may disagree with our
definition. Our terminology better fits the English definition of the word
``atom'', and moreover agrees with its usual usage in set theoretic forcing.
\end{remark}

\begin{example}
\label{x-1}
The atoms of a power set Boolean algebra
$(\power(X),\cup,\cap,\setminus,\allowbreak\emptyset,X)$ are precisely the singletons.
Moreover, the relative atoms of some interval 
$[a,b]=\{x\subseteq X:a\subseteq x\subseteq b\}$ of $\power(X)$ 
are precisely elements of the form $a\cup\{\xi\}$ for some $\xi\in b\setminus a$. 
Thus any interval of a power set Boolean algebra is atomic.
\end{example}

\begin{example}
\label{x-2}
The atoms of a product lattice of the form $\prod_{i\in I}\alpha_i$ 
where each $\alpha_i$ is an ordinal, 
are precisely the members of the form $\xi\cdot\chi_i$ for some $i\in I$ 
and~$0\ne\xi<\alpha_i$. Moreover, if $x\in\prod_{i\in I}\alpha_i$ and
$Y\subseteq\prod_{i\in I}\alpha_i$ is \emph{directed}, i.e.~every two elements of
$Y$ has a common upper bound in $Y$, then the union of intervals $\bigcup_{y\in
  Y}[x,y]$ is a sublattice. Its atoms are members of the form $x+\xi\cdot\chi_i$ for
some $i\in I$ and $0\ne\xi<\alpha_i$. These sublattices are thus atomic.
\end{example}

\begin{example}
\label{x-3}
The category algebra $\cat(X)$ (cf.~\Section\ref{sec:category-algebras}) 
of a Hausdorff space is atomless.
\end{example}

The following results, proposition~\ref{p-8} through corollary~\ref{l-4}, 
are consequences of the fact
(lemma~\ref{u-4}) that simultaneously order preserving and
reflecting maps between quasi orders
come close to being isomorphisms onto their range. Indeed they are
\emph{preisomorphisms} onto their range according to the terminology
of~\cite{MR1365749}.

\begin{prop}
\label{p-8}
Let $(O,\le)$ and $(Q,\altle)$ be quasi orders.
If $\sigma:O\to Q$ is both order preserving and reflecting then $p< q$ iff
$\sigma(p)\altl\sigma(q)$ for all $p,q\in O$.
\end{prop}

\begin{prop}
\label{p-6}
If $\sigma:O\to Q$ is an order preserving and reflecting map between two quasi
orders then the set of minimal elements is
mapped to the set of all relatively minimal elements of the range,
i.e.~$\sigma[O^0]=\ran(\sigma)^0$. 
\end{prop}

\begin{prop}
\label{p-7}
If $\sigma:O\to Q$ is both order preserving and reflecting, then the image of $O^+$
is the positive part of the suborder $\ran(\sigma)$,
i.e.~$\sigma[O^+]=\ran(\sigma)^+$.
\end{prop}

We let $O\div\asym$ denote the antisymmetric quotient, 
i.e.~the equivalence classes modulo $p\sim_\asym q$ if $p\le q$ and $q\le p$. 
Recall that for any quasi order this yields a poset, where the ordering $[p]\le [q]$
if $p\le q$ is well defined since $\sim_\asym$ is a congruence for the quasi order. 

\begin{lem}
\label{u-4}
Let $(O,\le)$ and $(Q,\altle)$ be quasi orders.
Suppose $\sigma$ is an order preserving and reflecting map 
between $(O,\le)$ and $(Q,\altle)$. 
Then $\bar\sigma:O\div\asym\to Q\div\asym$ is well defined by
\begin{equation*}
  \bar\sigma([p])=[\sigma(p)]
\end{equation*}
and is an embedding.
\end{lem}
\begin{proof}
$\bar\sigma$ is well defined because $\sigma$ is order preserving. 
Since $O\div\asym$ is a poset, it remains to show that
$\bar\sigma$ is order preserving and reflecting. But since $\sigma$ is, $[p]\le [q]$
iff $p\le q$ iff $\sigma(p)\le\sigma(q)$ iff
$\bar\sigma([p])=[\sigma(p)]\le[\sigma(q)]=\bar\sigma([q])$. 
\end{proof}

\begin{cor}
\label{l-4}
Let $\sigma:O\to Q$ be both order preserving and reflecting. Then the image of the
atoms are the relative atoms of the image, i.e.~$\sigma[\atom(O)]=\atom(\ran(\sigma))$.
\end{cor}
\begin{proof}
By lemma~\ref{u-4}, $\bar\sigma[\atom(O\div\asym)]=\atom(\ran(\bar\sigma))$. 
The proof is completed by noting that the atoms of the antisymmetric quotient consist
of the equivalence classes of atoms.
\end{proof}
\subsection{Continuity}
\label{sec:continuity}

We examine continuity phenomena for order homomorphisms. 
Then some consequences (proposition~\ref{p-3}--corollary~\ref{l-15}) are
deduced relevant to computing the range of continuous homomorphisms. 

Recall that a subset $A$ of a quasi order is \emph{directed} if it is nonempty and 
every two elements of $A$ have a common upper bound in $A$. In Domain Theory, 
the fundamental topology on a poset is the \emph{Scott topology} whose closed sets
consist of all lower sets that are closed under directed suprema (i.e.~the suprema
of directed subsets), cf.~\cite{MR1365749},\cite{MR1426367}. 
The supremum of a directed set $A$ is normally written as~$\bigsqcup A$. 

\begin{defn}
We call a function $\sigma:P\to Q$ between two posets \emph{continuous} if it is
topologically continuous for the Scott topologies on $P$ and $Q$. It is
\emph{cocontinuous} if it is topologically continuous for the Scott topologies of
the duals of the posets $P$ and $Q$. It is \emph{dually continuous} if it is both
continuous and cocontinuous. 
\end{defn}

\begin{lem}
\label{l-25}
Let $\sigma:P\to Q$ be a function between two posets. Then the following are
equivalent.
\textup{
\begin{enumerate}[leftmargin=*, label=(\alph*), ref=\alph*, widest=b]
\item\label{item:22} \textit{$\sigma$ is continuous.}
\item\label{item:23} \textit{The conjunction of}:
  \begin{enumerate}[leftmargin=*, label=(\arabic*), ref=\arabic*, widest=2]
  \item\label{item:24} \textit{$\sigma$ is a homomorphism \tu(i.e.~order preserving\tu),}
  \item\label{item:25}  \textit{$\sigma$ preserves directed suprema, 
      i.e.~$\sigma\bigl(\bigsqcup A\bigr)=\bigsqcup_{a\in A}\sigma(a)$ 
      whenever $A\subseteq P$ is directed and $\bigsqcup A$ exists.}
  \end{enumerate}
\end{enumerate}}
\end{lem}
\begin{proof}
See e.g.~\cite{MR1365749}. 
\end{proof}

\begin{example}
\label{x-10}
Of course, homomorphisms need not be continuous. Indeed, if $\U$ is a nonprincipal
ultrafilter on $\N$ then $\sigma:\pN\to\pN$ given by $\sigma(A)=\emptyset$ if
$A\notin\U$ and $\sigma(B)=\N$ if $B\in\U$ is a Boolean algebra homomorphism, 
yet it is clearly noncontinuous. 
\end{example}

The following characterization of preservation of nonempty suprema, 
when the domain is a semilattice, is well known. 

\begin{lem}
\label{l-26}
Let $S$ be a semilattice and $Q$ a poset. 
Then the following are equivalent for any $\sigma:S\to Q$.
\textup{
  \begin{enumerate}[leftmargin=*, label=(\alph*), ref=\alph*, widest=b]
  \item\label{item:26} \textit{$\sigma$ is continuous and join-preserving.}
  \item\label{item:27}  \textit{$\sigma$ preserves arbitrary nonempty suprema,
      i.e.~$\sigma\bigl(\spr A\bigr)=\spr_{a\in A}\sigma(a)$ whenever $A\subseteq S$
      is nonempty and $\spr A$ exists.}
  \end{enumerate}}
\end{lem}
\begin{proof}
Assume that $\sigma:S\to Q$ is join-preserving and
continuous with respect to the Scott topologies on $S$ and $Q$.
Suppose $A\subseteq S$ is nonempty and $\spr A$ exists. 
Let $\F$ be the set of all nonempty finite subsets of $A$. 
Then $\bigl\{\spr F:F\in\F\bigr\}$ is a directed subset of $S$ 
with $\bigsqcup_{F\in\F}\spr F=\spr A$,
and thus by lemma~\ref{l-25}(\ref{item:23}\ref{item:25}), 
\begin{equation}
  \label{eq:126}
  \sigma\left(\spr A\right)=\sigma\left(\bigsqcup_{F\in\F}\spr F\right)
  =\bigsqcup_{F\in\F}\sigma\left(\spr F\right)
  =\bigsqcup_{F\in\F}\spr_{a\in F}\sigma(a)
  =\spr_{a\in A}\sigma(a).
\end{equation}
Conversely, if $\sigma$ preserves nonempty suprema then in particular it preserves
directed suprema, and thus is continuous by lemma~\ref{l-25}. 
\end{proof}

\begin{cor}
\label{o-24}
A semilattice homomorphism is continuous iff it preserves non\-empty suprema.
\end{cor}

\begin{cor}
\label{o-25}
A meet semilattice homomorphism is cocontinuous iff it preserves nonempty infima.
\end{cor}
\begin{proof}
By corollary~\ref{o-24} since the dual of a meet semilattice is a (join) semilattice.
\end{proof}

\begin{cor}
\label{o-27}
A lattice homomorphism is dually continuous iff it preserves both nonempty suprema
and infima. 
\end{cor}

\begin{remark}
\label{r-9}
The essential notion for us is preservation of nonempty suprema, as in
lemma~\ref{l-26}(\ref{item:27}). We were tempted to simply define continuity as this
property, but this seemed a bit artificial. Category theoretic continuity does
not quite capture this notion either. In the category of partial orders, all colimits
are coproducts, which are suprema. Thus a mapping is cocontinuous in the category
theoretic sense iff it preserves arbitrary suprema, and not just nonempty ones.

However, all of our example posets are semilattices and as we shall see, the
embeddings under consideration are all semilattice homomorphisms. 
Thus by corollary~\ref{o-24}, for our purposes the two notions coincide.
\end{remark}

\begin{thm}
\label{l-12}
Every order embedding between two posets with a preregular range preserves all 
nonempty suprema. In particular, it is continuous.
\end{thm}
\begin{proof}
Let $(P,\le)$ and $(Q,\altle)$ be two posets, and $\sigma:P\to Q$ an embedding with
preregular range. 
Take $\emptyset\ne A\subseteq P$ such that $\spr A$ exists. 
Then $\sigma\bigl(\spr A\bigr)=\spr_{a\in A}^{\ran(\sigma)}\sigma(a)$ since
$(P,\le)\cong (\ran(\sigma),\altle\nolinebreak)$ via $\sigma$. And by preregularity, 
$\spr_{a\in A}^{\ran(\sigma)}\sigma(a)=\spr_{a\in A}\sigma(a)$. This proves that
$\sigma$ preserves nonempty suprema, and thus it is continuous by lemma~\ref{l-25}.
\end{proof}

\begin{cor}
\label{o-1}
Every order embedding between two posets 
with preregular range preserves nonempty infima.
In particular, it is cocontinuous.
\end{cor}
\begin{proof}
By theorem~\ref{l-12}, 
since any order embedding is also an embedding between the dual posets. 
\end{proof}

Preserving all suprema amounts to preserving all nonempty suprema, 
and mapping the minimum element of the domain, if it has one, 
to the minimum of the codomain. 
We do not want to impose the requirement of preserving minimums in our discourse,
but we will note the following.

\begin{cor}
\label{o-22}
Every order embedding between two posets with regular range preserves arbitrary
suprema and infima, and thus is continuous and cocontinuous in the category
theoretic sense \tu(cf.~remark~\tu{\ref{r-9})}. 
\end{cor}

As mentioned, we will focus on semilattices.

\begin{cor}
\label{o-23}
Every order embedding between two semilattices with preregular range is a continuous
semilattice embedding \tu(and as such, preserves nonempty suprema\tu). 
\end{cor}
\begin{proof}
Theorem~\ref{l-12} and lemma~\ref{l-26}. 
\end{proof}

\begin{cor}
\label{o-12}
Every order embedding between two lattices with preregular range is 
a dually continuous lattice embedding \tu(and thus preserves nonempty suprema and
infima\tu).
\end{cor}

\begin{cor}
\label{o-20}
Every order embedding between two lattices with convex range is a dually continuous
lattice embedding.
\end{cor}
\begin{proof}
Lemma~\ref{l-6} and corollary~\ref{o-12}. 
\end{proof}

Let us point out the obvious, that order isomorphisms are automatically continuous.

\begin{prop}
\label{p-50}
Every order isomorphism between two posets is dually 
Scott continuous and both continuous
and cocontinuous in the categorical sense.
\end{prop}

The remainder of this section is concerned with obtaining regularity properties 
of the range of some homomorphism. 
Sometimes the computation of the range is trivial:

\begin{prop}
\label{p-3}
Suppose that $(O,\le)$ is a quasi order that has both a minimum and 
maximum  element, say $a$ and $b$, respectively.
If $\sigma:O\to Q$ is a quasi order homomorphism 
between $O$ and some quasi order $(Q,\altle)$ with convex range 
then $\ran(\sigma)=[\sigma(a),\sigma(b)]=\{q\in Q:\sigma(a)\altle q\altle \sigma(b)\}$. 
\end{prop}

\begin{lem}
\label{l-17}
Every continuous join-preserving order embedding between a complete semilattice and a
poset has upwards boundedly order closed range.
\end{lem}
\begin{proof}
Let $\sigma:L\to P$ be as in the hypothesis.
Suppose $\emptyset\ne A\subseteq\ran(\sigma)$ is bounded by 
$p\in\ran(\sigma)$ and $\spr^P A$ exists. Then since $\sigma\inv(p)$ is an upper
bound of $\sigma\inv[A]$, $q=\spr\sigma\inv[A]$ exists. 
Thus $\spr^P A=\sigma(q)\in\ran(\sigma)$ as required, by lemma~\ref{l-26}. 
\end{proof}

Let us   also point out a version of the preceding lemma 
that holds for homomorphisms that are not necessarily embeddings.

\begin{lem}
\label{l-23}
Every continuous lattice homomorphism from a complete semilattice into a lattice has
upwards boundedly order closed range.
\end{lem}
\begin{proof}
Letting $\sigma:L\to M$ be as in the hypothesis, suppose $\emptyset\ne
A\subseteq\ran(\sigma)$ is bounded by $p\in\ran(\sigma)$ and $\spr A$ exists.
Choose $c\in L$ such that $\sigma(c)=p$, and for each $a\in A$, 
choose $b_a\in L$ such that $\sigma(b_a)=a$.
Since $L$ is a complete semilattice, we can let $r=\spr_{a\in A}(b_a\land c)$.
Then $\sigma(r)=\spr_{a\in A}\sigma(b_a\land c)
=\spr_{a\in A}(\sigma(b_a)\land p)
=\spr_{a\in A}\sigma(b_a)
=\spr A$, as required.
\end{proof}

\begin{lem}
\label{l-18}
Every cocontinuous meet semilattice homomorphism between a complete semilattice 
and a poset has downwards order closed range.
\end{lem}
\begin{proof}
Let $\sigma:L\to P$ be as in the hypothesis.
Suppose $\emptyset\ne A\subseteq\ran(\sigma)$ and $\ifm A$ exists in $P$. 
Observe that $p=\ifm\sigma\inv[A]$ exists by as
$L$ is a complete semilattice (cf.~\Section\ref{sec:terminology}).
Hence $\sigma(p)=\ifm A$ by corollary~\ref{o-25}, 
and $\sigma(p)\in\ran(\sigma)$ as wanted.
\end{proof}

\begin{cor}\label{o-21}
Every order embedding between a complete semilattice and a poset with preregular
range, in fact has a boundedly order closed range.
\end{cor}
\begin{proof}
Every order embedding between two posets with preregular range preserves nonempty
suprema and infima, and in particular is dually continuous, by theorem~\ref{l-12}
and corollary~\ref{o-1}.
Thus result follows immediately from lemmas~\ref{l-17} and~\ref{l-18}. 
\end{proof}

\begin{defn}
\label{d-6}
A map $f:O\to Q$ between two quasi orders is said to \emph{preserve boundedness} 
if whenever $A\subseteq O$ is bounded (cf.~\Section\ref{sec:terminology}), 
so is $f[A]\subseteq Q$. More generally, a partial function $f:O\pfn Q$ is said to
\emph{preserve boundedness in~$O$} 
if $f[A]$ is bounded whenever $A\subseteq\dom(f)$ is bounded in $O$. 
On the other hand, a function $f:O\to Q$ is said 
to \emph{preserve unboundedness}
if whenever $A\subseteq O$ is unbounded, so is $f[A]\subseteq Q$. 
\end{defn}

\begin{example}
The embedding $i:\omega\to\omega+1$ (with $i(n)=n$) does \emph{not} preserve
unboundedness. 
\end{example}

\begin{cor}
  \label{l-15}
Every order embedding between a complete semilattice and a poset 
preserving unboundedness, and  with preregular range, 
in fact has an order closed range.
\end{cor}
\begin{proof}
Let $\sigma:L\to P$ be an order embedding preserving unboundedness 
with preregular range, where $L$ is a complete semilattice. 
Note that it is in fact established in the proof of corollary~\ref{o-21} that
$\sigma$ has downwards order closed range.
Now suppose $\emptyset\ne A\subseteq\ran(\sigma)$ and $\spr A$ exists in $P$. 
Since $A$ is bounded, so must be $\sigma\inv[A]$ by preservation of unboundedness. 
Therefore, $p=\spr\sigma\inv[A]$ exists. 
But now we have $\spr A=\sigma(p)\in\ran(\sigma)$ as required, 
by theorem~\ref{l-12}.  
\end{proof}
\subsection{Bases of lattices}
\label{sec:bases-lattices}

We turn our attention to various notions of density for subsets of lattices, 
with the goal of finding one suitable for extending order embeddings. 
The technical definition of strong interval predensity is introduced; and we prove
(in lemma~\ref{l-21}) that for a substantial class of commutative monoids,
subsemigroups forming a basis are strongly preinterval dense. 
The motivation for the various
definitions will become clear in the next section \Section\ref{sec:algebraic-domains}.


\begin{defn}
\label{d-11}
A subset $D$ of a quasi order $(O,\le)$ is \emph{dense} 
if every $p\in O^+$ has a $d\le p$ in $D^+$. 
\end{defn}

\begin{notation}
\label{notn:down}
For a quasi order $(O,\le)$, $Q\subseteq O$ and $p\in O$, 
we let $Q_p$ denote the down set $\{q\in Q:q\le p\}$.
\end{notation}

\begin{defn}
A subset $D$ of a poset $(P,\le)$ is \emph{join dense} if every $p\in P$
satisfies
\begin{equation}
  \label{eq:125}
  p=\spr D_p.
\end{equation}
\end{defn}

\begin{defn}
\label{d-4}
A subset $D$ of a poset $P$ is  called \emph{interval predense} if every $p<q$ in $P$
has a $d\in D$ with $d\nleq p$ and $d\le q$. 
\end{defn}

\begin{prop}
\label{p-23}
For any subset of a poset, \emph{join dense} $\to$ \emph{interval predense} $\to$ 
\emph{dense}.
\end{prop}

The following should be compared with the fact that for a Boolean algebra $(B,\le)$,
a subset $D\subseteq B$ is dense iff it is join dense 
(see e.g.~\cite[Ch.~2, \Section4]{MR991565}). 

\begin{lem}
\label{p-14}
Let $L$ be a meet semilattice. 
Then the following are equivalent for all $D\subseteq L$.
\tu{
  \begin{enumerate}[leftmargin=*, label=(\alph*), ref=\alph*, widest=b]
  \item\label{item:20} \textit{$D$ is interval predense.}
  \item\label{item:21} \textit{$D$ is join dense.}
  \end{enumerate}}
\end{lem}
\begin{proof}
By proposition~\ref{p-23}, we need to prove~\eqref{item:20}$\to$\eqref{item:21}.
Fix $p\in L$, and suppose that $q$ is an upper bound of $D_p$. 
We need to show that $p\le q$. 
Supposing to the contrary, $p\land q<p$ and thus there exists $d\nleq
p\land q$ in $D_p$. But this contradicts  the fact that $p\land q$ is an upper bound
of $D_p$.
\end{proof}

The results of section~\Section\ref{sec:infin-distr-laws} are applied here to show
that dense subsemigroups of complete semilattice monoids are join-dense.

\begin{thm}
\label{l-3}
Let $(M,+)$ be a commutative cancellative complete semilattice monoid. 
Then every dense subsemigroup of $M$ is join-dense.
\end{thm}
\begin{proof}
Suppose $D\subseteq M$ is a dense subsemigroup. Take $p\in M$.
The supremum $q=\spr D_p$ exists 
because we are dealing with a complete semilattice,
and obviously $q\le p$; hence, there is an $r\in M$ such that $q+ r=p$. 
It remains to show that $r=0$. Supposing to the contrary, by denseness
there exists $s\in D^+$ with $s\le r$. 
By corollary~\ref{o-10} we can use infinite right $(+,\lor)$-distributivity, to obtain
\begin{equation}
  \label{eq:59}
  q+ s=\spr(D_p+ s).
\end{equation}
However, $D_p+ s\subseteq D$ since $D$ is a subsemigroup. Thus in fact
$D_p+ s\subseteq D_p$ because $d+ s\le q+ s$ for all $d\in D_p$ by
proposition~\ref{p-18}, and $q+ s\le p$ by proposition~\ref{p-4}. 
Therefore, $q+ s\le q$ by~\eqref{eq:59}, which implies $s=0$ by cancellativity and
lemma~\ref{p-15}. However, $0\notin D^+$ (see remark~\ref{r-5}), a contradiction.
\end{proof}

\begin{defn}
We call a subposet $A$ of a poset $(P,\le)$ \emph{flat} if there exists $p\in P$
such that $a\land b=p$ for all $a\ne b$ in $A$. A lattice $(L,\le)$ is called
\emph{flat-complete} if $\spr A$ exists for every flat $A\subseteq L$. 
\end{defn}

\begin{example}
\label{x-17}
Since $(\N,\le)$ is a chain, it is trivially flat-complete because a flat subset of
chain has at most two elements. Thus for any index set $I$, $(\irri I,\le)$ is
flat-complete because flat-completeness is preserved under products. 
\end{example}

The notion of being a complete semilattice is incomparable with being
flat-complete. 

\begin{example}
\label{x-16}
We consider two sublattices of the complete semilattice $[0,\infty)^\N$ 
with the product order, which is also flat-complete being a product of chains. 
The sublattice $(\rationals\cap[0,\infty))^\N$ is not a complete semilattice since
the nonnegative rationals form a lattice that is not a complete semilattice, 
but it is flat-complete.
On the other hand the family of $B([0,\infty)^\N)$ of bounded sequences 
of nonnegative reals is clearly a complete semilattice, but not flat-complete.
\end{example}

\noindent However, in the case of Boolean algebras the following is well known (and
easily proved). 

\begin{prop}
\label{p-37}
Every flat-complete Boolean algebra is a complete Boolean algebra. 
\end{prop}

\begin{defn}
\label{d-5}
Let $(L,\le)$ be a pointed lattice (cf.~\Section\ref{sec:terminology}). 
A \emph{basis} for $L$ is a meet subsemilattice
$B\subseteq L$ such that every $a\in L$ has a family $A\subseteq B$ satisfying
\begin{enumerate}[leftmargin=*, label=(\roman*), widest=ii]
\item\label{item:6} $A$ is pairwise incompatible, 
i.e.~$b\land c=0$ for all $b\ne c$ in $A$,
\item\label{item:19} $\spr A=a$. 
\end{enumerate}
\end{defn}

\begin{remark}
In set theoretic terminology pairwise incompatible families are called 
\emph{antichains}; however, this disagrees with its  usage elsewhere, 
as a pairwise incomparable set (e.g.~\cite{MR1902334},\cite{MR1429390}). 
\end{remark}

\begin{lem}
\label{l-24}
Every dense meet subsemilattice of a complete Boolean algebra is basis. 
\end{lem}
\begin{proof}
Let $D$ be dense. Fixing $a\in D$, let $A\subseteq D$ be a maximal pairwise
incompatible family below $a$. Then $\spr A$ exists by completeness, 
and we cannot have $\spr A<a$ or else there exists $b\le a-\spr A$ in $D^+$, 
and thus $b\notin A$ contradicting its maximality.
\end{proof}

\begin{example}
\label{x-19}
Clopen sets are in particular regular open sets.
Thus the collection of equivalence classes of the clopen subsets of $X$ forms a
sublattice of $\cat(X)$.
Let $X$ be a zero dimensional topological space. Then the family of clopen sets 
identifies with a dense subset of $\cat(X)$, and therefore it forms a basis by
lemma~\ref{l-24}. 
\end{example}

\begin{defn}
\label{d-8}
We will say that a subset $D$ of a lattice $(L,\le)$ 
is \emph{strongly interval predense} 
if for every $p<q$ in $L$ there exists $d\in D$ such that 
\begin{enumerate}[leftmargin=*, label=(\roman*), ref=\roman*, widest=iii]
\item $d\le q$,
\item\label{item:28} $d\land p<d$,
\item\label{item:5} $d\land p\in D$.
\end{enumerate}
This is equivalent to adding condition~\eqref{item:5} to interval predensity. 
\end{defn}

Being join-dense does not entail being strongly interval predense.

\begin{example}
\label{x-15}
Let $\varTheta$ be some initial segment of the ordinals with its usual linear ordering
$\in$. Then $(\varTheta,\in)$ is a complete semilattice, and is moreover a complete
lattice when it is closed (e.g.~$\varTheta=\omega+1$). 
Let $S_\varTheta$ be the set of all successor ordinals~in~$\varTheta$. 
Then $S_\varTheta$ is evidently join-dense and thus interval predense, 
but so long as $\omega+1\in\varTheta$, it is not strongly interval
predense, because $\omega<\omega+1$ but there is no successor ordinal satisfying
both~\eqref{item:28} and~\eqref{item:5} for $p=\omega$. 
\end{example}

Strong interval predensity does not entail preregularity.

\begin{example}
\label{x-18}
Consider the complete Boolean algebra $(\power(X),\subseteq)$ where $X$ is some
infinite set. Then $\Fin(X)$, the set of all finite subsets of $X$ forms a
sublattice, and thus $\Fin(X)\cup\{X\}$ is a bounded sublattice.
It is strongly interval predense, because for any $a\subset b$ in $\power(X)$,
picking $x\in b\setminus a$, $\{x\}\in\Fin(X)\cup\{X\}$, $\{x\}\subseteq b$ and
$\{x\}\cap a=\emptyset\subset\{x\}$ is in $\Fin(X)$.
However, $\Fin(X)\cup\{X\}$ is not preregular, because taking any infinite
$y\subseteq X$, $\spr^{\Fin(X)\cup\{X\}}\Fin(y)=X$ whereas $\spr^{\power(X)}\Fin(y)=y$.
\end{example}

\begin{lem}
\label{l-21}
Let $(M,+)$ be a commutative cancellative complete semilattice mo\-noid
that is moreover \jid. Then every subsemigroup that is a basis is also strongly
interval predense. 
\end{lem}
\begin{proof}
Suppose that $B$ is a subsemigroup that is a basis. 
Take $a<b$ in $M$, and let $c\in M$ satisfy $a+c=b$. Hence $c\ne 0$.
Since $B$ is a basis there is a pairwise incompatible family $A\subseteq B$
with supremum $a$. Fix any $d\in A$.
First we consider the case $c\land d=0$.
Since $B$ is dense there is a $\bar d\in B^+_c$. 
Then $d+\bar d\in B$  because $B$ is a subsemigroup. 
And theorem~\ref{l-7} implies that $M$ is $(+,\lorand)$-distributive, 
and hence $d+\bar d=d\lor\bar d$ by lemma~\ref{p-33}.
Now  $d+\bar d$ witnesses that $B$ is strongly interval predense
because $d+\bar d\le a+c=b$, 
$(d+\bar d)\land a=(d\lor\bar d)\land a=(d\land a)\lor(\bar d\land a)=d\in B$ by
distributivity, and $d<d+\bar d$.

Otherwise, when $c\land d\ne 0$, we choose $\bar d\in B_{c\land d}^+$. 
Then $d+\bar d\in B$, $d+\bar d\le a+c=b$, and using \jid,
$(d+\bar d)\land a=(d+\bar d)\land\spr A=\spr_{p\in A}(d+\bar d)\land p
=(d+\bar d)\land d=d$ by lemma~\ref{l-20}, as required.
\end{proof}
\subsection{Extensions of embeddings}
\label{sec:algebraic-domains}

In this section we are concerned with extending some homomorphism 
from a join-dense subset of a lattice to the entire lattice, 
while preserving additional properties. Two of the main results are
theorem~\ref{o-28} which allows for the extension of a lattice embedding of a
strongly interval predense basis to a continuous lattice embedding of the entire
lattice, and theorem~\ref{o-19} which allows us to do the same while also
maintaining a convex range. We also introduce the notion of the $P$-Scott topology
on a subposet of some poset $(P,\le)$. 

\subsubsection{The subposet topology}\label{sec:subposet-topology}
At this juncture we need to look more closely at the topology on a subposet $Q$ of
$(P,\le)$. The point is that the Scott topology on $Q$ as given by its ordering
$\le$, may differ from the subspace topology that $Q$ inherits from the Scott
topology on $P$, and moreover the two topologies may be incomparable. 
In fact, we are interested in a third topology on $Q$. 

\begin{defn}
\label{d-9}
Let $(P,\le)$ be a poset and $Q\subseteq P$. We say that $A\subseteq Q$ is
\emph{$P$-closed under directed suprema} if for every directed $A\subseteq Q$,
if $a=\bigsqcup^P A$ exists and is in $Q$, then $a\in A$. Similarly, $A\subseteq Q$ is
said to be \emph{upwards $P$-order closed}  if: 
for all $\emptyset\ne B\subseteq A$, if $a=\spr^P B$ exists and is in $Q$
then $a\in A$. The dual notions are called \emph{$P$-closed under downwards directed
  infima} and
\emph{downwards $P$-order closed}, respectively.
\end{defn}

\begin{prop}
\label{p-49}
If $F\subseteq P\supseteq Q$ is closed under directed suprema 
\tu(downwards directed infima\tu), 
then $F\cap Q$ is $P$-closed under directed suprema 
\tu(downwards directed infima\tu).
\end{prop}

\begin{prop}
\label{p-42}
If $F\subseteq P\supseteq Q$ is upwards \tu(downwards\tu) order closed in $P$,
then $F\cap Q$ is upwards \tu(downwards\tu) $P$-order closed 
as a subposet of $(Q,\le)$.
\end{prop}

\begin{prop}
\label{p-48}
If $A\subseteq Q\subseteq P$ is closed under directed suprema 
\tu(downwards directed infima\tu) in the poset $(Q,\le)$, 
then it is also $P$-closed under directed suprema \tu(downwards directed infima\tu).
\end{prop}
\begin{proof}
By proposition~\ref{p-28}.
\end{proof}

\begin{prop}
\label{p-39}
If $A\subseteq Q\subseteq P$ is upwards \tu(downwards\tu) order closed 
as a subposet of $(Q,\le)$, 
then it is also upwards \tu(downwards\tu) $P$-order closed. 
\end{prop}
\begin{proof}
By proposition~\ref{p-28}. 
\end{proof}

Notice that if $Q$ is a preregular subset of $P$, then $P$-closed coincides with
closed. Thus, e.g.~%

\begin{prop}
\label{p-44}
Let $Q$ be a preregular subset of $P$. 
Then $A\subseteq Q$ is closed under directed suprema iff it is $P$-closed under
directed suprema. 
\end{prop}

\begin{defn}
Let $(P,\le)$ be a poset and $Q\subseteq P$. 
Then the family of complements of lower subsets of $Q$ that are $P$-closed under
directed suprema is a topology on $Q$; 
we call it the \emph{$P$-Scott topology} on $Q$. 
A function $\sigma:Q\to R$, where $R$ is some poset, is called
\emph{$P$-continuous} if it is a continuous function with respect to the $P$-Scott
topology on $Q$ and the Scott topology on $R$. 
\end{defn}

\begin{prop}
\label{p-38}
The $P$-Scott topology on a subposet $Q\subseteq P$ 
is a common refinement of the its Scott topology and its subspace topology.
\end{prop}
\begin{proof}
The $P$-Scott topology refines the Scott topology by proposition~\ref{p-48},
and it refines the subspace topology by proposition~\ref{p-49}. 
\end{proof}

Note that for $Q=P$, the $P$-Scott topology is just the Scott topology. 
More generally:

\begin{prop}
\label{p-40}
If $Q\subseteq P$ is preregular then the Scott topology on $Q$ coincides with its
$P$-Scott topology.
\end{prop}
\begin{proof}
By proposition~\ref{p-44}. 
\end{proof}

Corresponding to lemma~\ref{l-25} we have:

\begin{prop}
\label{p-46}
Let $Q$ be a subposet of $(P,\le)$, and $(R,\altle)$ some poset. 
Then any function $\sigma:Q\to R$ is $P$-continuous iff it is a homomorphism
that preserves directed suprema in $P$, i.e.~$\sigma\bigl(\bigsqcup^P
A\bigr)=\bigsqcup_{a\in A}\sigma(a)$ whenever $A\subseteq Q$ is directed and
$\bigsqcup^P A$ exists and is in $Q$. 
\end{prop}

And corresponding to lemma~\ref{l-26}:

\begin{prop}
\label{p-47}
Suppose $(S,\le)$ is a semilattice and $Q\subseteq S$ is a subsemilattice, and
suppose $(R,\altle)$ is a poset. Then any function $\sigma:Q\to R$ is $S$-continuous
and join-preserving iff it preserves nonempty suprema in $P$,
i.e.~$\sigma\bigl(\spr^P A\bigr)=\spr_{a\in A}\sigma(a)$ 
whenever $A\subseteq Q$ and $\spr^P A$ exists and is in $Q$. 
\end{prop}

The following lemma shows that for complete semilattices,
convexity of the range of an embedding follows from its convexity when suitably
restricted to an interval predense subset of its domain.

\begin{lem}
\label{l-16}
Let $L$ and $M$ be complete semilattices, with $M$ satisfying the \jid,
and let $\sigma:L\to M$ be a continuous semilattice embedding.
Suppose $0\in D\subseteq L$ and $E\subseteq M$ 
are both join-dense \tu(equivalently, interval predense\tu), $E$ is a  subsemilattice,
and that $\sigma[D]\subseteq E$ is a convex subset of $E$. 
Then $\ran(\sigma)$ is convex. 
\end{lem}
\begin{proof}
Given $p\le q$ in $\ran(\sigma)$, 
we choose $p\le r\le q$ and prove that $r\in\ran(\sigma)$. Let
\begin{equation}
  \label{eq:99}
  s=\spr\{c\in\sigma[D]:c\le r\},
\end{equation}
and note that since $0\in D$,
$p\in\ran(\sigma)$ implies that the set on the right is nonempty, 
namely $\sigma(0)\le p\le r$ is in $\sigma[D]$.
Observe that $s\in\ran(\sigma)$ because $\ran(\sigma)$ is upwards boundedly order
closed by lemma~\ref{l-17}. Therefore it suffices to prove $r=s$.

Suppose to the contrary that $s<r$. 
Since $D$ is join-dense,
$\sigma\inv(q)=\spr\{d\in\nobreak D:d\le\sigma\inv(q)\}$,
and thus by corollary~\ref{o-24}, $q=\spr\{c\in\sigma[D]:c\le q\}$. Now by the~\jid, 
\begin{equation}
  \label{eq:101}
  r=r\land q=\spr\{r\land c: c \in\sigma[D]\text{, }c\le q\}.
\end{equation}
Hence there exists $c\in\sigma[D]$ such that $s\land c=s\land (r\land c)< r\land c$.
Since $E$ is interval predense, there exists $d\in E$ such that $d\le r\land c$
but 
\begin{equation}
  d\nleq s\land c.\label{eq:100}
\end{equation}
We now have $d'=d\lor \sigma(0)\in E$ since $E$ is a subsemilattice, 
and thus as $\sigma(0)\le d'\le c$, $d'\in\sigma[D]$ by convexity. 
But then~\eqref{eq:99} says that $d'\le s$, which would imply $d\le s\land c$, 
contradicting~\eqref{eq:100}.
\end{proof}

Preservation of boundedness (cf.~definition~\ref{d-6}) 
allows one to extend continuous order 
homomorphisms from join-dense subsets.
However, embeddings do not necessarily extend. 
For example consider the continuous partial embedding $j:\omega+2\pfn \omega+1$ 
where $j(n)=n$ for $n<\omega$ and $j(\omega+1)=\omega$; 
the unique extension to a continuous homomorphism on
$\omega+2$ maps $\omega$ to $\omega$ and thus is not an injection. 

\begin{lem}
\label{l-19}
Let $L$ be a \jid\ lattice, and let $M$ be a complete semilattice.
Suppose $D\subseteq L$ is a join-dense meet subsemilattice.
Then every function $\sigma:D\to M$ that preserves nonempty suprema in $L$ 
\tu(cf.~proposition~\tu{\ref{p-47})},
and preserves boundedness in $L$,
has an extension to a continuous semilattice homomorphism on $L$, 
that is uniquely determined on $L^+$.
And this extension is moreover a lattice homomorphism when $\sigma$ is 
meet-preserving and $M$ is~\jid. 
\end{lem}
\begin{proof}
Let $\sigma:D\to M$ be as specified in the hypothesis.
We can define $\bar\sigma:L\to M$ by
\begin{equation}
  \label{eq:72}
  \bar\sigma(p)=\spr_{d\in D_p}\sigma(d)
\end{equation}
because the supremum exists as $\sigma[D_p]$ is bounded.
$\bar\sigma(d)=\sigma(d)$ for all $d\in D$ because $\sigma$ is order preserving, 
and it is also clear that $\bar\sigma$ is order preserving. 
Suppose $\emptyset\ne A\subseteq L$ and $\spr A$ (i.e.~$\spr^L A$) exists. 
Clearly  $\spr_{p\in A}\bar\sigma(p)\le\bar\sigma\bigl(\spr A\bigr)$; 
while every $d\in D_{\spr\!A}$ 
satisfies $\spr_{p\in A}^L(d\land D_p)=d\land\spr^L_{p\in A}D_p=d\land \spr A=d$ 
since $L$ is \jid\ and $D$ is join-dense, 
and thus $\spr_{p\in A}\bar\sigma(p)=\spr_{p\in A}\spr_{d'\in D_p}\sigma(d')
\ge\spr_{p\in A}\spr_{d'\in D_p}\sigma(d\land d')=\sigma(d)$
since $D$ is a meet subsemilattice and since $\sigma$ preserves nonempty suprema in
$L$. This proves the other inequality,  establishing that $\bar\sigma$ preserves
nonempty suprema. Thus $\bar\sigma$ is a continuous semilattice homomorphism by
corollary~\ref{o-24}. 

\begin{claim}
\label{c-12}
If $\sigma$ is meet-preserving and $M$ is \jid\ then $\bar\sigma$ is
meet-preserving, and thus is a lattice homomorphism.
\end{claim}
\begin{proof}
Using the fact that $D_{p\land q}=D_p\land D_q$ 
since $D$ is a meet subsemilattice for the second equality, and the fact that
$\sigma$ is meet-preserving and proposition~\ref{p-27} for the third, 
\begin{equation}
\label{eq:123}
\begin{split}
\bar\sigma(p\land q)&=\spr_{d\in D_{p\land q}}\sigma(d)\\
&=\spr_{c\in D_p}\spr_{d\in D_q}\sigma(c\land d)\\
&=\Biggl(\spr_{c\in D_p}\sigma(c)\Biggr)\land\Biggl(\spr_{d\in D_q}\sigma(d)\Biggr)\\
&=\bar\sigma(p)\land\bar\sigma(q).\\[-19pt]
\end{split}
\end{equation}
\end{proof}

For all $p\in L$, if $D_p\ne\emptyset$ then equation~\eqref{eq:72} must hold for any
continuous extension of $\sigma$ to a semilattice homomorphism, 
proving uniqueness on $L^+$.
\end{proof}

\begin{remark}
\label{r-7}
We do not necessarily have a unique extension to a pointed lattice 
$L$ unless $0$ is already in $D$.
For example, $\N^+=\{1,2,\dots\}$ is join-dense in $\N$, but $n\mapsto n+1$ has two
different extensions to a continuous lattice homomorphism. Of course, this freedom
at $0$ disappears if we require that all suprema are preserved (i.e.~category
theoretic continuity). 
\end{remark}

In our intended applications, we are extending from a join-dense sublattice, and we
can then state the preceding result in the following simplified form.

\begin{cor}
\label{o-26}
Let $L$ be a \jid\ lattice, and $M$ be a complete semilattice. 
Suppose $D\subseteq L$ is a join-dense sublattice. 
Then every $L$-continuous semilattice homomorphism $\sigma:D\to M$, 
that preserves boundedness in $L$, has an
extension to a continuous semilattice homomorphism on $L$, 
uniquely determined on $L^+$. This extension is
moreover a lattice homomorphism when $\sigma$ is a lattice homomorphism and $M$ 
is~\jid. 
\end{cor}
\begin{proof}
By lemma~\ref{l-19} and proposition~\ref{p-47}.
\end{proof}

We use flat-completeness to ensure preservation of boundedness.

\begin{lem}
\label{l-9}
Let $L$ be a lattice that is \jid, and let $M$ be a flat-complete lattice.
Suppose $B\subseteq L$ is a basis. 
Then every meet-preserving order homomorphism between $B$
and $M$, that preserves nonempty suprema in $L$,
preserves boundedness in $L$.
\end{lem}
\begin{proof}
Let $\sigma:B\to M$ be as specified by the hypothesis.
Given $0\ne p\in L$ we must show that $\sigma[B_p]$ is bounded in $M$.
Find a pairwise incompatible $A\subseteq B_p$ such that $\spr A=p$. 
Since $\sigma$ is meet-preserving, 
$\sigma(a)\land\sigma(b)=\sigma(0)$ for all $a\ne b$ in $A$, and thus $\sigma[A]$ is flat.
Therefore $\spr_{a\in A}\sigma(a)$ exists by flat-completeness.
It now suffices to prove that $\spr_{a\in B_p}\sigma(a)=\spr_{a\in A}\sigma(a)$.
But for any $c\in B_p$, $c=c\land p=\spr_{a\in A}^L c\land a$ by the~\jid, and thus
as $A\ne\emptyset$,
by preservation of nonempty suprema in $L$,
$\sigma(c)=\spr_{a\in A}\sigma(c\land a)\le\spr_{a\in A}\sigma(a)$, as
wanted. 
\end{proof}

Combining corollary~\ref{o-26} with the preceding lemma, we obtain the following
result.

\begin{cor}
\label{o-30}
Let $L$ be a lattice, and $M$ be a complete and flat-complete semilattice, with both
\jid. 
Suppose $B\subseteq L$ is a sublattice that forms a basis. 
Then every $L$-continuous lattice homomorphism $\sigma:B\to M$ has an
extension to a continuous lattice homomorphism on $L$, unique on $L^+$.
\end{cor}
\begin{proof}
Note that bases are obviously join-dense. 
\end{proof}

Lattice homomorphisms with codomain a preregular subset of $M$ can also be extended.

\begin{cor}
\label{o-31}
Let $L$ be a lattice, and $M$ be a complete and flat-complete semilattice, with both
\jid. Suppose $B\subseteq L$ is a sublattice forming a basis, and $E\subseteq M$ is
a preregular subset of $M$. Then every $L$-continuous order homomorphism
$\sigma:B\to E$ that is
both join and meet-preserving has an extension to a continuous lattice
homomorphism $\bar\sigma:L\to M$, unique on $L^+$. 
\end{cor}
\begin{proof}
Letting $\sigma':B\to M$ be the same function as $\sigma$ but with codomain $M$,
$\sigma'$ is a lattice homomorphism because $E$ is preregular. 
And by propositions~\ref{p-38} and~\ref{p-40}, the Scott topology on $E$ refines its
subspace topology, and thus $\sigma'$ is also $L$-continuous. 
Now the result is immediate from corollary~\ref{o-30}. 
\end{proof}

The notion of a strongly interval predense subset (cf.~definition~\ref{d-8}) was
introduced to allow  boundedness preserving embeddings to be extended.

\begin{cor}
\label{o-16}
Let $L$ be a lattice and $M$ be a complete semilattice with both satisfying the \jid.
Suppose $D\subseteq L$ is a strongly interval predense meet subsemilattice.
Then every meet subsemilattice embedding of $D$ into $M$ that preserves nonempty
suprema in $L$, and that preserves boundedness in $L$,
has an extension to a continuous lattice embedding of $L$ into $M$, 
uniquely determined on $L^+$.
\end{cor}
\begin{proof}
$D$ is join-dense by lemma~\ref{p-14}.
Hence by lemma~\ref{l-19}, 
there is a continuous lattice homomorphism $\bar\sigma:L\to M$
extending~$\sigma$, uniquely determined by equation~\eqref{eq:72} on $L^+$.
It remains to verify that $\bar\sigma$ is an embedding. 

By proposition~\ref{p-34}, it suffices to show that $\bar\sigma$ is
strictly order preserving. Take $p<q$ in $L$, and suppose towards a contradiction
that $\bar\sigma(p)=\bar\sigma(q)$. 
Since $D$ is strongly interval predense, 
there exists $d\in D_q$ such that $d\land p<d$ and $d\land p\in D$.
But then since $\sigma$ is an embedding, we have 
\begin{equation}
\label{eq:124}
\begin{split}
\bar\sigma(d\land p)&=\sigma(d\land p)\\
&<\sigma(d).
\end{split}
\end{equation}
However, by assumption, $\bar\sigma(d\land p)=\bar\sigma(d)\land\bar\sigma(p)
=\sigma(d)\land\bar\sigma(q)=\sigma(d)$, contradicting~\eqref{eq:124}.
\end{proof}

Now we come to the main result of this section, 
that embeddings into a complete flat-complete semilattice 
can be extended from a strongly interval predense basis.

\begin{thm}
\label{o-28}
Let $L$ be a lattice and $M$ a complete semilattice that is flat-complete, with both
\jid. Suppose $B\subseteq L$ is a strongly interval predense basis. 
Then every meet subsemilattice embedding of $B$ into $M$ that preserves nonempty
suprema in $L$, has an extension to a continuous lattice embedding of $L$ into $M$,
unique on $L^+$. 
\end{thm}
\begin{proof}
Let $\sigma:B\to M$ be as hypothesized. 
Then $\sigma$ preserves boundedness in $L$ by lemma~\ref{l-9}.
The result now follows from corollary~\ref{o-16}. 
\end{proof}

We obtain a very natural simplified form when $B$ is a sublattice.

\begin{cor}
\label{o-29}
Let $L$ be a lattice and $M$ a complete semilattice that is flat-complete, with both
\jid. Suppose $B\subseteq L$ is a sublattice that forms a strongly interval
predense basis. 
Then every $L$-continuous lattice embedding of $B$ into $M$ has an extension to a
continuous embedding of $L$ into $M$, unique on $L^+$.
\end{cor}
\begin{proof}
By theorem~\ref{o-28} and proposition~\ref{p-47}.
\end{proof}

Embeddings with preregular range satisfy the hypothesis of the main theorem.

\begin{cor}
\label{u-10}
Let $L$ be a \jid\ lattice, and let $M$ be a \jid\ complete semilattice that is also
flat-complete.
Suppose $B\subseteq L$ is a strongly interval predense basis.
Then every order embedding of $B$ into $M$ with preregular range
has an extension to a continuous lattice embedding of $L$ into $M$, 
unique on $L^+$.
\end{cor}
\begin{proof}
Let $\sigma:B\to M$ be an order embedding with convex range.
By theorem~\ref{l-12}, $\sigma$ preserves arbitrary nonempty suprema, and thus in
particular preserves nonempty suprema in $L$. And corollary~\ref{o-1} implies that
$\sigma$ is meet-preserving. 
Therefore $\sigma$ extends to a continuous lattice embedding of $L$ into $M$,
uniquely on $L^+$, by theorem~\ref{o-28}.
\end{proof}

We shall want to extend an embedding from a subsemigroup forming a basis.

\begin{cor}
\label{o-33}
Let $(L,+)$ be a commutative cancellative \jid\ complete semilattice monoid,
and let $M$ be a \jid\ complete flat-complete semilattice. Suppose $B\subseteq L$ is
a subsemigroup that forms a basis.
Then every meet subsemilattice embedding of $B$ into $M$ preserving nonempty 
suprema in $L$, has an extension to a continuous lattice embedding of $L$ into $M$,
unique on $L^+$. 
\end{cor}
\begin{proof}
By theorem~\ref{o-28}, because $B$ is strongly interval predense by lemma~\ref{l-21}.
\end{proof}

The following result allows one to extend a continuous embedding from a
basis of a complete semilattice to the whole semilattice while maintaining a convex
range. 

\begin{thm}
\label{o-19}
Let $L$ and $M$ be complete semilattices that are \jid, with $M$ flat-complete.
Suppose $0\in B\subseteq L$ is a strongly interval predense basis 
and $E\subseteq M$ is a join-dense preregular sublattice. 
Then every order embedding of $B$ into $E$ with convex range \tu(in $E$\tu) has a
unique extension to a continuous lattice embedding $\sigma:L\to M$, 
and moreover $\ran(\sigma)$ is convex. 
\end{thm}
\begin{proof}
Let $\sigma:B\to E$ be an order embedding with convex range.
Since $E$ is a lattice, $\ran(\sigma)$ is preregular in $E$ by lemma~\ref{l-6}.
And thus by proposition~\ref{p-54}, $\ran(\sigma)$ is a preregular subset of $M$. 
Now applying corollary~\ref{u-10}, 
$\sigma$ has an extension, say $\bar\sigma:L\to M$, to a continuous lattice
embedding of $L$ into $M$, uniquely determined on $L^+$. 
Hence $\bar\sigma$ is in fact the unique extension of
$\sigma$ to a continuous lattice homomorphism by remark~\ref{r-7}. 
And as the hypotheses of lemma~\ref{l-16} are all satisfied, $\ran(\bar\sigma)$
is convex. 
\end{proof}

\begin{cor}
\label{o-32}
Let $(L,+)$ be a commutative cancellative complete semilattice monoid,
and let $M$ be a complete flat-complete semilattice, with both \jid.
Suppose $B\subseteq L$ is a submonoid forming a basis 
and $E\subseteq M$ is a join-dense preregular sublattice. 
Then every order embedding $\sigma:B\to E$ with convex range has a
unique extension to a continuous lattice embedding of $L$ into $M$ with convex range.
\end{cor}
\begin{proof}
Theorem~\ref{o-19} and lemma~\ref{l-21}.
\end{proof}

Our intended application uses the following corollary.

\begin{cor}
\label{o-34}
Let $(L,+)$ be a commutative cancellative complete semilattice monoid, and let
$(M,+)$ be a commutative cancellative complete and flat-complete semilattice monoid,
with both \jid. Suppose $B\subseteq L$ is a submonoid forming a basis and
$E\subseteq M$ is a dense submonoid closed under subtraction.
Then every order embedding $\sigma:B\to E$ with convex range has a unique extension
to a continuous lattice embedding of $L$ into $M$ with convex range. 
\end{cor}
\begin{proof}
By lemma~\ref{l-27}, $E$ is regular and in particular preregular, and thus is also a
sublattice (corollary~\ref{o-35}). Now the result is an immediate consequence of
corollary~\ref{o-32}. 
\end{proof}
\section{The partial orders}
\label{sec:partial-orders}

A precise description is obtained of the partial order embeddings with convex range
within the following classes of posets: power set algebras; powers of $\N$; category
algebras; Baire functions with codomain a power of $\N$, modulo almost always
equality; and
continuous functions with codomain a power of $\N$, in
subsections~\ref{sec:new-subsection}--\ref{sec:cont-funct-into}, respectively.
In \Section\ref{sec:further-directions}, these embeddings are discussed for the
following classes: Baire functions with codomain a power set algebra, modulo almost
always equality; measurable functions with respect to some measure space, modulo
almost everywhere equality; the Boolean algebra $\pnfin$; and the lattice quotient
$\irrfin$ of the irrationals ordered by eventual dominance. 

Embeddings with convex range are considered exclusively. The results here fail
otherwise as demonstrated in example~\ref{x-11}.
\subsection{Power set algebras}
\label{sec:new-subsection}

We begin by considering order embeddings between power set algebras. 

\begin{thm}
\label{p-12}
Order embeddings $\sigma$ from
$(\power(X),\allowbreak \subseteq)$
into $(\power(Y),\subseteq)$ with convex range consist precisely of maps of the form 
\begin{equation}
  \label{eq:48}
  \sigma(a)=h[a]\cup b\espc\tu{for all $a\subseteq X$},
\end{equation}
for some injection $h:X\to Y$ and some $b\subseteq Y$ with $h[X]\cap b=\emptyset$.  
All of these  maps are moreover dually continuous lattice embeddings.
Furthermore, $h$ is a bijection iff 
$\sigma$ is a Boolean algebra isomorphism
between $\power(X)$ and $\power(Y)$. 
\end{thm}
\begin{proof}
We have $\ran(\sigma)=[\sigma(\emptyset),\sigma(X)]$ by proposition~\ref{p-3}.
Thus for every singleton, 
$\sigma(\{x\})=\sigma(\emptyset)\cup\{y_x\}$ 
for some $y_x\in\sigma(X)\setminus\sigma(\emptyset)$, 
by corollary~\ref{l-4} and example~\ref{x-1}. By corollary~\ref{o-20}, $\sigma$ is
dually continuous. In particular, continuity implies by lemma~\ref{l-26} that
$\sigma(a)=\sigma(\emptyset)\cup\{y_x:x\in a\}$ for all $a\subseteq X$. 
Thus the desired $h$ is defined by
\begin{equation}
  \label{eq:11}
  h(x)=y_x,
\end{equation}
with $b=\sigma(\emptyset)$.
Conversely, any $\sigma$ of the above form is easily seen to be an embedding
whenever $h$ is an injection and $h[X]\cap b=\emptyset$. 
In the case where $h$ is a bijection, $h[X]=Y$ implies $b=\emptyset$ as $h[X]\cap
b=\emptyset$, 
and thus $\ran(\sigma)=\power(Y)$ implies $\sigma$ is moreover an
isomorphism. 
\end{proof}

\begin{cor}
\label{u-1}
Every map $\sigma:\power(X)\to\power(Y)$ satisfying\tu{: 
\begin{enumerate}[leftmargin=*, label=(\alph*), ref=\alph*, widest=b]
\item\label{item:3} \textit{$a\subseteq b$ iff $\sigma(a)\subseteq\sigma(b)$, 
for all $a,b\subseteq X$,}
\item \label{item:4} \textit{$\sigma(a)\subseteq c\subseteq\sigma(b)$ 
implies $c\in\ran(\sigma)$, for all $a,b\subseteq X$ and $c\subseteq Y$,}
\end{enumerate}}
\noindent is of the form
\begin{equation*}
  \sigma(a)=h[a]\cup\sigma(\emptyset)\espc\tu{for all $a\subseteq X$,}
\end{equation*}
for some injective $h:X\to Y$. 
\end{cor}
\begin{proof}
Theorem~\ref{p-12} is applied to the posets $(\power(X),\subseteq)$ and
$(\power(Y),\subseteq)$, as conditions~\eqref{item:3} and~\eqref{item:4} translate
to $\sigma$ being an embedding with convex range. 
\end{proof}

\begin{cor}
\label{o-11}
The range of any embedding $\sigma:\power(X)\to\power(Y)$ with convex range is the
interval Boolean algebra $[\sigma(\emptyset),\sigma(X)]$, 
and $i\inv\circ\sigma$ is in fact a Boolean algebra
isomorphism, where $i\inv$ any left inverse of the inclusion function
$i:[\sigma(\emptyset),\sigma(X)]\to\power(Y)$. 
\end{cor}

The following result is well known, e.g.~it is mentioned in~\cite[pp.~171]{MR1623206}.

\begin{cor}
\label{o-5}
Every quasi order isomorphism $\sigma:\power(X)\to\power(Y)$ is of the form
$\sigma(a)=h[a]$ where $h:X\to Y$ is a bijection. Thus order automorphisms of
any power set algebra $\power(X)$ are always given by 
the image of a permutation of $X$. 
\end{cor}

\begin{example}
\label{x-11}
Consider the two and three element sets $\{0,1\}$ and $\{0,1,2\}$. 
Then let $\sigma:\power(\{0,1\})\to\power(\{0,1,2\})$ be defined
by $\sigma(\emptyset)=\emptyset$, $\sigma(\{0\})=\{0\}$, $\sigma(\{1\})=\{1\}$ and
$\sigma(\{0,1\})=\{0,1,2\}$. It is an order embedding that cannot be expressed as
the image of any function $h$ on $\{0,1\}$. 
\end{example}
\subsection{Powers of $\N$}
\label{sec:powers-n}

Now we examine order embeddings between powers of $\N$. These are 
the most important for purposes of this series of papers. Let $I$ be an index set. 
We have seen in example~\ref{x-9} that the product poset $(\irri I,\le)$ is in fact
a lattice monoid, and moreover is a complete semilattice satisfying a number of
properties including infinite $+$-distributivity over $\lor$. 
Note that the lattice operations are
given by $x\lor y=\max\{x,y\}$ and $x\land y=\min\{x,y\}$ 
where the $\max$ and $\min$ are taken coordinatewise, 
e.g.~$\max\{f,g\}(i)=\max\{f(i),g(i)\}$ for all $i\in I$. 

\begin{notation}
For $x\in\prod_{i\in I}X_i$ and $y\in\prod_{j\in J}Y_j$ we let $x\bigext y$ denote
the image of $(x,y)$ under the natural association
between $\bigl(\prod_{i\in I}X_i\bigr)\times\bigl(\prod_{j\in J}Y_i\bigr)$ 
and $\prod_{i\in I\djun J}Z_i$ where $Z_i=X_i$ for $i\in I$ and $Z_j=Y_j$ for $j\in J$.
Similarly for functions, 
i.e.~we write $f\bigext g$ for $x\mapsto f(x)\bigext g(x)$. 
\end{notation}

\begin{defn}
\label{d-3}
Let $I$ be an index set.
For a function $g:J\to I$ ($J\subseteq I$), the \emph{projection by $g$} is element
$\pi_g:\irri I\to \irri J$ defined by
\begin{equation*}
  \pi_g(x)=x\circ g.
\end{equation*}
\end{defn}

\begin{prop}
\label{p-29}
$\pi_g$ is a continuous lattice homomorphism. 
It is an epimorphism iff $g$ is an injection, and it is a
monomorphism iff $g$ is onto.
\end{prop}

\begin{notation}
\label{notn:irri}
For $i\in I$, we let $\chi_i\in \irri I$ denote the characteristic function of $i$,
i.e.~$\chi_i(j)$ is $1$ if $j=i$ and $0$ if $j\ne i$. We write $\zero_I\in\irri I$
for element that is constantly equal to $0$, i.e.~$\zero_I(i)=0$ for all $i\in I$.
More generally, for $n\in\N$ write $\mathbf n_I$ for the element constantly equal to
$n$.
\end{notation}

\noindent Clearly, for every $x\in\irri I$,
\begin{equation}
  \label{eq:111}
  x=\spr_{i\in I}x(i)\cdot\chi_i.
\end{equation}

\begin{thm}
\label{l-8}
The order embeddings $\sigma$ from
$(\irri I,\allowbreak \le)$
into $(\irri J,\allowbreak\le)$ with convex range consist of maps of the form 
\begin{equation}
  \label{eq:108}
  \sigma=\pi_g\bigext\zero_{J\setminus \dom(g)}+y
\end{equation}
for some bijection $g:K\to I$ where $K\subseteq J$ and some $y\in\irri J$, and all
such maps are dually continuous lattice embeddings. Thus $\sigma$ is the
sum of a monoid homomorphism for coordinatewise addition and a constant. 
\end{thm}
\begin{proof}
We know that every order embedding of $\irri I$ into $\irri J$ with convex range is
a dually continuous lattice homomorphism by corollary~\ref{o-20}.
Put $y=\sigma(\zero_I)$. 
From corollary~\ref{l-4} and example~\ref{x-2}, we have that each
$\chi_i$ is mapped to $k\cdot \chi_{j_i}+y$ for some $j_i\in J$. Applying
proposition~\ref{p-6} with $\chi_{i}\in(\irri I\setminus\{\zero_I\})^0$, by
convexity we see that $k=1$. Arguing similarly by induction, we obtain
\begin{equation}
  \label{eq:69}
  \sigma(k\cdot\chi_i)=k\cdot \chi_{j_i}+y\espc\text{for all $k\in\N$}.
\end{equation}
Thus $g:K\to I$ defined by $K=\{j_i:i\in I\}$ and
\begin{equation}
  \label{eq:5}
  g(j_i)=i
\end{equation}
is a bijection satisfying
$\sigma(k\cdot\chi_{i})=k\cdot\chi_{j_i}+y
=\pi_g(k\cdot\chi_i)\bigext\zero_{J\setminus K}+y$. 
Thus the continuity of $\sigma$, equation~\eqref{eq:111}, proposition~\ref{p-29} and
infinite $(+,\lor)$-distributivity yield
\begin{equation}
  \label{eq:112}
  \begin{split}
  \sigma(x)&=\sigma\left(\spr_{i\in I}x(i)\cdot\chi_i\right)
  =\spr_{i\in I}\sigma\bigl(x(i)\cdot\chi_i\bigr)
  =\spr_{i\in I}\bigl(x(i)\cdot\chi_{j_i}+y\bigr)\\
  &=\spr_{i\in I}\Bigl(\pi_g\bigl(x(i)\cdot\chi_i\bigr)\bigext\zero_{J\setminus
    K}+y\Bigr)
  =\left(\spr_{i\in I}\pi_g\bigl(x(i)\cdot\chi_i\bigr)\bigext\zero_{J\setminus K}\right)+y\\
  &=\pi_g(x)\bigext\zero_{J\setminus K}+y.
  \end{split}
\end{equation} 
Conversely, it is clear that
whenever $g$ is a bijection and $y\in\irri J$ is arbitrary, $\sigma$ as above
defines an embedding with convex range. 
\end{proof}

\begin{cor}
\label{o-2}
Every map $\sigma:\irri I\to \irri J$ satisfying\tu{: 
\begin{enumerate}[label=(\alph*), ref=\alph*, leftmargin=*, widest=b]
\item\label{item:8} \textit{$x\le y$ iff $\sigma(x)\le\sigma(y)$, for all $x,y\in\irri I$,}
\item \label{item:9}\textit{$\sigma(x)\le z\le\sigma(y)$ implies $z\in\ran(\sigma)$, 
                                        for all $x,y\in\irri I$ and $z\in \irri J$,}
\end{enumerate}}
\noindent is of the form
\begin{equation*}
  \sigma=\pi_g\bigext\zero_{J\setminus \dom(g)}+\sigma(\zero_I)
\end{equation*}
for some bijection $g:K\to I$ where $K\subseteq J$. 
\end{cor}
\begin{proof}
Theorem~\ref{l-8} applies as
conditions~\eqref{item:8} and~\eqref{item:9} translate to $\sigma$ being an
embedding with convex range.
\end{proof}

\begin{cor}
\label{o-13}
Any embedding $\sigma:\irri I\to\irri J$ with convex range has 
$\ran(\sigma)={\irri K}\bigext\zero_{J\setminus K}+y$ for some $K\subseteq J$ 
and $y\in\irri J$. 
\end{cor}
\begin{proof}
By theorem~\ref{l-8} because $\ran(\pi_g)=\irri K$ since $g$ is injective.
\end{proof}

\begin{cor}
\label{o-6}
Every order embedding $\sigma$ of $\irri I$ into a downwards closed subset of $\irri J$ is
of the form
\begin{equation*}
  \sigma=\pi_g\bigext\zero_{J\setminus\dom(g)}
\end{equation*}
for some bijection $g:K\to I$ where $K\subseteq J$, and thus has $\ran(\sigma)={\irri
K}\bigext\zero_{J\setminus K}$. In particular, it is a monoid homomorphism. 
\end{cor}

\begin{cor}
\label{o-3}
Every quasi order isomorphism $\sigma:\irri I\to\irri J$ is of the form $\pi_g$ where
$g:J\to I$ is a bijection. Thus order automorphisms of $\irri I$ are all of the form
$\pi_g$ where $g$ is a permutation of $I$. 
\end{cor}
\subsection{Category algebras}
\label{sec:category-algebras}

To start we review the basic concepts.
Let $X$ be a topological space. We let $\meager(X)$ denote the $\sigma$-ideal of
meager subsets of $X$ (i.e.~countable unions of nowhere dense sets). 
Recall that $X$ is a \emph{Baire space} if it has no nonempty open meager sets.
And recall that $B\subseteq X$ has the \emph{Baire property} 
if it can be approximated by an open set; 
that is there is an open $U$ such that $B\diff U\in\meager(X)$. 
We write $\bp(X)$ for the $\sigma$-algebra of subsets of $X$ with the Baire
property. 
The \emph{category algebra} of $X$ is the quotient algebra $\bp(X)\div\meager(X)$,
i.e.~equivalence classes of $\bp(X)$ modulo $A\diff B\in\meager(X)$ ordered by
$[A]\le[B]$ if $A\setminus B\in\meager(X)$. It is a complete Boolean algebra, with
$[A]\lor[B]=[A\cup B]$, $[A]\land[B]=[A\cap B]$ and $-[A]=[A^\complement]$.
We shall denote the category algebra by $\cat(X)=\bp(X)\div\meager(X)$. 
Recall that a subset $G\subseteq X$ is
\emph{regular open} if it is equal to the interior of its closure; symbolically,
$G=\overbarg G\interior$, i.e.~$\overbarg G$ denotes the topological closure of $G$
while $H\interior$ denotes the interior of $H$. 
We write $\ro(X)$ for the regular open algebra of $X$.  
It is a complete Boolean with $G\lor H=\overline{G\cup H}\interior$, 
$G \land H=G\cap H$ and $-G={\overbar G}{}^\complement$.
A topological space is \emph{semiregular} if its regular open sets form a base 
for its topology.

We shall make use of the following topological facts.

\begin{thm}[Baire Category Theorem]
\label{u-6}
Every locally compact Hausdorff space is a Baire space.
\end{thm}

\begin{prop}
\label{p-9}
If $\U$ is a collection of pairwise disjoint open sets, and $M_U\subseteq U$ is
meager for all $U\in\U$, then $\bigcup_{U\in\U}M_U$ is meager. 
\end{prop}
\begin{proof}
The main observations are that for any open set $U$, $\overbarg U\setminus U$ is
nowhere dense, and that the collection of nowhere dense sets forms an ideal. 
\end{proof}

\begin{prop}
\label{p-1}
Let $X$ be a topological space. Then it has a largest open meager set $U_X$. 
And $G\mapsto[G]$ is a Boolean algebra isomorphism 
between $\ro(X\setminus\overbarg U_X)$ and $\cat(X)$.
Thus if $X$ is a Baire space 
then its category algebra is isomorphic to its regular open algebra. 
\end{prop}
\begin{proof}
Use proposition~\ref{p-9}. 
\end{proof}

\begin{prop}
\label{p-11}
If $X$ is regular then it is semiregular, 
and in particular every nonempty open subset of $X$ 
contains a nonempty regular open subset. 
\end{prop}

\begin{defn}
\label{d-12}
We shall call $Y$ a \emph{Baire compactification} of a topological space~$X$, 
if
\begin{enumerate}[label=(\roman*), ref=\roman*, leftmargin=*, widest=iii]
\item $Y$ is a compact Hausdorff space,
\item $X$ topologically embeds into $Y$, i.e.~there is a map $f:X\to Y$ such that $f$
is a homeomorphism between $X$ and $f[X]$,
\item\label{item:1} $f[X]$ is comeager in $Y$.
\end{enumerate}
\end{defn}

Thus condition~\eqref{item:1} is a strengthening of the usual notion of a
compactification where it is only required that $f[X]$ is dense.

\begin{example}
\label{x-6}
Locally compact Hausdorff spaces have a Baire compactification 
since they in fact have a one point compactification.
\end{example}

\begin{prop}
\label{p-57}
Any topological space with a Baire compactification is a Baire space.
\end{prop}
\begin{proof}
Straight from definition~\ref{d-12}. 
\end{proof}

\begin{lem}
\label{l-1}
Every Polish space \tu(i.e.~separable completely metrizable space\tu) 
has a metric Baire compactification. 
\end{lem}
\begin{proof}
See e.g.~\cite{MR1321597}.
\end{proof}

\begin{example}
\label{x-5}
The irrationals $\irr$ with the product topology, which are homeomorphic to the
irrationals of the real line,  is a Polish space.
\end{example}

\begin{prop}
\label{p-2}
If $Y$ is a Baire compactification of $X$ then $\cat(Y)$ is isomorphic
as a Boolean algebra to $\cat(X)$. 
\end{prop}
\begin{proof}
Identifying $X$ as a comeager subset of $Y$, $[B]\mapsto[B\cap X]$ is an isomorphism
between $\bp(Y)\div\meager(Y)$ and $\bp(X)\div\meager(X)$.
\end{proof}

\begin{thm}
\label{u-3}
Let $X$ be a regular space, and suppose $Y$ has a Baire compactification. 
The order embeddings $\sigma:\cat(X)\to\cat(Y)$ 
with convex range consist of all maps of the form
\begin{equation}
\label{eq:52}
  \sigma([B])=[f[B]]\lor[A]\espc\tu{for all $B\in\bp(X)$},
\end{equation}
for some Baire measurable $f:X\to Y$ and $A\in\bp(Y)$ such that
\tu{
\begin{enumerate}[label=(\alph*), ref=\alph*, leftmargin=*, widest=b]
\item\label{item:15} \textit{$f[B]\in\bp(Y)$ for all $B\in\bp(X)$},
 \item  \label{eq:26} 
\textit{$M\in\meager(X)$ iff $f[M]\in\meager(Y)$ for all $M\subseteq X$},
\item \label{item:11}\textit{$f[B]\cap f[B^\complement]\in\meager(Y)$ for all
    $B\in\bp(X)$},
\item\label{item:10} $f[X]\cap A\in\meager(Y)$.
\end{enumerate}}
\noindent 
Moreover, any such map is a dually continuous lattice embedding. 
Furthermore, $f[X]$ is comeager iff $\sigma$ is a Boolean algebra isomorphism.
\end{thm}
\begin{proof}
We can simplify by assuming that $X$ is a Baire space by proposition~\ref{p-1},
since clearly $\bp(X\setminus \overbarg U_X)\div\meager(X\setminus\overbarg U_X)
\cong\bp(X)\div\meager(X)$. 
Let $\bar\sigma:\cat(X)\to\cat(Y)$ 
be an order embedding with convex range. 
Then $\bar\sigma':\cat(X)\to\cat(Y)$ given by 
$\bar\sigma'(a)=\bar\sigma(a)-\bar\sigma(0_{\cat(X)})$ 
has a downwards closed range, 
and it will suffice to prove that it is of the form
\begin{equation}
  \label{eq:53}
  \bar\sigma'([B])=[f[B]]\espc\text{for all $B\in\bp(X)$}
\end{equation}
for a Baire function $f:X\to Y$ satisfying~\eqref{item:15}--\eqref{item:11},
because then taking $[A]=\bar\sigma(0_{\cat(X)})$, condition~\eqref{item:10} follows
from~\eqref{eq:53} and~\eqref{item:15} and~\eqref{eq:26} since $\bar\sigma$ is 1--1.
By going to a Baire compactification via proposition~\ref{p-2}, 
we can assume that $Y$ is compact and Hausdorff.
Note that $Y$ is a Baire space by the Baire Category Theorem (theorem~\ref{u-6}). 
Hence, by proposition~\ref{p-1}, we can define $\sigma:\ro(X)\to\ro(Y)$ by
$[\sigma(G)]=\bar\sigma'([G])$ and thus obtain
an embedding with downwards closed range.

For each $z\in X$, let
\begin{equation}
  \label{eq:23}
  \G_z=\{G\in\ro(X):G\ni z\},
\end{equation}
and
\begin{equation}
  \label{eq:24}
  K_z=\bigcap_{G\in\G_z}\overline{\sigma(G)}.
\end{equation}
Since $\G_z$ is closed under finite intersections, $K_z\ne\emptyset$ by compactness.
Now choose $f:X\to Y$ so that $f(z)\in K_z$  for all $z\in X$.

\begin{claim}
\label{c-4}
$G\diff f\inv[\sigma(G)]\in\meager(X)$ for all $G\in\ro(X)$.
\end{claim}
\begin{proof}
Let $\V=\{H\in\ro(X):\overline{\sigma(H)}\subseteq \sigma(G)\}$. 
Since $z\in H$ and $\overline{\sigma(H)}\subseteq G$ 
imply $f(z)\in K_z\subseteq \overline{\sigma(H)}\subseteq \sigma(G)$,  
and since $\sigma(H)\subseteq\sigma(G)$ implies $H\subseteq G$ by order reflection, 
\begin{equation}
\label{eq:27}
\bigcup\V\subseteq G\cap f\inv[\sigma(G)].
\end{equation}
On the other hand, letting $\W=\{H\in\ro(X):
\sigma(H)\subseteq \sigma(G)^\complement\}$,  
since $z\in H$ and $\overline{\sigma(H)}\subseteq \sigma(G)^\complement$ 
imply $f(z)\notin \sigma(G)$, and since $\sigma(H)\subseteq \sigma(G)^\complement$
implies $H\cap G=\emptyset$ by order preservation, we obtain
\begin{equation}
  \label{eq:28}
  \bigcup\W\subseteq (G\cup f\inv[\sigma(G)])^\complement.
\end{equation}
Therefore, as $\bigl(G\diff f\inv[\sigma(G)]\bigr)\cap\bigl(\bigcup\V\cup\bigcup\W\bigr)
=\emptyset$, it remains
to show that $\bigcup\V\cup\bigcup\W$ is comeager. But if this were not the case,
then there would exist $H\in\ro(X)^+$ such that
\begin{equation}
\label{eq:8}
H\cap\bigcup\V\cup\bigcup\W=\emptyset. 
\end{equation}
And in particular $H\notin\W$, and thus
$\sigma(H)\cap\sigma(G)\ne\emptyset$. Then since $Y$ is in particular regular,
there exists a nonempty $Z\in\ro(Y)$ with
$\overbar{Z}\subseteq\sigma(H)\cap\sigma(G)$ by proposition~\ref{p-11}. 
Now since the range of $\sigma$ is downwards closed, 
there exists a nonempty $E\in\ro(X)^+$ such that $\sigma(E)=Z$, and hence $E\in\V$. 
However, by order reflection $E\subseteq H$, contradicting~\eqref{eq:8}.
\end{proof}

\begin{claim}
\label{c-7}
$f\inv[H]\in\bp$ for all $H\in\ro(Y)$.
\end{claim}
\begin{proof}
We have
$f\inv[H]=f\inv[H\cap\sigma(X)]\cup f\inv[H\setminus\sigma(X)]$,
and thus $f\inv[H]\diff f\inv[H\cap\sigma(X)]\in\meager$, since
$f\inv[\sigma(X)]$ is comeager by claim~\ref{c-4}. It thus suffices to prove
$f\inv[H\cap\sigma(X)]\in\bp$. 
But by downwards closedness, there exists $G\in\ro(X)$ such that
$\sigma(G)=H\cap\sigma(X)$. Thus the result follows by applying claim~\ref{c-4}
to~$G$. 
\end{proof}

\begin{claim}
\label{c-8}
$f\inv[[\meager]]\subseteq\meager$, 
i.e.~$f\inv[M]\in\meager(X)$ for all $M\in\meager(Y)$.
\end{claim}
\begin{proof}
It suffices to prove that $f\inv[F]$ is a nowhere dense set 
for every nowhere dense $F\subseteq Y$. 
Let $U\subseteq X$ be a given nonempty set. We need to find a nonempty open
$V\subseteq U$ such that $f[V]\cap F=\emptyset$. 
By proposition~\ref{p-11}, there exists $G\subseteq U$ in $\ro(X)^+$. 
Then there exists $E\le\sigma(G)$ in $\ro(Y)^+$ such that 
$\overbar E\cap F=\emptyset$. 
And there exists $H\in\ro(X)^+$ such that $\sigma(H)=E$. 
Moreover, by order reflection, $H\subseteq G\subseteq U$. 
Now $f[H]\subseteq\overline{\sigma(H)}$ and
thus $f[H]\cap F=\emptyset$, as needed.
\end{proof}

\begin{claim}
\label{c-6}
$f[[\meager]]\subseteq\meager$.
\end{claim}
\begin{proof}
Let $F\subseteq X$ be nowhere dense. 
It suffices to  prove that for all $H\in\ro(Y)^+$, 
there exists $Z\subseteq H$ in $\ro(Y)^+$ such that $Z\cap f[F]=\emptyset$. 
Given said $H$, if $H\cap\sigma(X)=\emptyset$ then since
$f[X]\subseteq\overline{\sigma(X)}$, we already have $H\cap f[F]=\emptyset$. 
Otherwise, there exists $G\in\ro(X)^+$ such that $\sigma(G)\subseteq H$. 
Since $G\nsubseteq \overbar F$, 
by regularity  there exist $Z,W\in\ro(X)^+$ such that $Z\subseteq G$, 
$W\supseteq F$ and $W\cap Z=\emptyset$.
Clearly $f[F]\subseteq\overline{\sigma(W)}$, 
and since $\sigma(W)\cap\sigma(Z)=\emptyset$ 
which implies $\overline{\sigma(W)}\cap\sigma(Z)=\emptyset$, 
$f[F]\cap\sigma(Z)=\emptyset$. By order preservation, $\sigma(Z)\subseteq H$,
completing the proof. 
\end{proof}

\begin{claim}
\label{c-1}
$f[G]\diff \sigma(G)\in\meager(Y)$ for all $G\in\ro(X)$.
\end{claim}
\begin{proof}
First we note that $\sigma(G)\setminus f[G]\in\meager$. For supposing to the
contrary, there exists $H\subseteq\sigma(G)$ in $\ro(Y)^+$ such that 
\begin{equation}
\label{eq:49}
H\cap f[G]\in\meager.
\end{equation}
By the assumption on $\ran(\sigma)$, there is an $E$ such that $\sigma(E)=H$, 
and $E\subseteq G$ by order reflection. 
Now $E=(E\cap f\inv[H])\cup(E\setminus f\inv[H])$, 
and thus $E\cap f\inv[H]\notin\meager$ by claim~\ref{c-4}. 
But $E\cap f\inv[H]\subseteq f\inv[H\cap f[G]]\in\meager$ 
by~\eqref{eq:49} and claim~\ref{c-8}, a contradiction.

On the other
hand $f[G]\subseteq\overline{\sigma(G)}$ 
proving $f[G]\setminus\sigma(G)\in\meager$.
\end{proof}

By claims~\ref{c-7} and~\ref{c-8}, $f$ is Baire measurable. 
And by claim~\ref{c-1}, $\bar\sigma'([G])=[\sigma(G)]=[f[G]]$ for all $G\in\ro(X)$. 
Note that condition~\eqref{eq:26} is
equivalent to the conjunction of claims~\ref{c-8} and~\ref{c-6}. 
And now~\eqref{item:15} follows, and moreover 
the expression in~\eqref{eq:53} is well defined, 
i.e.~$[f[B]]$ does not depend on the choice of $B\in[G]$. 
As for condition~\eqref{item:11}, 
given $B\in\bp(X)$ it suffices to take the $G\in\ro(X)$ with $[G]=[B]$, 
and observe that $f[G^\complement]\diff f[\overbarg G{}^\complement]\in\meager(Y)$ 
by~\eqref{eq:26}, and thus
$[f[G]\cap\nolinebreak f[G^\complement]]=[f[G]]\land[f[\overbarg G{}^\complement]]=
[\sigma(G)]\land[\sigma(\overbarg G{}^\complement)]=[\sigma(G)]\land[-\sigma(G)]
=0_{\cat(Y)}$.

Conversely, suppose that $f:X\to Y$ and $A\in\bp(Y)$ satisfy the hypotheses. 
Then $\sigma$ determined by~\eqref{eq:52} makes sense 
and is well defined by~\eqref{item:15} and~\eqref{eq:26}, 
and convexity of the range follows from the Baire measurability of $f$.
It is order preserving because $[B]\le [C]$ implies that 
$f[B]\setminus f[C]\subseteq f[B\setminus C]\in\meager$ by~\eqref{eq:26}, 
which implies that $\sigma([B])-\sigma([C])=0_{\cat(Y)}$. 
And it is order reflecting because $\sigma([B])-\sigma([C])=0_{\cat(Y)}$
implies $f[B]\setminus f[C]\in\meager$ by~\eqref{item:10}, and this implies that
$f[B\setminus C]\in\meager$: For 
\begin{equation}
\label{eq:4}
f[B]\cap f[C]^\complement=f[B]\setminus f[C]
\subseteq f[B\setminus C]\subseteq f[B]\cap f[C^\complement]
\end{equation}
and $(f[B]\cap f[C^\complement])\setminus(f[B]\cap f[C]^\complement)\subseteq
f[C]\cap f[C^\complement]\in\meager$ by~\eqref{item:11}. 
Therefore, $B\setminus C\in\meager$ by~\eqref{eq:26} as required.
\end{proof}

In the context of Boolean algebra homomorphisms, it would be more natural to use
$\bar f\inv$ for some $\bar f:Y\to X$ rather than $f:X\to Y$ since this would allow
for homomorphisms that are not monomorphisms. We reformulate the preceding 
theorem in this manner.

\begin{cor}
\label{o-15}
Let $X$ be a regular space, and let $Y$ have a Baire compactification. Then order
embeddings of $\cat(X)$ into $\cat(Y)$ with convex range consist precisely of maps
of the form
\begin{equation}
  \label{eq:76}
  \sigma([B])=[f\inv[B]]\lor[A]\espc\tu{for all $B\in\bp(X)$},
\end{equation}
for some $Z\in\bp(Y)$ and some Baire measurable injection $f:Z\to X$ such that
\tu{
\begin{enumerate}[leftmargin=*, label=(\alph*), ref=\alph*, widest=b]
\item\label{item:2} \textit{$f[C]\in\bp(X)$ for all $C\in\bp(Z)$,}
\item\label{item:7} \textit{$M\in\meager(X)$ iff $f\inv[M]\in\meager(Y)$ for all
    $M\subseteq X$,}
\item\label{item:18} \textit{$Z\cap A\in\meager(Y)$.}
\end{enumerate}}
\noindent 
Furthermore, $Z$ is comeager iff $\sigma$ is a Boolean algebra isomorphism. 
\end{cor}
\begin{proof}
Supposing $\sigma$ is an order embedding, 
let $\bar f\in\baire(X,Y)$ and $A\in\bp(Y)$ satisfy the conclusion of
theorem~\ref{u-3}. Put $Z=\ran(\bar f)$. 
Choose a right inverse $f:Z\to X$ of $\bar f$. The fact that it is a right inverse entails
\begin{equation}
  \label{eq:73}
  f\inv[B]\subseteq\bar f[B]\espc\text{for all $B\subseteq X$}.
\end{equation}
And plugging in $B^\complement$ for $B$ yields
$\bar f[B^\complement]^\complement\subseteq f\inv[B]$ for all $B\subseteq X$.
Therefore
\begin{equation}
  \label{eq:75}
  \bar f[B]\setminus f\inv[B]\subseteq\bar f[B]\cap\bar
  f[B^\complement]\in\meager(Y)
  \espc\text{for all $B\in\bp(X)$}
\end{equation}
by condition~\eqref{item:11} for $\bar f$. It now follows from the conclusion
of theorem~\ref{u-3} that $\sigma$ satisfies equation~\eqref{eq:76} and that $f$ is
a Baire measurable function satisfying~\eqref{item:7} and~\eqref{item:18}.
Condition~\eqref{item:2} can be verified similarly.

The converse is completely straightforward.
\end{proof}

It seems likely to us that the
preceding corollary is known for the special case of Boolean algebra isomorphisms
(i.e.~$Z$ comeager), although we do not know of it in the literature, except in the
specific case of some explicit  Polish spaces $X$ and $Y$. 

\begin{cor}
\label{o-14}
Let $X$ be regular and $Y$ have a Baire compactification.
The range of any order embedding $\sigma:\cat(X)\to\cat(Y)$ with convex range is the
interval Boolean algebra $[\sigma(0_{\cat(X)}),\allowbreak\sigma(1_{\cat(X)})]$, 
and $i\inv\circ\sigma$ is a Boolean algebra isomorphism between $\cat(X)$ 
and the interval Boolean algebra, where $i$ is the inclusion map.
\end{cor}

In the second countable case we get a sharper description of $f$. 

\begin{cor}
\label{u-8}
Let $X$ be a second countable Hausdorff space, and let $Y$ be a Polish space.
The order embeddings $\sigma:\cat(X)\to\cat(Y)$ 
with convex range consist of all of the maps of the form
\begin{equation}
\label{eq:54}
  \sigma([B])=[f[B\cap C]]\lor [A]\espc\tu{for all $B\in\bp(X)$},
\end{equation}
for some comeager $C\subseteq X$, $f:C\to Y$ and $A\subseteq Y$ such that\tu{
  \begin{enumerate}[label=(\alph*), ref=\alph*, leftmargin=*, widest=b]
  \item \textit{$f$ is a topological embedding,}
  \item \label{item:16} $f[C]\in\bp(Y)$,
  \item\label{item:17} $f[C]\cap A\in\meager(Y)$.
  \end{enumerate}}
\end{cor}
\begin{proof}
Let $\sigma$ be an embedding with convex range.
Note that $X$ is regular, because $X$ is metrizable 
(by Urysohn's Metrization Theorem), and note that $Y$ has a Baire compactification
by lemma~\ref{l-1}.
By theorem~\ref{u-3}, there is a Baire measurable $f:X\to Y$ and $A\in\bp(Y)$ 
such that 
\begin{enumeq}
\item \label{eq:51} $\sigma([B])=[f[B]]\lor[A]$ for all $B\in\bp(X)$,
\item\label{item:12}$M\in\meager(X)$ iff $f[M]\in\meager(Y)$ for all $M\subseteq X$,
\item \label{item:13} $f[B]\cap f[B^\complement]\in\meager(Y)$ for all $B\in\bp(X)$,
\item \label{item:14} $f[X]\cap A\in\meager[Y]$.
\end{enumeq}
Since $Y$ is a Baire space by proposition~\ref{p-57}, we can
define $\bar\sigma:\ro(X)\to\ro(Y)$ by letting $\bar\sigma(G)$ be the unique member
of $\ro(Y)$ such that $[\bar\sigma(G)]=\sigma([G])-[A]$. Note that $\bar\sigma$ is
an embedding with downwards closed range because $[A]=\sigma(0_{\meager(X)})$
by~\eqref{eq:51}. 

Then equations~\eqref{eq:51} and~\eqref{item:14} together imply 
that $f[G]\diff\bar\sigma(G)\in\meager(Y)$ for all $G\in\ro(X)$, 
and equation~\eqref{item:13} implies that $f[G]\cap f[G^\complement]\in\meager(Y)$
for all $G\in\ro(X)$. 
As $f\inv[f[G]\diff\bar\sigma(G)]\in\meager(X)$ and
$f\inv[f[G]\cap f[G^\complement]]\in\meager(X)$ for all $G$ by~\eqref{item:12}, 
and as $X$ is second countable, we can find a comeager $C\subseteq X$ such that
\begin{align}
  \label{eq:7}
  f[C]\cap (f[G]\diff\bar\sigma(G))&=\emptyset,\\
  \label{eq:9}
  f[C]\cap f[G]\cap f[G^\complement]&=\emptyset,
\end{align}
for all $G\in\ro(X)$. Therefore,
\begin{equation}
  \label{eq:56}
  f[G\cap C]=\bar\sigma(G)\cap f[C]\espc\text{for all $G\in\ro(X)$},
\end{equation}
because $f[G\cap C]=f[G]\cap f[C]$ by~\eqref{eq:9}, and $f[G]\cap
f[C]=\bar\sigma(G)\cap f[C]$ by~\eqref{eq:7}.
By the semiregularity of $X$ this implies that $f\restriction C$ 
maps all opens sets to relative open
subsets of its range. And since $X$ is Hausdorff, every two distinct points in $C$
have disjoint regular open neighbourhoods, and thus their images under
$f\restriction C$ are disjoint
by~\eqref{eq:56} since $\bar\sigma$ is an embedding; hence, $f\restriction C$ is 1--1.
And by shrinking $C$, since $Y$ is in particular second countable, we 
can moreover arrange that $f\restriction C$ is continuous (by Baire's Theorem). 
Therefore $f\restriction C$ is a topological embedding, while~\eqref{item:16}
and~\eqref{item:17} are consequences of~\eqref{eq:51} and~\eqref{item:14},
respectively. 

Conversely, it is clear that any $C\subseteq X$, $f:C\to Y$ and $A$ 
satisfying the hypotheses
define an order embedding $\sigma$ with convex range via~\eqref{eq:54}.
\end{proof}
\subsection{Baire functions into powers of $\N$ modulo almost always equality}
\label{sec:baire-functions}

The results in this section will play a crucial role in the next paper~\cite{irrationals} of
this series. The main results are theorem~\ref{l-2} where a number of properties of
the monoid $\baireaa{X,\irri I}$ are established; theorem~\ref{u-2} where embeddings
of $\baireaa{X,\irri I}$ into $\baireaa{Y,\irri J}$ with convex range are
characterized for arbitrary topological spaces $X$ and $Y$; and theorem~\ref{o-8}
where these embeddings, for $X$ regular and $Y$ having a Baire compactification, are
described in terms of generalized projections (definition~\ref{d-13}).

Before preceding, we recall that a \emph{Baire measurable function} is
of course a function $f:X\to Y$ between two topological spaces such that
$f\inv(U)\in\bp(X)$ for every open $U\subseteq Y$. We often just call them
\emph{Baire functions} (not to be confused with the different concept of Baire class
$n$ functions) and we write $\baire(X,Y)$ for the family of Baire functions from $X$
into $Y$. 
And we write $\cat(X,Y)$ for the family of continuous functions $X$ into $Y$. 

Given some relation $R$ on a topological space $S$, and given a topological space
$X$, we would like to use $R$ to induce a relation on $\baire(X,S)$ in the sense of
``almost always'' in the space $X$. However, we
shall see that this is not as straightforward as it may seem. 

It will be constructive to express a given relation as a conjunction. Suppose then
that $\R$ is a family of binary relations on a fixed space $S$. 
We define a relation $\R_\aa$ on $\baire(X,S)$ by
\begin{equation}
  \label{eq:18}
  f\mathrel{\R_\aa} g\If \ifm_{R\in\R} \ulc f(z) \rel g(z)
  \text{ for almost all }z\in X\urc,
\end{equation}
where ``for almost all'' is interpreted as comeagerly
many. Thus we are generalizing
the relation $\ifm\R$ from $S$ to $\baire(X,S)$. However, the resulting
relation on the function space may depend on the representation when the family $\R$
is uncountable, as shown in example~\ref{x-8} below. Also note that the definition
in~\eqref{eq:18} can be generalized from binary relations to arbitrary $n$-place
relations.

\begin{prop}
\label{p-31}
If each $R\in\R$ is a quasi ordering of $S$, 
then $\R_\aa$ is a quasi ordering of $\baire(X,S)$.
\end{prop}

\begin{example}
\label{x-7}
We consider the relation $\le$ on $\irri I$. 
For each $i\in I$, 
let $R_i$ be the relation satisfying $x\rel_i y$  iff $x(i)\le y(i)$. 
Then ${\le}={\ifm_{i\in I}R_i}$. Putting $\R=\{R_i:i\in I\}$ we obtain the following
relation on $\baire(X,\irri I)$:
\begin{equation}
  \label{eq:20}
  f\mathrel{\R_\aa} g\Iff f(z)(i)\le g(z)(i)
  \text{ for almost all $z\in X$}\espc\text{for all $i\in I$}.
\end{equation}
Let us denote $\R_\aa$ from equation~\eqref{eq:20} as $\leaa$. 
We repeat this construction for the equality relation. 
Thus we let $\S=\{S_i:i\in I\}$ where $x\mathrel{S_i} y$ iff $x(i)=y(i)$, and let
$\eqaa$ denote the relation $\S_\aa$. These are both quasi orderings of the family
of Baire functions by proposition~\ref{p-31}. 
\end{example}

These are the `correct' representations of $\le$ and $=$, 
because of the following connection
with set theoretic forcing. 
We shall not attempt to explain the forcing concepts
here (see e.g.~\cite{MR597342}). Noting that $\eqaa$ is an equivalence relation on
$\baire(X,\irri I)$, we observe that there are Borel representatives of the Baire
measurable functions.

\begin{lem}
\label{l-10}
Let $X$ be a topological space. Then $\baireaa{\irri I}$ is equal to $\borelaa{\irri I}$. 
\end{lem}
\begin{proof}
By our representation $\S$ of equality, it suffices to observe that every Baire
function $f:X\to\N$ has a Borel function $g_f:X\to\N$ with $[f]\eqaa[g_f]$. 
Then given $[F]\in\baireaa{\irri I}$ we can define $G:X\to\irri I$ by
$G(z)(i)=g_{F(\cdot)(i)}(z)$ for all $z\in X$ and $i\in I$, which yields a Borel
function $G$ with $[G]\eqaa[F]$. 
\end{proof}

Note too that:

\begin{prop}
\label{p-30}
If $f\in\baire(X,\irri I)$ and $g\eqaa f$ then $g\in\baire(X,\irri I)$. 
\end{prop}

Given $f\in\baire(X,\irri I)$, we let $f'\in\borel(X,\irri I)$ satisfy $[f']\eqaa[f]$. 
The advantage of Borel sets and Borel functions is that they have an explicit
description that allows them to be interpreted in new larger ``universes'' extending
the present one (e.g.~forcing extensions), that may contain e.g.~new real numbers.
Now we make three assumptions on the topological space $X$. 
Firstly, we assume that it is given by an explicit definition, 
by which we of course mean a formula of set theory, possibly with a parameter. 
This is true of any set $X$, but the point is that we might
consider a formula $\varphi(x)$ defining the real line $\reals$ for example, so that
when we go to a forcing extension the same formula gives new reals. 
Secondly, we assume the \emph{defined} space $X$ is Hausdorff, 
i.e.~in any model extending our
present universe and satisfying enough of the axioms of set theory, the space $X$
defined by the formula $\varphi(x)$ is Hausdorff. 
Then for any sufficiently \emph{generic} filter $G$ on the poset
$\bigl(\borel(X)\setminus\meager(X),\subseteq\bigr)$, which produces the same
forcing extension as $\cat(X)$,
$\bigcap G$ has at most one element when the intersection is interpreted in the
forcing extension by $\cat(X)$. Thirdly, 
we require that the defined space satisfies some topological `completeness' so that
$\bigcap G\ne\nobreak\emptyset$; 
this is where we want to use a definition of a topological
space rather than a fixed space, because a fixed space may lose its completeness
property in a larger universe. We will not explicitly formulate this notion of
completeness, but we note that complete metrizability is sufficient although there
is no need to restrict ourselves to metric spaces. 
Under these three assumptions, we let $\dot g$ denote the single element of $\bigcap
G$ in the forcing extension by $\cat(X)$. Back to the function $f$, we write $\dot
f$ for $f'(\dot g)$ where the Borel function $f'$ is being interpreted in the
forcing extension by $\cat(X)$. 

Now the ``forcing theorem'' applied to the forcing notion $\cat(X)$ and 
elements of~$\irri I$ states that 
for every formula $\varphi(x_0,\dots,x_{n-1})$ of the language of set theory, 
for all $f_0,\dots,f_{n-1}\in\baire(X,\irri I)$,
\begin{multline}
  \label{eq:114}
  \cat(X)\forces\varphi(\dot f_0,\dots,\dot f_{n-1})\\\IFF 
  \varphi(f_0(z),\dots,f_{n-1}(z))\text{ for almost all $z\in X$}.
\end{multline}
In particular,  for all $f,g\in\baire(X,\irri I)$,
\begin{align}
  \label{eq:21}
  \cat(X)&\forces\dot f=\dot g\IFF f\eqaa g, \\
  \label{eq:22}
  \cat(X)&\forces\dot f\le\dot g\IFF f\leaa g.
\end{align}
In set theoretic terminology, every member of $\baire(X,\irri I)$ `names' a member
of $\irri I$ lying in the forcing extension by $\cat(X)$, and visa versa. Thus
equations~\eqref{eq:21} and~\eqref{eq:22} state that $f$ and $g$ name the same
object iff $f\eqaa g$, and $f$ names a member of $\irri I$ below the member named by
$g$ iff $f\leaa g$. This correspondence has been known for a long time; for example,
the measure theoretic analogue for measure algebras (discussed in
\Section\ref{sec:l0mu-irri-idiveqae}) appears in 1967 (\cite{MR0218233}).
Indeed $\cat(X)$ is (isomorphic as a forcing notion to) 
the original forcing notion used by Cohen
(\cite{MR0157890},\cite{MR0159745}) to prove that the Continuum Hypothesis is
independent of the usual axioms of mathematics.

On the other hand, we consider the `obvious' representation of equality to be
`incorrect'. Considering the singleton $\{=\}$, we obtain following generalization
of equality to the space $\baire(X,\irri I)$ of Baire functions:
\begin{equation}
  \label{eq:30}
  f\mathrel{\{=\}_\aa} g\Iff f(z)= g(z)\espc\text{for almost all $z\in X$.}
\end{equation}
If the index set $I$ is countable then $\{=\}_\aa$ is equivalent to $\eqaa$ because
$\meager(X)$ is a $\sigma$-ideal. However, if the index set is of larger cardinality
then these relations may differ, as in the following example~\ref{x-8}. Also
lemma~\ref{l-10} may fail with $\{=\}_\aa$. 

\begin{example}
\label{x-8}
Define $f\in\baire(\reals,\irri \reals)$ by
\begin{equation}
  \label{eq:31}
  f(z)(x)=
    \begin{cases}
      1,&\text{if $x=z$,}\\
      0,&\text{if $x\ne z$.}
    \end{cases}
\end{equation}
Then $f\eqaa \zero_\reals$; 
however, it is not the case that $f\mathrel{\{=\}_\aa}\zero_\reals$,
because $f(z)\ne\zero_\reals$ for all $z\in\reals$. 
\end{example}

Perhaps the most basic aspect of $\eqaa$ is that it gives a product order.

\begin{prop}
\label{p-26}
$[f]\mapsto([f(\cdot)(i)]:i\in I)$ defines an isomorphism 
between\linebreak $\baireaa{X,\irri I}$ and $(\baireaa{X,\N})^I$.
\end{prop}

Observe that for all $f,g\in\baire(X,\irri I)$, $f\leaa g$ and $g\leaa f$ iff
$f\eqaa g$, and thus the poset $(\baireaa{X,\irri I},\leaa)$ is the antisymmetric
quotient of $(\baire(X,\irri I),\allowbreak\leaa\nobreak)$ 
(cf.~\Section\ref{sec:homom-with-regul}). 
Also $(\baire(X,\irri I),\le)$, where $f\le g$ means $f(z)\le g(z)$ 
for all $z\in X$, is clearly a lattice where $f\lor g=\max\{f,g\}$ and
$f\land g=\min\{f,g\}$ are taken pointwise:
\begin{equation}
  \label{eq:10}
  \max\{f,g\}(z)=\max\{f(z),g(z)\}\And\min\{f,g\}(z)=\min\{f(z),g(z)\}
\end{equation}
for all $z\in X$ (cf.~\Section\ref{sec:powers-n}). 
Furthermore, $(\baire(X,\irri I),+)$, where addition is taken pointwise, is clearly a
commutative monoid (cf.~example~\ref{x-9}). In fact, it is a lattice monoid where
the monoid ordering $\le$ is the above lattice order. 
And it is easily verified that $\eqaa$ is a congruence for this monoid
(cf.~\Section\ref{sec:terminology}), and
thus the quotient $\baireaa{X,\irri I}$ is a quotient monoid with 
the monoid operation 
\begin{equation}
  \label{eq:113}
  [f]+[g]=[f+g].
\end{equation}
We can also generalize the Baer--Specker group to the quotient group
$\baireaa{X,\intgr^I}$. Then $\baireaa{X,\irri I}$ embeds as a monoid in
$\baireaa{X,\intgr^I}$, and in particular it is a cancellative monoid. 

We use forcing to prove that the quotient is a complete semilattice satisfying all
of the distributivity laws of $\irri I$.

\begin{thm}
\label{l-2}
Let $X$ be a definable Hausdorff and \tu`complete\tu' topological space 
\tu(e.g.~$X=\reals$\tu).  
Then $(\baireaa{X,\irri I},+)$ is a commutative cancellative complete
semilattice monoid, that is \jid, \mid, infinitely $+$-distributive over $\lor$
and~$\land$ and flat-complete. 
The monoid order is $\leaa$, and the lattice operations are given by
\begin{equation}
  \label{eq:12}
  [f]\lor[g]=[\max\{f,g\}]\And[f]\land[g]=[\min\{f,g\}].
\end{equation}
And we have the forcing expressions
\begin{alignat}{4}
  \label{eq:40}
  \cat(X)&\forces\dot h&&=\max\{\dot f,\dot g\}&&\IFF [h]&&=[f]\lor[g],\\
  \label{eq:41}
  \cat(X)&\forces\dot h&&=\min\{\dot f,\dot g\}&&\IFF [h]&&=[f]\land[g],\\
  \label{eq:115} 
  \cat(X)&\forces\dot h&&=\dot f+\dot g&&\IFF [h]&&=[f]+[g].
\intertext{Furthermore, equations~\eqref{eq:40} and~\eqref{eq:41} generalize to 
families $\F\subseteq\allowbreak\baire(X,\irri I)$ as}
  \label{eq:128}
  \cat(X)&\forces\dot h&&=\spr_{f\in\F}\dot f&&\IFF [h]&&=\spr_{f\in\F}[f],\\
  \label{eq:129}
  \cat(X)&\forces\dot h&&=\ifm_{f\in\F}\dot f&&\IFF [h]&&=\ifm_{f\in\F}[f],
\end{alignat}
respectively.
\end{thm}
\begin{proof}
Equation~\eqref{eq:115} follows from the forcing theorem~\eqref{eq:114} with
$\varphi(x_0,x_1,x_2)=\ulc x_0=x_1+x_2\urc$. First we verify that $\leaa$ coincides
with the monoid order. 
Recall that the monoid ordering of $\irri I$ coincides with the usual ordering
(example~\ref{x-9}). Thus using~\eqref{eq:22} and~\eqref{eq:115},
\begin{equation}
\label{eq:116}
\begin{split}
f\leaa h&\IFF \cat(X)\forces \dot f\le\dot h \\
&\IFF \cat(X)\forces \exists g\in\irri I\spc\dot h=\dot f+\dot g\\
&\IFF \exists g\in\baire(X,\irri I)\spc \cat(X)\forces\dot h=\dot f+\dot g\\
&\IFF \exists g\spc [h]=[f]+[g]\\
&\IFF [f]\len{(\baireaa{X,\irri I},+)} [h].
\end{split}
\end{equation}

Next we prove equation~\eqref{eq:40}. 
Recall that $\max$ and $\min$ are the lattice operations on $(\irri I,\le)$. 
Suppose $\cat(X)\forces\dot h=\max\{\dot f,\dot g\}$. 
Then $[f],[g]\leaa[h]$ by~\eqref{eq:22}, and if $[h']\geaa[f],[g]$ then
$\cat(X)\forces\dot h'\ge\dot h$ by~\eqref{eq:22}, and thus $[h]\leaa[h']$ again
by~\eqref{eq:22}, proving $[h]=[f]\lor[g]$. Conversely, suppose $[h]=[f]\lor[g]$.
Then $\cat(X)\forces\dot h\ge\max\{\dot f,\dot g\}$ by~\eqref{eq:22}. 
And if $\cat(X)\forces\dot h'\ge\max\{\dot f,\dot g\}$, 
then $\cat(X)\forces\dot h\le\dot h'$ because $[h]\leaa[h']$, 
completing the
proof that $\cat(X)\forces\dot h=\max\{\dot f,\dot g\}$. The same proof generalizes
to arbitrary $\F$ as in equation~\eqref{eq:128}.
Equations~\eqref{eq:41} and~\eqref{eq:129} hold dually. 
And the lattice operations~\eqref{eq:12} follow from~\eqref{eq:40} and~\eqref{eq:41}
and the forcing theorem~\eqref{eq:114}, 
e.g.~with $\varphi(x_0,x_1,x_2)=\ulc x_0=\max\{x_1,x_2\}\urc$. 

Now we prove that $\baireaa{X,\irri I}$ is a complete semilattice.
Suppose $\F\subseteq\baireaa{X,\irri I}$ is a bounded subfamily. 
Then letting $\F'\subseteq\baire(X,\irri I)$ be a selection of representatives, 
$\cat(X)\forces\ulc\dot\F'$ is bounded$\urc$ by~\eqref{eq:22}. 
Then since $\irri I$ is a complete semilattice, 
there exists $g\in\baire(X,\irri I)$ such that
\begin{equation}
  \label{eq:42}
  \cat(X)\forces\dot g=\spr\dot\F'.
\end{equation}
Thus $[g]=\spr\F$ follows from~\eqref{eq:22}.

To prove the \jid, suppose $f\in\baire(X,\irri I)$, 
$\F\subseteq \baireaa{X,\irri I}$ and $\spr\F$ exists. 
Letting $\F'\subseteq\baire(X,\irri I)$ be a selection of representatives, 
$\cat(X)\forces\spr\dot\F'$ exists. And then since $\irri I$ satisfies the \jid, 
$\cat(X)\forces\dot f\land \spr\dot\F'=\spr(\dot f\land\dot\F')$. 
Therefore, $f\land\spr\F=\spr(f\land\F)$. 

The \mid\ property holds dually.  And infinite $+$-distributivity over $\lor$ and
$\land$ is a consequence of theorem~\ref{l-7} (although it 
can be transferred from $\irri I$ similarly to the \jid). 
Flat-completeness transfers from $\irri I$ similarly, as $\irri I$ is 
flat-complete (example~\ref{x-17}). 
\end{proof}

\begin{remark}
\label{r-6}
Theorem~\ref{l-2} can be proved in a routine matter without using forcing, and
indeed a drawback of the forcing method is the need for the definable topological
requirements on $X$, i.e.~the theorem is true for any topological space $X$. 
However, the reason that it is possible to give a straightforward direct proof of
e.g.~completeness is that
the supremum of a bounded family has a simple definition in terms of the family. 
Indeed there are very similar situations, involving for example the eventual
dominance order $\lefnt$ on $\irri I$ (see~\Section\ref{sec:irrfin-lefnt}) 
rather than the product order, where the supremum 
has no reasonably simple definition and the forcing translation becomes
essential (see e.g.~\cite{MR1621933}). 
\end{remark}

\begin{remark}
\label{r-1}
If the quotient over $\{=\}_\aa$ is taken instead, then the quotient lattice
over the relation $\{\le\}_\aa$ is not in general a complete semilattice.
\end{remark}

\begin{remark}
\label{r-10}
Note that equation~\eqref{eq:12} does not generalize from two functions to arbitrary
families of functions. For a counterexample, consider the family $\{f_z:z\in X\}$
where $f_z(z)=\one_I$ and $f_z(y)=\zero_I$ for all $y\ne z$. 
\end{remark}

\begin{notation}
\label{notn:basic}
Let $X$ be some fixed topological space, with largest open meager set $U_X$
(cf.~proposition~\ref{p-1}). 
For $a\in\cat(X)$, we let $\bar a$ denote the regular open set
$B\in\ro(X\setminus\overbarg U_X)$ such that $a=[B]$. 
Let $I$ be some fixed index set.
For $a\in\cat(X)$ and $x\in\irri I$ with $x\ne\zero_I$, 
we let $\basel a x$ denote the member
$[f]$ of $\baireaa{X,\irri I}$ where $f:X\to\irri I$ is given by
\begin{equation}
\label{eq:66}
  f(z)=
\begin{cases}
x&\text{if $z\in \bar a$,}\\
\zero_I&\text{if $z\notin\bar a$.}
\end{cases}
\end{equation}
Let $\<a,x\>^+$ denote $x$. 
For each $i\in I$, we let $\pi_i:(\irri I)^X\to \irri X$ denote the pointwise projection
\begin{equation}
  \pi_i(f)(z)=f(z)(i)\espc\text{for all $z\in X$}.
\end{equation}

Now assume that $S$ is some fixed topological space.
Suppose that $\varphi(x_0,\dots,x_{n-1},\allowbreak y_0,\dots,y_{m-1})$ is a formula  
and $a_0,\dots,a_{m-1}$ are some parameters such that 
$\bigl\{z\in\nobreak X:\varphi\bigl(f_0(z),\dots,f_{n-1}(z),a_1,\dots,a_{m-1}\bigr)\bigr\}\in\bp(X)$
whenever $f_0,\dots,f_{n-1}$ are Baire measurable functions from $X$ into $S$.
Then we denote the equivalence class 
\begin{equation}
  \label{eq:117}
  [z\in X:\varphi(f_0(z),\dots,f_{n-1}(z),a_1,\dots,a_{m-1})]\in\cat(X).
\end{equation}
by $\lbrak\varphi(\dot f_0,\dots,\dot f_{n-1},a_1,\dots,a_{m-1})\rbrak$. 
In forcing terminology, this is the Boolean truth value of the sentence
$\varphi(\dot f_0,\dots,\dot f_{n-1},a_0,\dots,a_{m-1})$. 
The following examples illustrate the usefulness of this notation.

\begin{example}
  \label{x-4}
Suppose $f:X\to\irri I$ is Baire measurable. Then for all $i\in I$ and~$n\in\N$,
\begin{equation}
  \label{eq:118}
  \|\dot f(i)=n\|=[z\in X:f(z)(i)=n].
\end{equation}
Suppose $K:X\to\power(Y)$ is a Baire function where $Y$ is some fixed set. Then for
all~$y\in Y$,
\begin{equation}
  \label{eq:119}
  \|y\in\dot K\|=[z\in X:y\in K(z)].
\end{equation}
\end{example}

For $[f]\in\baireaa{X,\irri I}$, we denote
\begin{equation}
  \pos([f])=\|\dot f\ne\zero_I\|.
\end{equation}
We can also write $\pi_i([f])$ for $[\pi_i(f)]$ when dealing with
equivalence classes modulo~$\eqaa$. Note that
\begin{equation}
  \label{eq:65}
  \pos(\pi_i([f]))=\pos(\basel{1_{\cat(X)}}{\chi_i}\land[f]).
\end{equation}
\end{notation}

Also note that
\begin{equation}
  \label{eq:80}
  [f]=\spr_{i\in I}\spr_{n=0}^\infty n\cdot\<\| \dot f(i)=n\|,\chi_i\>,
\end{equation}
because one obtains $[f]=\spr_{i\in I}\pi_i([f])\cdot\chi_i$ from equation~\eqref{eq:128} and
the fact that by equation~\eqref{eq:111},
\begin{equation}
  \label{eq:120}
  \cat(X)\forces\dot f=\spr_{i\in I}\dot f(i)\cdot\chi_i=\spr_{i\in I}\dot\pi_i(f)\cdot\chi_i
\end{equation}
(or alternatively use proposition~\ref{p-26});
and since clearly $[\pi_i(f)\cdot\chi_i]=\spr_{n=0}^\infty n\cdot\<\|\dot f(i)=\nobreak
n\|,\allowbreak \chi_i\>$. 

\begin{defn}
We call members of $\baireaa{X,\irri I}$ of the form $\<a,n\cdot\chi_i\>$
($a\in\cat(X)$, $n\in\N$, $i\in I$) \emph{basic elements}.
\end{defn}

\begin{prop}
\label{p-35}
The basic elements form a basis for $(\baireaa{X,\irri I},\leaa)$.
\end{prop}
\begin{proof}
By equation~\eqref{eq:80}, 
because $n\cdot\<\|\dot f(i)=n\|,\chi_i\>\land 
n^*\cdot\<\|\dot f(i^*)=n^*\|,\chi_{i^*}\>=[\zero]$ 
whenever either $i\ne i^*$ or $n\ne n^*$.
\end{proof}

\begin{prop}
\label{p-25}
$(\baireaa{X,\irri I})_{\<a,\chi_i\>}\cong\cat(X)_a$ 
\tu(cf.~notation~\tu{\ref{notn:down}}\tu).
\end{prop}
\begin{proof}
$\<b,\chi_i\>\mapsto b$ is an isomorphism.
\end{proof}

We generalize notation~\ref{notn:irri}.

\begin{notation}
\label{notn:irrigen}
For a function $K:Y\to\power(J)$, we let $\zero_{K-}$ denote the element of
$\cts(\irri{K-})$ given by $\zero_{K-}(z)=\zero_{K(z)}$ for all $z\in Y$. 
For a fixed index set $I\subseteq J$, when $Y$ is implicitly understood,
we write $\zero_{I-}$ for the function in $(\irri I)^Y$ constantly equal to
$\zero_I$; in this case, we may write $\zero-$ when the index set is also understood.
\end{notation}

\begin{thm}
\label{u-2}
If $\sigma$ is an order embedding of $\baireaa{X,\irri I}$ into
$\baireaa{Y,\irri J}$ with convex range then $\sigma$ is a dually continuous
lattice homomorphism determined by 
$\sigma(\basel a{\chi_i})$ \tu($a\in\cat(X)$, $i\in I$\tu). 
Indeed it is of the form
\begin{equation}
\label{eq:14}
   \sigma([f])=\spr_{i\in I}\spr_{n=0}^\infty n
   \cdot \sigma(\basel{\| \dot f(i)=n\|}{\chi_i})+[h]
\end{equation}
for some $[h]\in\baireaa{Y,\irri J}$. And for each $i\in I$, 
$\sigma(\basel a {\chi_i})$ is determined by $a^i_j\in\cat(X)$ \tu{($j\in J$)} where 
\begin{align}
\label{eq:35}
a^i_j\land a^i_{j^*}&=0_{\cat(X)}\espc\tu{for all $j\ne j^*$,} \\
\label{eq:36}
\spr_{j\in J}a^i_j&=1_{\cat(X)},
\end{align}
and isomorphisms
$\tau^i_j:\cat(X)_{a^i_j}\to\cat(Y)_{\tau(a^i_j)}$ \tu($j\in J$\tu) via\tu:
\begin{equation}
  \label{eq:13}
  \sigma(\basel a {\chi_i})=\spr_{j\in J}\basel{\tau^i_j(a\land a^i_j)}{\chi_j}+[h],
\end{equation}
subject to the constraint
\begin{equation}
  \label{eq:15}
  \tau^i_j(a^i_j)\land\tau^{i^*}_j( a^{i^*}_j)=0_{\cat(Y)}\espc\tu{for all $i\ne i^*$}.
\end{equation}
Conversely, any such $\tau^i_j$ \tu($i\in I$, $j\in J$\tu) and $h$ determine an
order embedding $\sigma$ defined by equations~\eqref{eq:14} and~\eqref{eq:13}. 
\end{thm}
\begin{proof}
Note that we can assume that $\sigma([\zero_{I-}])=[\zero_{J-}]$ 
by considering the homomorphism $\sigma'=\sigma-\sigma([\zero_{I-}])$ instead.
That $\sigma$ is a dually continuous lattice homomorphism is by corollary~\ref{o-20}.

\begin{claim}
\label{c-2}
If $[g]$ is a basic element, say  $\basel a {n\cdot\chi_i}$ for some $n$, $i$ and $a$,
then $\sigma([g])\leaa\mathbf n$.
\end{claim}
\begin{proof}
The proof is by induction on $n=0,1,2,\dots$. 
For $n=1,2,\dots$, suppose $[g]=\basel a {n\cdot\chi_i}$ for some $i$ and $a$, and 
by theorem~\ref{l-2} and remark~\ref{r-6} 
we can let $[h_n]$ be the supremum of all $[f]\leaa[g]$ such
that $\sigma([f])\leaa\mathbf n$.
By continuity, $\sigma([h_n])\leaa\mathbf n$.
Assuming the claim fails for $g$, $[h_n]\laa[g]$, 
and thus there exists a $b_n\le a$ such that 
\begin{equation}
  \label{eq:6}
  [h_n]\land\basel {b_n}{n\cdot\chi_i}=\basel {b_n}{k\cdot\chi_i}
\end{equation}
for some $k=0,\dots,n-1$. Then $\sigma(\basel{b_n}{n\cdot\chi_i})$ is not almost always
at most $\mathbf n$, and thus there is a 
basic $[p]\leaa\sigma(\basel{b_n}{n\cdot\chi_i})$ 
such that $p^+=m\cdot\chi_j$ for some $m> n$ and some $j$. 
But then there exists a basic $[q]\leaa[p]$ where $q^+=n\cdot\chi_j$. 
Letting $[f]$ be the element mapped to $[q]$, $[f]\leaa\basel{b_n}{n\cdot\chi_i}$
by order reflection. However, $[f]\leaa[h_n]$ by the definition of $h_n$, and thus
$[f]\leaa\basel{b_n}{k\cdot\chi_i}$ by~\eqref{eq:6}. 
Now we have obtained a contradiction
because by the induction hypothesis $\sigma([f])\leaa \mathbf k$,
and thus $\sigma([f])\ne[q]$. 
\end{proof}

Consider for each $i\in I$ and $j\in J$,
\begin{equation}
  \label{eq:57}
  b^i_j=\pos\bigl(\basel{1_{\cat(Y)}}{\chi_j}
  \land\sigma(\basel{1_{\cat(X)}}{\chi_i})\bigr)
  =\pos(\pi_j(\sigma(\basel{1_{\cat(X)}}{\chi_i}))).
\end{equation}
Since $\sigma$ is injective and $\ran(\sigma)$ is downwards closed, 
there is a unique $a^i_j\in\cat(X)$ such that
\begin{equation}
  \label{eq:70}
  \sigma(\<a^i_j,\chi_i\>)=\<b^i_j,\chi_j\>,
\end{equation}
and thus by proposition~\ref{p-25} we can define an isomorphism $\tau^i_j:\cat(X)_{a^i_j}\to\cat(Y)_{b^i_j}$ by
\begin{equation}
  \label{eq:71}
  \sigma(\<a,\chi_i\>)=\<\tau^i_j(a),\chi_j\>.
\end{equation}
Now  since claim~\ref{c-2} implies that $\sigma(\basel{1_{\cat(X)}}{\chi_i})\leaa\one$,
equation~\eqref{eq:70} implies equation~\eqref{eq:36}. 
And then equation~\eqref{eq:13}, with $h=\zero-$, 
follows from~\eqref{eq:36},~\eqref{eq:71} and continuity.
Note too that~\eqref{eq:35} holds because
$\<b_j,\chi_j\>\land\<b_{j^*},\chi_{j^*}\>$ for $j\ne j^*$.

\begin{claim}
\label{c-3}
For all $g$,
$\pos(\pi_j(\sigma([g])))=\pos(\pi_j(\sigma([n\cdot g])))$ for all $n=1,2,\dots$, 
and all $j\in J$.
\end{claim}
\begin{proof}
Since $\sigma$ is order preserving and $n\ge 1$,  $\sigma([g])\le\sigma([n\cdot g])$. 
It thusly suffices to show that 
$\pos(\pi_j(\sigma([n\cdot g])))-\pos(\pi_j(\sigma([g])))=0_{\cat(Y)}$.
Supposing to the contrary that the difference $a$ is nonzero, 
then $\basel{a}{\chi_j}\leaa\sigma([n\cdot g])$.
And then there is an $[f]\leaa [n\cdot g]$ such that $\sigma([f])=\basel{a}{\chi_j}$. 
However, 
$[f]\land[g]\ne[\zero-]$ and thus
$\sigma([f])\land\sigma([g])\eqaa\sigma([f]\land[g])\ne[\zero-]$, contradicting the fact
that $\basel a{\chi_j}\land \sigma([g])=[\zero-]$.
\end{proof}

This establishes by induction on $n=0,1,2,\dots$ that 
\begin{equation}
\label{eq:29}
\sigma(\basel a{n\cdot\chi_i})=n\cdot\sigma(\basel a{\chi_i}).
\end{equation} 
To see this, first note that by continuity and~\eqref{eq:36}, 
we may assume that $a\le a^i_j$ for some $j$.
By claim~\ref{c-2}, $\sigma(\basel a{n\cdot\chi_i})\leaa \mathbf n$, 
while by~\eqref{eq:71} and claim~\ref{c-3}, 
$\pos(\pi_{j^*}(\sigma(\basel a{n\cdot\chi_i})))=0_{\cat(Y)}$ for all $j^*\ne j$, and
$\pos(\pi_j(\sigma(\basel a{n\cdot\chi_i})))=\tau^i_j(a)$, and therefore 
$\sigma(\basel a{n\cdot\chi_i})\leaa\basel{\tau^i_j(a)}{n\cdot\chi_j}
=n\cdot\sigma(\basel a{\chi_i})$. 

As for the opposite inequality, supposing towards a contradiction that 
$\sigma(\basel a{n\cdot\chi_i})\laa\allowbreak\basel{\tau^i_j(a)}{n\cdot\chi_j}$, 
there exists $b\le\tau^i_j(a)$ such that
\begin{equation}
  \label{eq:25}
  \sigma(\basel a{n\cdot\chi_i})\land \basel{b}{n\cdot\chi_j}
  \eqaa \basel b{k\cdot\chi_j}
\end{equation}
for some $k=0,\dots,n-1$. 
Then letting $c$ be the element such that $\tau^i_j(c)=b$, 
$\sigma(\basel c{k\cdot\chi_i})=k\cdot\sigma(\basel c{\chi_i})
=k\cdot\basel b{\chi_j}=\basel b{k\cdot\chi_j}$ by the induction
hypothesis; however, 
$\sigma(\basel c{n\cdot\chi_i})\leaa\sigma(\basel a{n\cdot\chi_i})$ 
and we know that 
$\sigma(\basel c{n\cdot\chi_i})\leaa\basel b{n\cdot\chi_j}$ because we have
proved the other inequality of~\eqref{eq:29}. 
Thus $\sigma(\basel c{n\cdot\chi_i})\leaa\basel b{k\cdot\chi_j}$ by~\eqref{eq:25},
contradicting the fact that $\sigma$ is a monomorphism.

Equations~\eqref{eq:80} and~\eqref{eq:29} 
and continuity establish~\eqref{eq:14}, 
and~\eqref{eq:15} easily follows from the injectivity of $\sigma$.

The converse should be clear. 
\end{proof}

\begin{cor}
\label{o-4}
Order isomorphisms between $\baireaa{\irri I}$ and $\baireaa{\irri J}$ are precisely
those homomorphisms satisfying~\eqref{eq:14}--\eqref{eq:15}, $[h]=[\zero_{J-}]$
and the equation
\begin{equation}
  \label{eq:19}
  \spr_{i\in I}\tau_j^i(a_j^i)=1_{\cat(Y)}\espc\tu{for all $j\in J$},
\end{equation}
or equivalently $\ifm_{j\in J}\spr_{i\in I}\tau^i_j(a^i_j)=1_{\cat(Y)}$. 
\end{cor}

We would like to  generalize theorem~\ref{l-8} to the realm of Baire functions.
We begin by imposing topological requirements to apply the results of 
\Section\ref{sec:category-algebras}, in order to obtain a preliminary result 
(corollary~\ref{o-9}).

\begin{notation}
\label{notn:nodot}
In the context of notation~\ref{notn:basic}, we remove the dot when we do not want
to go to the equivalence class modulo $\meager(X)$, and thus
$\|\varphi(f_0,\dots,f_{n-1},a_0,\dots,\allowbreak a_{m-1})\|
=\{z\in X:\varphi(f_0(z),\dots,f_{n-1}(z),a_0,\allowbreak \dots,\allowbreak a_{m-1})\}$. 
Thus for example, for any $f:X\to\irri I$, $\|f(i)=n\|=\{z\in X:f(z)(i)=n\}$. 
Also, when we take a subset $B\subseteq X$ as opposed to a member of $\cat(X)$,
the notation $\<B,x\>$ ($x\in\irri I$) refers to the function $f:X\to \irri I$ such
that $f(z)=x$ for $z\in B$, and $f(z)=\zero_{I-}$ when $x\notin B$. 
\end{notation}

\begin{cor}
\label{o-9}
Let $X$ be a regular space, and suppose $Y$ has a Baire compactification.
Order embeddings $\sigma:\baireaa{X,\irri I}\to\baireaa{Y,\irri J}$ 
with convex range are completely determined by $B_j\in\bp(Y)$, 
$H_j\in\baire(B_j,X)$ such that
\begin{gather}
  \label{eq:64}
  H\inv_j[[\meager(X)]]\subseteq\meager(Y),\\
  \label{eq:84}
  H\inv_j[M]\cap \|g_j=i\|\in\meager(Y)\impls M\cap H_j[\|g_j=i\|]\in\meager(X)
\intertext{for every $i\in I$ and $M\subseteq X$,}
  \label{eq:91}
   H_j[\|g_j=i\|]\cap H_{j^*}[\|g_{j^*}=i\|]=\emptyset
  \espc\tu{for all $j\ne j^*$ in $J$, and all $i\in I$},\\
  \label{eq:85}
  \spr_{j\in J}\bigl[H_j[\|g_j=i\|]\bigr]{}_{\cat(X)}=1_{\cat(X)}
  \espc\tu{for all $i\in I$},
\end{gather}
and $g_j\in\cts(B_j,I)$ \tu{($j\in J$)} via $\sigma([f])=[\rho_f]+\sigma([\zero_{I-}])$ 
where 
\begin{equation}
  \label{eq:32}
  \rho_f(z)(j)=
\begin{cases}
f\bigl(H_j(z)\bigr)\bigl(g_j(z)\bigr)&z\in B_j,\\
0& z\notin B_j,
\end{cases}
\end{equation}
for all $j\in J$. 
\end{cor}
\begin{proof}
Let $a^i_j$ and $\tau^i_j$ ($i\in I$, $j\in J$) be as given to us by
theorem~\ref{u-2}. For each $i$ and~$j$, by corollary~\ref{o-15}, there exists a
relatively comeager $C^i_j\subseteq\bar\tau^i_j(a^i_j)$ and  a Baire
measurable injection $f^i_j:C^i_j\to\bar a^i_j$ satisfying 
\begin{equation}
  \label{eq:61}
  \tau^i_j([B])=[(f^i_j)\inv[B]]\espc\text{for all $B\in\bp(\bar a^i_j)$}
\end{equation}
and conditions~\eqref{item:2}--\eqref{item:18} of the conclusion of the
corollary (note that $A=\emptyset$ there since $\tau^i_j$ is an isomorphism).
Note that $\bar \tau^i_j(a^i_j)\cap \bar \tau^{i^*}_j(a^{i^*}_j)=\emptyset$ 
for all $i\ne i^*$ by~\eqref{eq:15}. 
For each $j\in J$, put $B_j=\bigcup_{i\in I}C^i_j$, and note that
$B_j\in\bp(Y)$ by disjointness and proposition~\ref{p-9}. 
Define $g_j:B_j\to I$ by
\begin{equation}
  \label{eq:62}
  g_j(z)=\text{the unique $i\in I$ such that $z\in C^i_j$.}
\end{equation}
Then define $H_j:B_j\to X$ by
\begin{equation}
  \label{eq:63}
  H_j(z)=f^{g_j(z)}_j(z).
\end{equation}
Then $H_j$ is Baire measurable since each $f^i_j$ is  and by disjointness.
Noting that $\|g_j=i\|=C^i_j$, 
equations~\eqref{eq:64} and~\eqref{eq:84} hold by condition~\eqref{item:7}.
Equation~\eqref{eq:91} holds by~\eqref{eq:35}, because each
$H_j[\|g_j=i\|]=\ran(f^i_j)\subseteq\bar a^i_j$.
And equation~\eqref{eq:85} holds by~\eqref{eq:36}, because each
$[H_j[\|g_j=i\|]]_{\cat(X)}=[\ran(f^i_j)]_{\cat(X)}=a^i_j$. 

\begin{claim}
\label{c-10}
$f\leaa f'$ implies $\rho_f\leaa \rho_{f'}$.
\end{claim}
\begin{proof}
Fixing $j\in J$, we must show that $\rho_f(z)(j)\le \rho_{f'}(z)(j)$ for almost all $z\in Y$. 
Now $\{z:\rho_f(z)(j)\nleq \rho_{f'}(z)(j)\}=\bigcup_{i\in I}D_i$ where
$D_i=\{z\in C^i_j:f(H_j(z))(i)\nleq f'(H_j(z))(i)\}$. And it follows from
equation~\eqref{eq:64} that each of the $D_i$'s is meager. But $\{C^i_j:i\in I\}$
is a pairwise disjoint family of relatively open subsets of $B_j$, 
and thus $\bigcup_{i\in I}D_i$ is meager by proposition~\ref{p-9}.
\end{proof}

The preceding claim entails that $f\eqaa f'$ implies $\rho_f\eqaa \rho_{f'}$, and thus
$[f]\mapsto[\rho_f]$ is well defined. Also note that indeed $\rho_f$ is Baire measurable:
By lemma~\ref{l-10} there is a Borel $f'\eqaa f$, and $f'\circ H$ is thus Baire
measurable being the composition of a Borel function with a Baire function.
The Baire measurability of $\rho_{f'}$ now easily follows from the continuity of the
$g_j$'s. Hence $\rho_f\eqaa \rho_{f'}$ implies $\rho_f\in\baire(Y,\irri J)$ by
proposition~\ref{p-30}. 

\begin{claim}
\label{c-11}
Whenever $\F\subseteq\baireaa{X,\irri I}$ has a supremum $[f']$, 
$\spr_{[f]\in\F}[\rho_f]=[\rho_{f'}]$.
\end{claim}
\begin{proof}
Let $[f']=\spr\F$. Then $[\rho_f]\leaa[\rho_{f'}]$ for all $[f]\in\F$, by
claim~\ref{c-10}. And $[\rho_{f'}]$ is moreover the least upper bound: For
if $\rho_{f'}\nleaa h'$ then there exists $j\in J$ such that not almost all $z\in Y$
satisfy $\rho_{f'}(z)(j)\le h'(z)(j)$. And then there is an $i\in I$ such that 
 $C=\{z\in C^i_j:\rho_{f'}(z)(j)\nleq h'(z)(j)\}\notin\meager(Y)$. And we can find a
 nonmeager $D\subseteq C$ and $n\in \N$ such that $h'(z)(j)=n$ for all $z\in D$.
Since $g_j(z)=i$ and $H_j(z)=f^i_j(z)$ for all $z\in C^i_j$, 
we have
\begin{equation}
  \label{eq:77}
  f'\bigl(f^i_j(z)\bigr)(i)=\rho_{f'}(z)(j)>n\espc\text{for all $z\in D$}.
\end{equation}
Note that condition~\eqref{item:7} for $f^i_j$ implies that $f^i_j[D]\notin\meager(X)$. 
Therefore, there exists $[f]\in\F$ such that 
\begin{equation}
  \label{eq:74}
  W=\{y\in f^i_j[D]:f(y)(i)> n\}\notin\meager(X).
\end{equation}
Now $\rho_{f}(z)(j)=f\bigl(f^i_j(z)\bigr)(i)>n=h'(z)(j)$ for all $z\in (f^i_j)\inv[W]$,
as $(f^i_j)\inv[W]\subseteq D$ since $f^i_j$ is an injection. 
Thus $\rho_f\nleaa h'$ because $(f^i_j)\inv[W]\notin\meager(Y)$ by
condition~\eqref{item:7}. 
\end{proof}

\begin{sublem}
\label{slem:1}
  $\displaystyle [\rho_f]=\spr_{i\in I}\spr_{n=0}^\infty n\cdot [\rho_{\<\|f(i)=n\|,\chi_i\>}]$.
\end{sublem}
\begin{proof}
Note that
\begin{equation}
  \label{eq:81}
  \rho_{f+f'}=\rho_{f}+\rho_{f'}.
\end{equation}
Thus $\rho_{n\cdot f}=n\cdot \rho_f$, and this establishes the sublemma because
$[\rho_f]=\spr_{i\in I}\spr_{n=0}^\infty \allowbreak
[\rho_{\<\|f(i)=n\|,n\cdot\chi_i\>}]$ by
equation~\eqref{eq:80} and claim~\ref{c-11}. 
\end{proof}

Assume without loss of generality that $\sigma([\zero_{I-}])=[\zero_{J-}]$.
To verify that $\sigma([f])=[\rho_f]$ for all $[f]$, it suffices by sublemma~\ref{slem:1}
and theorem~\ref{u-2} to take $[f]$ of the form $\basel {a}{\chi_i}$ for some $a\le a^i_j$, 
$i\in I$ and $j\in J$. Let $f$ be the representative as in~\eqref{eq:66}.
Then $\rho_f(z)(j^*)=0$ unless $g_{j^*}(z)=i$ and 
$H_{j^*}(z)\in \bar a$, in which case $\rho_f(z)(j^*)=1$.
Note that $g_{j^*}(z)=i$ implies $H_{j^*}(z)=f^i_{j*}(z)\in \bar a^i_{j^*}$, and thus
$H_{j^*}(z)\in\bar a\subseteq\bar a^i_j$ implies $j^*=j$ by~\eqref{eq:35}.
Therefore, $[\rho_f]=\basel{[(f^i_j)\inv[\bar a]]}{\chi_j}=\basel{\tau^i_j(a)}{\chi_j}$ 
by~\eqref{eq:61}. 
Now by~\eqref{eq:13}, we have shown $[\rho_f]=\sigma([f])$.
\end{proof}

\begin{remark}
\label{r-2}
For second countable spaces we can improve to open and continuous $H_j$ ($j\in J$). 
More precisely, if $X$ is a second countable Hausdorff space and $Y$ is
Polish, then we can obtain the result of corollary~\ref{o-9} such that the
$H_j:B_j\to X$ are moreover open and continuous.
The proof is the same as for corollary~\ref{o-9}, but with the stronger hypothesis
we can use corollary~\ref{u-8} instead of corollary~\ref{o-15} 
to obtain topological embeddings
$f^i_j:B^i_j\to\bar\tau^i_j(a^i_j)$ for relatively comeager 
$B^i_j\subseteq\bar a^i_j$ ($i\in I$, $j\in J$). Now we let $C^i_j=\ran(f^i_j)$.
Then $H_j$ is defined instead by
\begin{equation}
  \label{eq:78}
  H_j(z)=\bigl(f^{g_j(z)}_j\bigr)\inv(z)
\end{equation}
and is in fact open and continuous. The rest of the proof goes through the same.
\end{remark}

Next we want to generalize the notion of a projection 
by a function (cf.~definition~\ref{d-3}). This requires a generalization of
functions with codomain a power of~$\N$. 

\begin{notation}
\label{notn:product}
Let $J$ be an index set, and let $X$ and $Y$ be topological spaces. 
Suppose $K:Y\to\power(J)$. Then $X^{K-}$ denotes the function $z\mapsto X^{K(z)}$, 
and thus $\prod X^{K-}$ ($=\prod_{z\in Y}X^{K(z)}$) denotes the family of all functions
$f$ with domain  $Y$ such that
\begin{equation}
  \label{eq:121}
  f(z)\in X^{K(z)}\espc\text{for all $z\in Y$}.
\end{equation}
Suppose $f\in\prod X^{K-}$. For each $j\in J$, we let $f^j$ denote the function
\begin{equation}
  \label{eq:122}
  f^j:\|j\in K\|\to X
\end{equation}
where $f^j(z)=f(z)(j)$ for all $z\in\|j\in K\|$. We let $\baire(X^{K-})$ denote the
family of all $f\in\prod X^{K-}$ such that $f^j\in\baire(\|j\in K\|,X)$ for all
$j\in J$, and we let $\cts(X^{K-})$ denote the family of all $f\in\prod X^{K-}$ 
such that $f^j\in\cts(\|j\in K\|,X)$ for all $j\in J$. 
\end{notation}

Note that this notation is of course a generalization, i.e.~if
$I\subseteq J$ is some fixed index set then letting $K$ be the function on $Y$ constantly
equal to $I$, $\prod X^{K-}=(X^I)^Y$.

\begin{defn}
\label{d-13}
Fix index sets $I$ and $J$. 
For each function $K:Y\to\power(J)$, 
and each $g\in \prod I^{K-}$ and $H\in \prod X^{K-}$, 
we define $\pi_{g,H}:(\irri I)^X\to\prod \irri{K-}$  by
\begin{equation}
  \label{eq:34}
  \pi_{g,H}(f)(z)(j)=f\bigl(H(z)(j)\bigr)\bigl(g(z)(j)\bigr)
  \espc\text{for all $j\in K(z)$, for all $z\in Y$}.
\end{equation}
\end{defn}

Note that this is indeed a generalization of the projections. 
For if $f\in(\irri I)^X$ is constant with respect to $z\in X$, say $f(z)=x$,
then $H$ is irrelevant and
\begin{equation}
  \label{eq:39}
  \pi_{g,H}(f)(z)=x\circ g(z)\espc\text{for all $z\in Y$}.
\end{equation}

\begin{thm}
\label{o-8}
Let $X$ be a regular space and suppose $Y$ has a Baire compactification.
The order embeddings from $\baireaa{X,\irri I}$ into $\baireaa{Y,\irri J}$ with
convex range consist of maps of the form
\begin{equation}
  \label{eq:33}
  \sigma([f])=[\pi_{g,H}(f)\bigext\zero_{J\setminus K-}]+[h]
  \espc\tu{for all $[f]\in\baireaa{X,\irri I}$}
\end{equation}
for some $K\in\baire(Y,\power(J))$, $g\in\cts(I^{K-})$,
$H\in\baire(X^{K-})$ such that
\begin{gather}
  \label{eq:17}
  (H^j)\inv[M]\in\meager(Y)\espc\tu{for all $M\in\meager(X)$ and all $j\in J$},\\
  \label{eq:86}
  (H^j)\inv[M]\cap\|g^j=i\|\in\meager(Y)\impls M\cap H^j[\|g^j=i\|]\in\meager(X)
  \intertext{for every $i\in I$ and $M\subseteq X$,}
  \label{eq:87}
  H^j[\|g^j=i\|]\cap H^{j^*}[\|g^{j^*}=i\|]=\emptyset
  \espc\tu{for all $j\ne j^*$ in $J$, and all $i\in I$},\\
  \label{eq:82}
  \spr_{j\in J}\bigl[H^j[\|g^j=i\|]\bigr]{}_{\cat(X)}=1_{\cat(X)}
  \espc\tu{for all $i\in I$},
\end{gather}
and some $[h]\in\baireaa{Y,\irri J}$.
Moreover, the range is given by
\begin{equation}
  \label{eq:88}
  \ran(\sigma)={\bigl(\baire(\irri{K-})\bigext\zero_{J\setminus K-}
    +h\bigr)}\div{\eqaa}.
\end{equation}
\end{thm}
\begin{proof}
Let $\sigma$ be such an order embedding with convex range.
Let $B_j\in\bp(Y)$, $H_j\in\baire(B_j,X)$,
$g_j\in\cts(B_j,I)$ ($j\in J$) and $\rho$ be as given by corollary~\ref{o-9}.
Define $K:Y\to \power(J)$ by
\begin{equation}
  \label{eq:3}
  K(z)=\{j\in J:z\in B_j\}.
\end{equation}
Clearly $K\in\baire(Y,\power(J))$. Define $g\in\cts(I^{K-})$ by
\begin{equation}
  \label{eq:16}
  g(z)(j)=g_j(z)\espc\text{for all $j\in K(z)$, for all $z\in Y$},
\end{equation}
so that $g^j=g_j$ for all $j\in J$. 
And we can define $H\in\baire(X^{K-})$ by $H(z)(j)=H_j(z)$ 
for all $j\in K(z)$ and all $z\in Y$, so that $H^j=H_j$ for all $j$. 
Equations~\eqref{eq:17}--\eqref{eq:82} 
will be satisfied by~\eqref{eq:64}--\eqref{eq:85}, respectively.  Then
\begin{equation}
  \rho_f(z)(j)=f\bigl(H_j(z)\bigr)\bigl(g_j(z)\bigr)=f\bigl(H(z)(j)\bigr)\bigl(g(z)(j)\bigr)
  \espc\text{for all $j\in K(z)$},\label{eq:38}
\end{equation}
and $\rho_f(z)(j)=0$ for all $j\in J\setminus K(z)$. Hence
$\rho_f={\pi_{g,H}}\bigext\zero_{J\setminus K-}$, and therefore, 
$\sigma([f])=[\pi_{g,H}(f)\bigext\zero_{J\setminus K-}]+\sigma([\zero_{I-}])$ 
for all $[f]\in\baireaa{X,\irri I}$. 

Conversely, suppose we are given $K$, $g$, $H$ and $h$ as in the corollary. 
Let $\sigma:\baireaa{X,\irri I}\to\baireaa{Y,\irri J}$ be as defined
in~\eqref{eq:33}. 
Equation~\eqref{eq:17} and the continuity of the $g^j$'s 
allows us to use the same argument 
as in the proof of claim~\ref{c-10} to prove that $f\leaa f'$ implies
$\pi_{g,H}(f)\leaa\pi_{g,H}(f')$, using $\{(g^j)\inv:i\in I\}$ here in place of
$\{C^i_j:i\in I\}$ there.  This establishes that $\sigma$ is a well defined
partial order homomorphism.
To prove that $\sigma$ is an embedding,
it remains to show that $\pi_{g,H}$ is order
reflecting.  Suppose then that $f\nleaa f'$. There is an $i\in I$ such that not
almost all $y\in X$ satisfy $f(y)(i)\le f'(y)(i)$. Thus by~\eqref{eq:82}, there
exists $j\in J$ such that
\begin{equation}
  \label{eq:83}
  \{y\in X:f(y)(i)>f'(y)(i)\}\cap H^j[\|g^j=i\|]\notin\meager(X).
\end{equation}
Therefore, $C=(H^j)\inv[y:f(y)(i)>f'(y)(i)]\cap\|g^j=i\|\notin\meager(Y)$ by
equation~\eqref{eq:86}. However, for every $z\in C$,
$\pi_{g,H}(f)(z)(j)=f\bigl(H^j(z)\bigr)\bigl(g^j(z)\bigr)
=f\bigl(H^j(z)\bigr)(i)>f'\bigl(H^j(z)\bigr)(i)=\pi_{g,H}(f')(z)(j)$, proving
$\pi_{g,H}(f)\nleaa\pi_{g,H}(f')$ as required.

For the converse, we still have yet to show that $\ran(\sigma)$ is convex.
However, since $\baire(\irri{K-})\bigext\zero_{J\setminus K-}$ is
evidently downwards closed, the family on the right hand side of
equation~\eqref{eq:88} is convex. Thus it will suffice to establish~\eqref{eq:88}.
Without loss of generality we may assume that $[h]=[\zero-]$. 
Take $h'\in\baire(\irri{K-})$. Using~\eqref{eq:87}, it is possible to
find $f:X\to\irri I$ such that for each $i\in I$,
\begin{equation}
  \label{eq:90}
  f\bigl(H^j(z)\bigr)(i)=h'(z)(j)\espc\text{for all $z\in\|g^j=i\|$},
\end{equation}
for all $j\in J$. Note that~\eqref{eq:82} guarantees that $f$ is Baire measurable. 
Since each $g^j$ is continuous, and since
equation~\eqref{eq:90} entails that 
$\pi_{g,H}(f)(z)(j)=f\bigl(H^j(z)\bigr)(i)=h'(z)(j)$ for almost all $z\in\|g^j=i\|$, for all
$i\in I$, we obtain $\pi_{g,H}(f)\eqaa h'$ as desired.
\end{proof}

\begin{remark}
\label{r-3}
When $X$ is a second countable Hausdorff space and $Y$ is Polish, 
we can moreover obtain $H\in\cts(X^{K-})$ with $H^j$ an open continuous
mapping for all ~$j\in J$, by remark~\ref{r-2}. 
\end{remark}

\begin{cor}
\label{o-17}
Let $X$ be regular and let $Y$ be a space with a Baire compactification. 
Any embedding $\sigma$ of $\baireaa{X,\irri I}$ into a downwards closed subset of
$\baireaa{Y,\irri J}$ is of the form
\begin{equation}
  \label{eq:89}
  \sigma([f])=[\pi_{g,H}(f)\bigext\zero_{J\setminus K-}]
\end{equation}
for some $K\in\baire(Y,\power(J))$, $g\in\cts(I^{K-})$ and
$H\in\baire(X^{K-})$. It is thus a monoid homomorphism. Moreover, its range is
the equivalence classes of the collection of functions 
$\baire(\irri{K-})\bigext\zero_{J\setminus K-}$. 
\end{cor}
\begin{proof}
That order embeddings with downwards closed ranges have the indicated form and the
indicated range, is an immediate consequence of theorem~\ref{o-8} 
since a downwards closed range implies $[h]=[\zero-]$. 
As for being a monoid homomorphism, observe that 
\begin{equation}
\pi_{g,H}(f+f')(z)=\pi_{g,H}(f)(z)+\pi_{g,H}(f')(z)\espc\text{for all $z$}.\qedhere
\end{equation} 
\end{proof}

\begin{cor}
\label{o-7}
Let $X$ be regular and let $Y$ be a space with a Baire compactification. 
The order isomorphisms between $\baireaa{X,\irri I}$ and $\baireaa{Y,\irri J}$ 
consist of maps of the form
\begin{equation}
  \label{eq:110}
  [f]\mapsto[\pi_{g,H}(f)]
\end{equation}
for some $g\in\cts(Y,I^{J})$ and $H\in\baire(Y,X^{J})$ such that
\begin{gather}
  \label{eq:37}
  (H^j)\inv[M]\in\meager(Y)\espc\tu{for all $M\in\meager(X)$ and all $j\in J$},\\
  \label{eq:92}
  (H^j)\inv[M]\cap\|g^j=i\|\in\meager(Y)\impls M\cap H^j[\|g^j=i\|]\in\meager(X)
  \intertext{for every $i\in I$ and $M\subseteq X$,}
  \label{eq:93}
  H^j[\|g^j=i\|]\cap H^{j^*}[\|g^{j^*}=i\|]=\emptyset
  \espc\tu{for all $j\ne j^*$ in $J$, and all $i\in I$},\\
  \label{eq:94}
  \spr_{j\in J}\bigl[H^j[\|g^j=i\|]\bigr]{}_{\cat(X)}=1_{\cat(X)}
  \espc\tu{for all $i\in I$}.
\end{gather}
\end{cor}
\begin{proof}
$g$ and $H$ are obtained from theorem~\ref{o-8}. Corollary~\ref{o-17} tells us that
$h$ vanishes, and that $K(z)=J$ for almost all $z$.
The function $g$ is continuous because $K$ is constant and the $g^j$'s are
continuous. 
\end{proof}
\subsection{Continuous functions into powers of $\N$}
\label{sec:cont-funct-into}

For any topological space $X$, $\cts(X,\irri I)\subseteq\baire(X,\irri I)$.
It is in fact a sublattice, 
i.e.~it is closed under the lattice operations on $\baire(X,\irri I)$,
of pointwise $\max$ and $\min$ (cf.~equation~\eqref{eq:10}).

\begin{prop}
\label{p-41}
$\cts(X,\irri I)$ is a sublattice of $(\baire(X,\irri I),\le)$.
\end{prop}
\begin{proof}
Take $f,g\in\cts(X,\irri I)$.
The coordinatewise mappings $\max,\min:\irri I\times\irri I\to \irri I$
(cf.~\Section\ref{sec:powers-n}) are both continuous, 
and $(f,g):X\to\irri I\times\irri I$ is continuous. Therefore,
$\max\{f,g\}={\max}\circ{(f,g)}$ and $\min\{f,g\}={\min}\circ(f,g)$ are both continuous. 
\end{proof}

Furthermore, $\cts(X,\irri I)$ is a submonoid of $(\baire(X,\irri I),+)$. In fact, 
noting that addition and subtraction are continuous operations on $\intgr^I$,
$\cts(X,\intgr^I)$ is closed under subtraction as in the proof of
proposition~\ref{p-41}, and therefore is a subgroup of the Abelian group
$(\baire(X,\intgr^I),+)$. Similarly:

\begin{prop}
\label{p-55}
$\cts(X,\irri I)$ is a submonoid of $(\baire(X,\irri I),+)$ closed under subtraction.
\end{prop}

\noindent Thus $\cts(X,\irri I)$ is essentially a `subgroup' of
$(\baire(X,\irri I),+)$ (see proposition~\ref{p-51}); indeed, note that
$\<\cts(X,\irri I)\>=\cts(X,\intgr^I)$.

Not only are members of $\cts(X,\irri I)$ representatives of members 
of $\baireaa{X,\irri I}$, but 
when $X$ is a Baire space, $\cts(X,\irri I)$ embeds into $\baireaa{X,\irri I}$.

\begin{prop}
\label{p-22}
If $X$ is a Baire space, then $f\mapsto[f]$ is an order embedding 
of $(\cts(X,\allowbreak\irri I),\allowbreak\le)$ into $(\baireaa{X,\irri I},\leaa)$,
and it is also both a monoid and a lattice embedding. 
\end{prop}
\begin{proof}
$f\mapsto[f]$ is an order embedding, 
because in a Baire space the comeager sets are dense.
It is a monoid homomorphism 
because, as was observed in~\Section\ref{sec:baire-functions},
$(\baireaa{X,\irri I},+)$ is a quotient monoid. 
And it is a lattice homomorphism by equation~\eqref{eq:12}.  
\end{proof}

\begin{lem}
\label{l-13}
Let $X$ be a zero dimensional space. 
Then $C(X,\irri I)$ is dense in the quasi order $(\baire(X,\irri I),\leaa)$.
\end{lem}
\begin{proof}
Take $f\in\baire(X,\irri I)$ with $[f]\ne[\zero_{I-}]$. 
Then for some $i\in I$, $f^i\eqnaa 0$ (i.e.~not $f^i\eqaa 0$), 
and thus there exists $n=1,2,\dots$ such
that $B=\{z\in X:f(z)(i)=n\}\notin\meager(X)$. Since $X$ is zero dimensional, there
exists a nonempty clopen set $C$ with $C\setminus B\in\meager(X)$.
Now $\<C,\chi_i\>\in \cts(X,\irri I)^+$ (cf.~notation~\ref{notn:nodot})
since $C$ is clopen, and clearly $\<C,\chi_i\>\leaa f$. 
\end{proof}

\begin{cor}
\label{o-36}
Let $X$ be a zero dimensional Baire space. 
Then $C(X,\irri I)$ densely embeds into $(\baire(X,\irri I),\leaa)$.
\end{cor}
\begin{proof}
Proposition~\ref{p-22} and lemma~\ref{l-13}. 
\end{proof}

\begin{example}
\label{x-13}
The irrationals $\irr$ form a zero dimensional Baire space, and indeed this is often
called the \emph{Baire space}. More generally, so is $\irri H$ for any $H$, and thus
$C(\irri H,\irri I)$ densely embeds into $\baireaa{\irri H,\irri I}$. 
\end{example}

The embedding is generally not onto. 

\begin{example}
\label{x-14}
Since the irrationals of the real line $\reals\setminus\rationals$ are homeomorphic
to our set of irrationals $\irr$,
$C(\reals\setminus\rationals,\reals\setminus\rationals)$ 
densely embeds into $\baireaa{\reals\setminus\rationals,\reals\setminus\rationals}$
(we can make sense of the latter using the fact that $\eqaa$ is the same as $\{=\}_\aa$ on
$\irr$). Now fix an irrational number $\sigma\in\reals$. Then a jump function
$f:\reals\setminus\rationals\to\reals\setminus\rationals$ where
\begin{equation}
  \label{eq:96}
  f(x)=\begin{cases}
    0,&\text{if $x\le\sigma$,}\\
    1,&\text{if $x>\sigma$,}
    \end{cases}
\end{equation}
is a Baire measurable function, but no function
$g:\reals\setminus\rationals\to\reals\setminus\rationals$ with $g\eqaa f$ is
continuous at $\sigma$. Therefore there is a function $f\in\baire(\irr,\irr)$ with no
member of $\cts(\irr,\irr)$ in its equivalence class modulo $\eqaa$. 
\end{example}

\begin{prop}
\label{p-56}
Let $X$ be a Baire space. Then $C(X,\irri I)$ identifies \tu(via $f\mapsto[f]$\tu)
with a submonoid of $(\baireaa{X,\irri I},+)$ that is closed under subtraction.
It also identifies with a sublattice of $(\baireaa{X,\irri I},\leaa)$.
\end{prop}
\begin{proof}
Propositions~\ref{p-55} and~\ref{p-22}. 
\end{proof}

When $X$ is a zero dimensional Baire space, $C(X,\irri I)$ identifies with a subset
of $\baireaa{X,\irri I}$ with several nice properties.

\begin{lem}
\label{l-14}
Let $X$ be a zero dimensional Baire space. 
Then $C(X,\irri I)$ identifies \tu(via $f\mapsto[f]$\tu) with a subset of the
lattice monoid $(\baireaa{X,\irri I},+)$ satisfying the following properties\tu:
It is regular, it forms a basis and it is strongly interval predense.
\end{lem}
\begin{proof}
First we show that $C(X,\irri I)$ is (i.e.~identifies with) a basis.
Since the basic elements form a basis by proposition~\ref{p-35}, 
it suffices to prove that each basic element $\<a,n\cdot\chi_i\>$ 
is the supremum of an antichain of $C(X,\irri I)$. 
And this is so because it is the supremum of an antichain of
elements of $C(X,\irri I)$ of the form $\<C,n\cdot\chi_i\>$ with $C\subseteq X$
clopen, by example~\ref{x-19}. 

We already know that $(\baireaa{X,\irri I},+)$ is a cancellative commutative monoid.
Therefore, $C(X,\irri I)$ is a regular subset by lemma~\ref{l-27}, because it is a
dense submonoid closed under subtraction by corollary~\ref{o-36} and
proposition~\ref{p-56}. 

In theorem~\ref{l-2} (see also remark~\ref{r-6}), it was established that
$\baireaa{X,\irri I}$ is moreover a complete semilattice monoid satisfying the \jid.
Therefore, since we have shown that $C(X,\irri I)$ is a submonoid forming a basis,
it is strongly interval predense by lemma~\ref{l-21}. 
\end{proof}

\begin{lem}
\label{u-5}
Let $X$ be a zero dimensional Baire space, and let $M$ be a complete semilattice
that is flat-complete and \jid.
Then every order embedding of $\cts(X,\irri I)$ into $M$ with preregular range 
uniquely extends to a continuous lattice embedding of
$\baireaa{X,\irri I}$ into $M$.
\end{lem}
\begin{proof}
We know that $\baireaa{X,\irri I}$ is a \jid\ lattice.
By lemma~\ref{l-14}, $C(X,\irri I)$ forms a strongly interval predense basis.
Hence every order embedding of $\cts(X,\irri I)$ into $M$ with preregular range has
an extension, unique on $C(X,\irri I)^+$, to a continuous lattice embedding
by corollary~\ref{u-10} (or alternatively, by corollary~\ref{o-33}).
And the extension is in fact unique because $C(X,\irri I)$, being a submonoid,
contains the zero of $\baireaa{X,\irri I}$ (cf.~remark~\ref{r-7}).
\end{proof}

\begin{thm}
\label{o-18}
Let $X$ and $Y$ be a zero dimensional Baire spaces.
Then every order embedding of $\cts(X,\irri I)$ into $\cts(Y,\irri J)$ 
with convex range has a unique extension to a continuous 
lattice embedding of $\baireaa{X,\irri I}$
into $\baireaa{Y,\irri J}$ with convex range.
\end{thm}
\begin{proof}
We know that $\baireaa{X,\irri I}$ is a cancellative commutative complete semilattice
monoid satisfying the \jid, and so is $\baireaa{Y,\irri J}$, and they are moreover
flat-complete by theorem~\ref{l-2}. Since $C(X,\irri I)$ is a submonoid forming a
basis and $C(Y,\irri J)$ is a dense submonoid closed under subtraction, the
conclusion is by corollary~\ref{o-34}.
\end{proof}

\begin{thm}
\label{u-7}
Let $X$ be a regular zero dimensional Baire space, 
and let $Y$ be a zero dimensional space with a Baire compactification. 
The order embeddings from $\cts(X,\irri I)$ into $\cts(Y,\irri J)$
with convex range consist of maps of the form 
\begin{equation}
  \label{eq:98}
  \sigma={\pi_{g,H}}\bigext\zero_{J\setminus K-}+h
\end{equation}
for some $K\in\cts(Y,\power(J))$, $g\in\cts(I^{K-})$,
$H\in\cts(X^{K-})$ such that
\begin{gather}
  \label{eq:102}
  (H^j)\inv[M]\in\meager(Y)\espc\tu{for all $M\in\meager(X)$ and all $j\in J$},\\
  \label{eq:103}
  (H^j)\inv[M]\cap\|g^j=i\|\in\meager(Y)\impls M\cap H^j[\|g^j=i\|]\in\meager(X)
  \intertext{for every $i\in I$ and $M\subseteq X$,}
  \label{eq:104}
  H^j[\|g^j=i\|]\cap H^{j^*}[\|g^{j^*}=i\|]=\emptyset
  \espc\tu{for all $j\ne j^*$ in $J$, and all $i\in I$},\\
  \label{eq:105}
  \spr_{j\in J}\bigl[H^j[\|g^j=i\|]\bigr]{}_{\cat(X)}=1_{\cat(X)}
  \espc\tu{for all $i\in I$},
\end{gather}
and some $h\in\cts(Y,\irri J)$. Moreover, the range is given by
\begin{equation}
  \label{eq:131}
  \ran(\sigma)=\cts(Y,\irri K)\bigext\zero_{J\setminus K-}+h.
\end{equation}
\end{thm}
\begin{proof}
Let $\sigma:\cts(X,\irri I)\to\cts(Y,\irri J)$ be an order embedding 
with convex range. 
Then, noting that $Y$ is a Baire space (proposition~\ref{p-57}), 
theorem~\ref{o-18} gives a unique embedding $\bar\sigma$ of
$\baireaa{X,\irri I}$ into $\baireaa{Y,\irri J}$ with convex range, such that
\begin{equation}
  \label{eq:97}
  \bar\sigma([f])=[\sigma(f)]\espc\text{for all $f\in\cts(X,\irri I)$.}
\end{equation}
Now we apply theorem~\ref{o-8} to obtain $K\in\baire(Y,\power(J))$, 
$g\in\cts(I^{K-})$, $H\in\baire(X^{K-})$ and $h\in\baire(Y,\irri J)$ satisfying
equations~\eqref{eq:102}--\eqref{eq:105}, and
$\bar\sigma=[{\pi_{g,H}}\bigext\zero_{J\setminus K-}+h]$. 
Now
 also equation~\eqref{eq:98} follows from~\eqref{eq:97}.
It remains to show that $K$ and $h$ are in fact continuous, and that
$H\in\cts(X^{K-})$. 

We see that $h$ is continuous by plugging $\zero_{I-}$ into~\eqref{eq:98}.
Thus, by equation~\eqref{eq:98} 
and closure under subtraction (proposition~\ref{p-55}), 
\begin{equation}
\label{eq:45}
{\pi_{g,H}(f)}\bigext\zero_{J\setminus K-}\text{ is continuous}
\espc\text{for all $f\in\cts(X,\irri I)$.} 
\end{equation}
We can use equation~\eqref{eq:39} to see that
$\pi_{g,H}(\one_{I-})={\one_{K-}}\bigext\zero_{J\setminus K-}$.
If $K$ is not continuous, then clearly neither is
${\one_{K-}}\bigext\zero_{J\setminus K-}$, contradicting~\eqref{eq:45} with
$f=\one_{I-}$. 

Suppose now that $H\notin\cts(X^{K-})$. Then for some $j\in J$, 
$H^j$ is noncontinuous at some $z_0\in\|j\in K\|$, 
say $U\ni H^j(z_0)$ is an open subset of $X$ 
and $H^j[V]\nsubseteq U$ for every open $V\ni z_0$. 
Since $g^j$ is continuous, there exists an open $V\ni z_0$ with 
$g^j(z)=i$ for all $z\in V$. 
Since $C(X,\irri I)$ is dense (corollary~\ref{o-36}), 
we can find a continuous $f:X\to\irri I$ such that
$f(H^j(z_0))(i)=1$ but $f(H^j(z))(i)=0$ for all $z\notin U$. Now
$\pi_{g,H}(f)\bigext\zero_{J\setminus K-}$ is
not continuous at $z_0$ contrary to~\eqref{eq:45}, 
because $\pi_{g,H}(f)(z_0)(j)=1$ 
but every open $V\supseteq V'\ni z_0$ has a $z\in V'$ with $H^j(z)\notin U$ and thus
$\pi_{g,H}(f)(z)(j)=0$. 

For the converse, it is clear that $\pi_{g,H}(f)$ is continuous for all $f\in\cts(X,\irri
I)$, whenever $K$ and is continuous, 
$g\in\cts(I^{K-})$ and $H\in\cts(X^{K-})$. And thus it follows from
theorem~\ref{o-8} that equation~\eqref{eq:98} determines an order
embedding. As for having convex range, this is automatic once
equation~\eqref{eq:131} is established. However, one can use the
properties~\eqref{eq:102}--\eqref{eq:105} to prove that
$\pi_{g,H}(f)\bigext\zero_{J\setminus K-}$ is noncontinuous whenever $f$ is
noncontinuous. Thus $\ran(\sigma)$ is equal to the intersection of $\cts(X,\irri I)$
with the set of representatives of $\ran(\bar\sigma)$. Now equation~\eqref{eq:131}
is an immediate consequence of equation~\eqref{eq:88}. 
\end{proof}
\subsection{Further directions}
\label{sec:further-directions}

We suggest some classes of partial orders (out of infinitely many possibilities)
for which we feel it would be interesting
to know what the partial order embeddings with convex range are. 

\subsubsection{${\baire(X,\power(S))}\div{\eqaa}$}
\label{sec:bairex-powersdiveqaa}

In \Section\ref{sec:new-subsection} we examined embeddings with convex range 
between arbitrary power set algebras. We expect that this can be generalized to the
quotient of Baire functions into power set algebras modulo almost always equality,
in a completely analogous manner to \Section\ref{sec:baire-functions}. 

Fixing some set $S$, the natural analogue of $\eqaa$ for $\power(S)$ is given by the
family $\{R_s:s\in S\}$ where $x \mathrel{R_s} y$ iff ($s\in x$ iff $s\in y$). And
$\subseteqaa$ is defined analogously. We should then be able to go through the same
analysis for embeddings between posets of the form
$(\baireaa{X,\power(S)},\subseteqaa)$, to obtain analogues of theorems~\ref{u-2}
and~\ref{o-8}. Also note that one obtains an analogous connection with set theoretic
forcing, e.g.~for $f,g\in\baire(X,\power(S))$, 
\begin{align}
  \cat(X)&\forces\dot f=\dot g\IFF f\eqaa g, \\
  \cat(X)&\forces\dot f\subseteq\dot g\IFF f\subseteqaa g.
\end{align}

\subsubsection{${L^0(\mu,\irri I)}\div{\eqae}$}
\label{sec:l0mu-irri-idiveqae}

Let $(X,\mu)$ be a measure space, and $S$ a topological space. 
For simplicity let us insist that the measure is totally finite
(i.e.~$\mu(X)<\infty$) or at least $\sigma$-finite.
Then we denote the family of all $\mu$-measurable functions from $X$ into $Y$ by
$L^0(\mu,S)$, i.e.~functions $f:X\to S$ such that $f\inv[U]$ is in the domain of the
measure $\mu$ for every open $U\subseteq S$. Then analogously to the families of
Baire measurable functions, for a set of relations $\R$ on $S$, we define a relation
$\R_\ae$ on $L^0(\mu,S)$ by
\begin{equation}
  \label{eq:106}
  f\mathrel{\R_\ae}g\If \ifm_{R\in\R}\ulc f(z)\rel g(z)\text{ for almost every }z\in
  X\urc,
\end{equation}
where ``for almost every'' is interpreted as the complement of a measure zero subset
of~$X$. Using the same sets of relations as in example~\ref{x-7}, 
we obtain $\leae$ and $\eqae$. 
And then we can form the quotient poset $({L^0(\mu,\irri I)}\div{\eqae},\leae)$. 
We have now obtained the measure theoretic analogue of the
quotient lattice of \Section\ref{sec:baire-functions}. We again have an analogue with
set theoretic forcing, but now our forcing notion is the complete Boolean algebra
$\random(\mu)=\dom(\mu)\div\nulls_\mu$ where $\nulls_\mu$ denotes the ideal 
of $\mu$-measure zero subsets of $X$ (i.e.~we are forcing with a \emph{measure
  algebra}; more precisely, $(\random(\mu),\nu)$, with $\nu([a])=\mu(a)$, is the
measure algebra of the measure space $(X,\mu)$). 
This is known as \emph{random forcing} and was invented by
Solovay~\cite{MR0265151} to prove the consistency of all subsets of the real line being
Lebesgue measurable with the usual axioms of mathematics 
minus the Axiom of Choice. 
By way of analogy, we have e.g.~for all $f,g\in L^0(\mu,\irri I)$,
\begin{align}
  \random(\mu)&\forces\dot f=\dot g\IFF f\eqae g,\\
  \random(\mu)&\forces\dot f\le\dot g\IFF f\leae g
\end{align}
(see e.g.~\cite{MR1784706}). 

Moreover, theorem~\ref{u-2} also describes the order embeddings between
${L^0(\mu,\irri I)}\div{\eqae}$ and ${L^0(\nu,\irri J)}\div{\eqae}$ with convex
range, by replacing $\cat(X)$ and $\cat(Y)$ with $\random(\mu)$ 
and $\random(\nu)$,
respectively. In fact, the proof of theorem~\ref{u-2} does not use properties
specific to $\cat(X)$ other than completeness; and this theorem can be generalized
to arbitrary quotients of families of functions measurable with respect to some
$\sigma$-ideal on a fixed set $X$, provided the ideal satisfies the \emph{countable
  chain condition}, and both the meager ideal and the ideal of $\mu$-measure zero
sets satisfy this property. 

The question here is whether we can obtain a measure theoretic analogue of
theorem~\ref{o-8}. 

\begin{question}
Supposing that $(X,\mu)$ and $(Y,\nu)$ are \tu`reasonable\tu' measure spaces, 
can we describe the order
embeddings of ${L^0(\mu,\irri I)}\div{\eqae}$ into ${L^0(\nu,\irri J)}\div{\eqae}$
with convex range, in an analogous way to theorem~\tu{\ref{o-8}}, using measurable
functions $K:Y\to\power(J)$, $g:Y\to I^{K-}$, $H:Y\to X^{K-}$ and
$h:Y\to\irri J$\tu?
\end{question}

\subsubsection{$(\pnfin,\subseteqfnt)$}
\label{sec:pnfin-subseteqfnt}

Consider now the power set $\pN$ of $\N$ quasi ordered by ``almost inclusion'',
i.e.~$a\subseteqfnt b$ if $a\setminus b$ is finite. The quotient
$(\pnfin,\subseteqfnt)$ over the equivalence relation $\eqfnt$ given by $a\eqfnt b$
if $a\subseteqfnt b$ and $b\subseteqfnt a$ (i.e.~$a\diff b$ is finite) is a Boolean
algebra of great importance in a number of areas of mathematics. 
In fact, there is a general program of research aimed at investigating the class of
partial orders that embed into $\pnfin$, to which the present paper is clearly relevant.
To give one example, this class of orders played a major role 
in the solution, by Solovay--Woodin (see~\cite{MR942216}), 
of a famous problem of Kaplansky~\cite{MR0031193} on automatic
continuity in Banach algebras.

The structure of the order automorphism group (equivalently, Boolean algebra
automorphisms) is independent of the usual axioms of mathematics, 
as was established by Shelah in his celebrated result~\cite{MR675955}:

\begin{thm}[Shelah]
  \label{u-9}
It is consistent that every order automorphism of $(\pnfin,\subseteqfnt)$ is
\emph{trivial}. 
\end{thm}

\noindent A \emph{trivial} automorphism of $(\pnfin,\subseteqfnt)$ is one of the form
\begin{equation}
  \label{eq:107}
  \sigma([a])=\bigl[h[a]\bigr]
\end{equation}
for some bijection $h:A\to B$ where $A,B\subseteq\N$ are both cofinite 
(see also below). This should be compared to corollary~\ref{o-5}. 
Note that it was already known that there
exist nontrivial automorphisms under $\ch$. On the other hand, 
it is a consequence of
Shelah's theorem that no nontrivial automorphism can have a simple (e.g.~Borel)
definition. 

We do not know if for example Shelah's result has been extended to say Boolean
algebra homomorphisms (or even monomorphisms or epimorphisms). 
Putting this in the framework of this paper, we obtain an even more general question.

\begin{question}
\label{q-1}
What are the order embeddings of $\pnfin$ into itself with convex range\tu? 
Is it consistent that they are all trivial\tu?
\end{question}

\noindent In question~\ref{q-1} the precise meaning of \emph{trivial} is a mapping 
$[a]\mapsto\bigl[h[a]\cup[b]\bigr]$ where $h:A\to \N$ is an injection with $A$
cofinite and $b\cap \ran(h)$ finite.  


\subsubsection{$(\irrfin,\lefnt)$}
\label{sec:irrfin-lefnt}

The \emph{eventual dominance} quasi ordering $\lefnt$ of the irrationals $\irr$ is
given by $x\lefnt y$ if $x(n)\le y(n)$ for all but finitely many $n\in\N$. The
quotient $(\irrfin,\lefnt)$ over the equivalence relation $\eqfnt$ given by $x\eqfnt
y$ if $x\lefnt y$ and $y\lefnt x$ (i.e.~$x(n)=y(n)$ for all but finitely many $n$)
is a lattice. Analogously to $\pnfin$ we can formulate the notion of a ``trivial''
order embedding.

\begin{defn}
An order endomorphism from $\irrfin$ into itself is called \emph{trivial} if it is of
the form
\begin{equation}
  \label{eq:109}
  \sigma([x])=[\pi_g(x)\bigext\zero_{{\N}\setminus{\dom(g)}}+y]
\end{equation}
for some finite--one $g:A\to \N$, and some $y\in\irr$. 
\end{defn}

\begin{prop}
Every such trivial mapping is indeed an order preserving map between
$(\irrfin,\lefnt)$ and itself. 
\end{prop}

\begin{prop}
A trivial endomorphism as in equation~\tu{\eqref{eq:109}} is an embedding 
iff $\ran(g)$ is cofinite.
\end{prop}

\begin{prop}
A trivial endomorphism has convex range iff it can be represented by 
equation~\tu{\eqref{eq:109}} with an injection $g$.
\end{prop}

\begin{prop}
A trivial endomorphism is an epimorphism iff it can be represented by 
equation~\tu{\eqref{eq:109}} with $A$ cofinite and $g$ an injection.
\end{prop}

\noindent Therefore, every trivial embedding of $\irrfin$ into itself with convex range 
is of the form~$[x]\mapsto[\pi_g(x)\bigext\zero_{{\N}\setminus{\dom(g)}}+y]$ 
where $g:A\to\N$ is an injection with cofinite range and $y\in\irr$. 
And it follows that every trivial automorphism is of the form 
$[x]\mapsto[\pi_g(x)\bigext\zero_{{\N}\setminus{\dom(g)}}]$ where $g:A\to B$ is a
bijection and $A$ and $B$ are both cofinite subsets of $\N$. 

We could not find any result in the literature for $\irrfin$ corresponding to
Shelah's theorem, nor did we attempt to construct a nontrivial endomorphism. 

\begin{question}
\label{q-2}
Is the existence of a nontrivial automorphism of $(\irrfin,\lefnt)$ consistent\tu?
If so, is it consistent that every automorphism of $\irrfin$ is trivial\tu?
\end{question}

Putting question~\ref{q-2} into the framework of embeddings with convex range, 
we obtain the following generalization.

\begin{question}
What are the order embeddings of $\irrfin$ into itself with convex range\tu? Is it
consistent that they are all trivial\tu?
\end{question}

\vspace{-30pt}
\bibliographystyle{amsalpha}
\bibliography{database}

\end{document}

%% file: macros.tex
\newcommand{\ch}{\mathrm{CH}}

\newcommand{\dom}{\operatorname{dom}}
\newcommand{\ran}{\operatorname{ran}}

\newcommand{\interior}{^\circ}
\newcommand{\ro}{\operatorname{RO}}

\newcommand{\power}{\P}
\newcommand{\lbrak}{\bigl\|}
\newcommand{\rbrak}{\bigr\|}

\newcommand{\inv}{^{-1}}

\newcommand\bigext{^\frown}
\newcommand{\djun}{\amalg}

\renewcommand{\div}{\mathbin{/}}

\newcommand{\diff}{\bigtriangleup}

\newcommand{\borel}{\operatorname{Borel}}
\newcommand{\cat}{\operatorname{Cat}}
\newcommand{\cts}{C}

\newcommand{\pfn}{\dashrightarrow}

\newcommand{\pos}{\operatorname{pv}}
\newcommand{\atom}{\operatorname{At}}

\newcommand{\ifm}{\bigwedge}
\newcommand{\spr}{\bigvee}

\newcommand{\lorand}{{\lor}\hspace{-4pt}{\land}}

\newcommand{\iprod}[2]{\<#1,#2\>}

\newcommand{\incompat}{\not\approx}
\newcommand{\rel}{\mathrel{R}}

\newcommand{\@@eq}[2]{=_{#1}^{#2}}
\newcommand{\@@neq}[2]{\ne_{#1}^{#2}}
\newcommand{\@@l}[2]{<_{#1}^{#2}}
\newcommand{\@@nl}[2]{\nless_{#1}^{#2}}
\newcommand{\@@le}[2]{\le_{#1}^{#2}}
\newcommand{\@@nle}[2]{\nleq_{#1}^{#2}}
\newcommand{\@@equiv}[2]{\equiv_{#1}^{#2}}

\newcommand{\@in}[1]{\in_#1}
\newcommand{\@eq}[1]{=_#1}
\newcommand{\@le}[1]{\le_#1}
\newcommand{\@nle}[1]{\nleq_#1}
\newcommand{\@ge}[1]{\ge_#1}

\newcommand{\u@le}[1]{\le^#1}
\newcommand{\u@l}[1]{<^#1}
\newcommand{\u@eq}[1]{=^#1}

\newcommand{\@lefnt}[1]{\@@le#1*}
\newcommand{\@nlefnt}[1]{\@@nle#1*}
\newcommand{\@eqfnt}[1]{\@@eq#1*}

\newcommand{\subseteqfnt}{\subseteq^*}

\newcommand{\lefnt}{\mathrel{\le^*}}

\newcommand{\eqfnt}{\mathrel{=^*}}

\newcommand{\len}[1]{\le_{#1}}
\newcommand{\gen}[1]{\@ge{{#1}}}
\newcommand{\altle}{\lesssim}
\newcommand{\altl}{\lnsim}

\newcommand{\str}{\mathrm{str}}

\newcommand{\asym}{\mathrm{asym}}
\newcommand{\trn}{\mathrm{tran}}
\newcommand{\negstr}{-\str}

\newcommand{\strlen}[1]{\@@le{#1}\str}
\newcommand{\nstrlen}[1]{\@@nle{#1}\str}
\newcommand{\strln}[1]{\@@l{#1}\str}
\newcommand{\nstrln}[1]{\@@nl{#1}\str}
\newcommand{\strequiv}[1]{\@@equiv{#1}\str}

\newcommand{\lenn}[2]{\@@le{#1}{#2}}
\newcommand{\nlenn}[2]{\@@nle{#1}{#2}}
\newcommand{\lnn}[2]{\@@l{#1}{#2}}
\newcommand{\nlnn}[2]{\@@nl{#1}{#2}}

\newcommand{\eqnn}[2]{\@@eq{#1}{#2}}
\newcommand{\neqnn}[2]{\@@neq{#1}{#2}}
\newcommand{\asymle}{\u@le\asym}
\newcommand{\asyml}{\u@l\asym}

\newcommand{\strtrnlen}[1]{\@@le{#1}{\str,\trn}}
\newcommand{\strtrnln}[1]{\@@l{#1}{\str,\trn}}
\newcommand{\trnle}{\u@le{\trn}}
\newcommand{\trnl}{\u@l{\trn}}
\newcommand{\trneq}{\u@eq{\trn}}
\newcommand{\negstrlen}[1]{\@@le{#1}\negstr}
\newcommand{\nnegstrlen}[1]{\@@nle{#1}\negstr}
\newcommand{\nnegstrln}[1]{\@@nl{#1}\negstr}
\newcommand{\negstreq}[1]{\@@eq{#1}\negstr}

\renewcommand{\aa}{\mathrm{aa}}
\newcommand{\eqaa}{\@eq\aa}
\newcommand{\inaa}{\@in\aa}
\newcommand{\laa}{<_\aa}
\newcommand{\leaa}{\@le\aa}
\newcommand{\nleaa}{\@nle\aa}
\newcommand{\geaa}{\@ge\aa}
\newcommand{\lefntaa}{\@lefnt\aa}
\newcommand{\nlefntaa}{\@nlefnt\aa}
\newcommand{\eqfntaa}{\@eqfnt\aa}
\newcommand{\subseteqaa}{\subseteq_{\aa}}

\newcommand{\naa}{\mathrm{naa}}
\newcommand{\eqnaa}{\@eq\naa}
\newcommand{\eqfntnaa}{\@eqfnt\naa}
\newcommand{\lenaa}{\@le\naa}
\newcommand{\lefntnaa}{\@lefnt\naa}

\newcommand{\eqae}{\@eq\ae}
\newcommand{\leae}{\@le\ae}

\newcommand{\forces}{\mskip 5mu plus 5mu\|\hspace{-2.5pt}{\textstyle \frac{\hspace{4.5pt}}{\hspace{4.5pt}}}\mskip 5mu plus 5mu}

\newcommand{\Iff}{\espc\mathrm{iff}\espc}
\newcommand{\If}{\espc\textrm{if}\espc}
\newcommand{\impls}{\espc\mathrm{implies}\espc}
\renewcommand{\And}{\espc\mathrm{and}\espc}

\renewcommand{\and}{\hespc\mathrm{and}\hespc}

\newcommand{\zero}{\mathbf{0}}
\newcommand{\one}{\mathbf{1}}
\newcommand{\Fin}{\mathrm{Fin}}
\newcommand{\pN}{\power(\N)}
\newcommand{\reals}{\mathbb R}

\newcommand{\N}{\mathbb N}

\newcommand{\rationals}{\mathbb Q}

\newcommand{\irrationals}{\mathbb N^{\hsp\mathbb N}}
\newcommand{\irri}[1]{\N^{\hsp#1}}

\newcommand{\nulls}{\mathcal N}
\newcommand{\meager}{\M}

\newcommand{\bp}{\mathrm{BP}}

\newcommand{\pnfin}{\pN\div\Fin}
\newcommand{\irrfin}{\irr\div\Fin}

\newcommand{\random}{\R}

\newcommand{\borelaa}[1]{\borel(#1)\div{\eqaa}}

\newcommand{\baire}{\operatorname{Baire}}

\newcommand{\baireaa}[1]{\baire(#1)\div{\eqaa}}

\DeclareFontFamily{U}{cmsy}{}
\DeclareFontShape{U}{cmsy}{m}{n}{<12> sfixed * [10] cmsy10 
<10> <9> <8> <7> <6> <5> sfixed * [10] cmsy10}{}
\DeclareSymbolFont{customtwo}{U}{cmsy}{m}{n} 
\DeclareMathSymbol{\sctn}{\mathord}{customtwo}{"78}

\DeclareFontFamily{U}{cmmi}{}
\DeclareFontShape{U}{cmmi}{m}{n}{<20> sfixed * [11] cmmib10 <12> sfixed * [10]
cmmi10 <10> <9> <8> sfixed * [6] cmmi6 <5> <6> <7> sfixed * [5] cmmi5}{}
\DeclareSymbolFont{custom}{U}{cmmi}{m}{n}
\DeclareMathSymbol{\rharpoon}{\mathord}{custom}{"2A}
\newlength{\widt}
\newlength{\widttwo}
\newlength{\hgt}

\newcommand{\<}{\langle}
\renewcommand{\>}{\rangle}

\newcommand{\irr}{\irrationals}

\newcommand{\intgr}{\mathbb{Z}}
\newcommand{\espc}{\quad}
\newlength{\@@hespc}
\newcommand{\hespc}{\hspace{\@@hespc}}
\newcommand{\overbar}[1]{\,\overline{\!{#1}}}
\newcommand{\overbarg}[1]{\hsp\overline{\hspb{#1}}}

\renewcommand{\ae}{\mathrm{ae}}

\newcommand{\ulc}{\ulcorner}
\newcommand{\urc}{\urcorner}

\newcommand{\tu}{\textup}

\newcommand{\F}{\mathcal F}
\newcommand{\G}{\mathcal G}

\newcommand{\M}{\mathcal M}

\renewcommand{\P}{\Pcal}
\newcommand{\Pcal}{\mathcal P}

\newcommand{\R}{\mathcal R}
\newcommand{\Section}{{\char'237}}
\renewcommand{\S}{\mathcal S}

\newcommand{\U}{\mathcal U}
\newcommand{\V}{\mathcal V}
\newcommand{\W}{\mathcal W}

\newcommand{\hsp}{\mspace{1.5mu}}
\newcommand{\hspb}{\mspace{-1.5mu}}
\newcommand{\spc}{\,\,\,}

\newcounter{saveenumi}




%% file: embeddings.bbl
\providecommand{\bysame}{\leavevmode\hbox to3em{\hrulefill}\thinspace}
\providecommand{\MR}{\relax\ifhmode\unskip\space\fi MR }
\providecommand{\MRhref}[2]{%
  \href{http://www.ams.org/mathscinet-getitem?mr=#1}{#2}
}
  \renewcommand{\MRhref}[2]{
    \href{http://www.ams.org/mathscinet-getitem?mr=#1}{#2}
  }
\renewcommand{\MR}[2]{\MRhref{#1}{#2} $\uparrow\hspb$}
\begin{thebibliography}{Weh92}

\bibitem[AJ94]{MR1365749} 
Samson Abramsky and Achim Jung, \emph{Domain theory}, Handbook of logic in
  computer science, Vol.\ 3, 
  Oxford Univ. Press, New York, 1994, pp.~1--168. \MR{MR1365749}{(97b:68121)}

\bibitem[Coh63]{MR0157890}
Paul~J. Cohen, \emph{The independence of the continuum hypothesis}, Proc. Nat.
  Acad. Sci. U.S.A. \textbf{50} (1963), 1143--1148. \MR{MR0157890}{(28 \#1118)}

\bibitem[Coh64]{MR0159745}
\bysame, \emph{The independence of the continuum hypothesis. {II}}, Proc. Nat.
  Acad. Sci. U.S.A. \textbf{51} (1964), 105--110. \MR{MR0159745}{(28 \#2962)}

\bibitem[DP02]{MR1902334}
B.~A. Davey and H.~A. Priestley, \emph{Introduction to lattices and order},
  second ed., Cambridge University Press, New York, 2002. \MR{MR1902334}{(2003e:06001)}

\bibitem[DR81]{MR620665}
Dwight Duffus and Ivan Rival, \emph{A structure theory for ordered sets},
  Discrete Math. \textbf{35} (1981), 53--118. \MR{MR620665}{(82f:06001)}

\bibitem[DW87]{MR942216}
H.~G. Dales and W.~H. Woodin, \emph{An introduction to independence for
  analysts}, London Mathematical Society Lecture Note Series, vol. 115,
  Cambridge University Press, Cambridge, 1987. \MR{MR942216}{(90d:03101)}

\bibitem[Eng97]{MR1429390}
Konrad Engel, \emph{Sperner theory}, Encyclopedia of Mathematics and its
  Applications, vol.~65, Cambridge University Press, Cambridge, 1997.
  \MR{MR1429390}{(98m:05187)}

\bibitem[Hir00a]{MR1621933}
James Hirschorn, \href{http://homepage.univie.ac.at/James.Hirschorn/research/tbf/tbf.html}{\emph{Towers of {B}orel functions}}, Proc. Amer. Math. Soc.
  \textbf{128} (2000), no.~2, 599--604. \MR{MR1621933}{(2000c:03040)}

\bibitem[Hir00b]{MR1784706}
\bysame, 
\href{http://homepage.univie.ac.at/James.Hirschorn/research/tmf/tmf.html}{\emph{Towers of measurable functions}}, Fund. Math. \textbf{164}
  (2000), no.~2, 165--192. \MR{MR1784706}{(2002i:03056)}

\bibitem[Hir06a]{irrationals}
\bysame, \emph{Characterizing the ordering of the irrationals by eventual
  dominance}, in progress, 2006. $\uparrow\hspb$

\bibitem[Hir06b]{pinning}
\bysame, \emph{Pinning quasi orders with their endomorphisms}, 
\href{http://arxiv.org/abs/math.RA/0612495}{arXiv:math.RA/0612495}, 2006.
\href{http://homepage.univie.ac.at/James.Hirschorn/research/pinning/pinning.html}{\texttt{http://homepage.univie.ac.at/James.Hirschorn/research/pinning/pinning.html}}
$\uparrow\hspb$

\bibitem[Kan69]{MR0285644}
Pl. Kannappan, \emph{Characterizing topology on an {A}belian semigroup by a
  functional equation}, Portugal. Math. \textbf{28} (1969), 97--101.
  \MR{MR0285644}{(44 \#2862)}

\bibitem[Kap49]{MR0031193}
Irving Kaplansky, \emph{Normed algebras}, Duke Math. J. \textbf{16} (1949),
  399--418. \MR{MR0031193}{(11,115d)}

\bibitem[Kec95]{MR1321597}
Alexander~S. Kechris, \emph{Classical descriptive set theory}, Graduate Texts
  in Mathematics, vol. 156, Springer-Verlag, New York, 1995. \MR{MR1321597
 }{(96e:03057)}

\bibitem[Kop89]{MR991565}
Sabine Koppelberg, \emph{Handbook of {B}oolean algebras. {V}ol. 1},
  North-Holland Publishing Co., Amsterdam, 1989, Edited by J. Donald Monk and
  Robert Bonnet. \MR{MR991565}{(90k:06002)}

\bibitem[Kun80]{MR597342}
Kenneth Kunen, \emph{Set theory}, Studies in Logic and the Foundations of
  Mathematics, vol. 102, North-Holland Publishing Co., Amsterdam, 1980, An
  introduction to independence proofs. \MR{MR597342}{(82f:03001)}

\bibitem[Sco67]{MR0218233}
Dana Scott, \emph{A proof of the independence of the continuum hypothesis},
  Math. Systems Theory \textbf{1} (1967), 89--111. \MR{MR0218233}{(36 \#1321)}

\bibitem[She82]{MR675955}
Saharon Shelah, \emph{Proper forcing}, Lecture Notes in Mathematics, vol. 940,
  Springer-Verlag, Berlin, 1982. \MR{MR675955}{(84h:03002)}

\bibitem[She98]{MR1623206}
\bysame, \emph{Proper and improper forcing}, second ed., Perspectives in
  Mathematical Logic, Springer-Verlag, Berlin, 1998. \MR{MR1623206}{(98m:03002)}

\bibitem[Smy92]{MR1426367}
M.~B. Smyth, \emph{Topology}, Handbook of logic in computer science, Vol.\ 1,
  Oxford Univ. Press, New York, 1992,
  pp.~641--761. \MR{MR1426367}{(97i:68131)}

\bibitem[Sol70]{MR0265151}
Robert~M. Solovay, \emph{A model of set-theory in which every set of reals is
  {L}ebesgue measurable}, Ann. of Math. (2) \textbf{92} (1970), 1--56.
  \MR{MR0265151}{(42 \#64)}

\bibitem[Weh92]{MR1190444}
Friedrich Wehrung, \emph{Injective positively ordered monoids. {I}, {II}}, J.
  Pure Appl. Algebra \textbf{83} (1992), no.~1, 43--82, 83--100. \MR{MR1190444}{(93k:06023)}

\end{thebibliography}
